\RequirePackage{fix-cm}
\documentclass[smallextended]{svjour3}       
\smartqed  
\usepackage{graphicx}
\usepackage{amsfonts}
\usepackage{comment}

\newtheorem{algorithm}{Algorithm}

\usepackage{mathptmx}      

\begin{document}

\title{The regularizing Levenberg-Marquardt scheme for history matching of petroleum reservoirs 
}


\author{Marco A. Iglesias    \and
     Clint Dawson.
}


\institute{Marco A. Iglesias \at
              University of Warwick \\
             \email{M.A.Iglesias-Hernandez@warwick.ac.uk
           \and
Clint Dawson \at
              The University of Texas at Austin \\
             \email{clint@ices.utexas.edu}           
}   
\and 
}

\date{Received: date / Accepted: date}

\maketitle

\begin{abstract}
In this paper we study a history matching approach that consists of finding stable approximations to the problem of minimizing the weighted least-squares functional that penalizes the misfit between the reservoir model predictions $G(u)$ and noisy observations $y^{\eta}$. In other words, we are interested in computing $u^{\eta}\equiv \arg\min_{u\in X}\frac{1}{2}\vert\vert \Gamma^{-1/2}(y-G(u))\vert\vert_{Y}^{2} $ where $\Gamma$ is the measurements error covariance, $Y$ is the observation space and $X$ is a set of admissible parameters. This is an ill-posed nonlinear inverse problem that we address by means of the regularizing Levenberg-Marquardt scheme developed in \cite{Hanke,Hanke2}. Under certain conditions on $G$, the theory of \cite{Hanke,Hanke2} ensures convergence of the scheme to stable approximations to the inverse problem. We propose an implementation of the regularizing Levenberg-Marquardt scheme that enforces prior knowledge of the geologic properties. In particular, the prior mean $\overline{u}$ is incorporated in the initial guess of the algorithm and the prior error covariance $C$ is enforced through the definition of the parameter space $X$. Our main goal is to numerically show that the proposed implementation of the regularizing Levenberg-Marquardt scheme of Hanke is a robust method capable of providing accurate estimates of the geologic properties for small noise measurements. In addition, we provide numerical evidence of the convergence and regularizing results predicted by the theory of \cite{Hanke,Hanke2} for a prototypical oil-water reservoir model. The performance for recovering the true permeability with the regularizing Levenberg-Marquardt scheme is compared against the more standard techniques for history matching proposed in \cite{Li,Tavakoli,svdRML,Oliver}. Our numerical experiments suggest that the history matching approach based on iterative regularization is robust and could potentially be used to improve further on various methodologies already proposed as effective tools for history matching in petroleum reservoirs
 \end{abstract}

\section{Introduction}\label{intro}

History matching is the process of modifying parameters (inputs) of a reservoir model so that the model response (output) matches production data. The adjusted parameters during the history matching process are often geologic properties of the subsurface whose lack of information gives rise to uncertainty in the predictions of the reservoir model. The parameters obtained by means of history matching are aimed to provide better predictions of the reservoir performance; these can potentially be used for optimal reservoir management and monitoring of the reservoir. Given the potential impact of the history matching process in the optimal production of hydrocarbons, numerous techniques have been proposed and widely investigated in the last decades. For a recent review of history matching techniques we refer the reader to \cite{OliverReview} . Standard approaches are presented in detail in the monograph \cite{Oliver}. 

Let us denote by $u$ the unknown parameter (geologic properties) in a reservoir whose dynamics are described with a parameter-to-output operator $G:X\to Y$ that maps the space of admissible parameters $X$ to the observation space $Y$. In this paper we study the history matching problem posed by the minimization of the weighted data misfit defined by 
\begin{eqnarray}\label{eq:1.1}
\Phi(u)\equiv \frac{1}{2}\vert\vert \Gamma^{-1/2}(y^{\eta}-G(u))\vert\vert_{Y}^{2} 
\end{eqnarray}
where $y^{\eta}$ is the production data and $\Gamma$ is the measurements error covariance. For reservoir modeling applications, the evaluation of the forward operator $G(u)$ involves the solution of highly nonlinear system of PDEs whose differential operators have spatially varying coefficients related to the geologic parameter $u$. Therefore, the operator $G$ is typically compact which from standard theory \cite{Engl} implies that the minimization of (\ref{eq:1.1}) is ill-posed in the sense of stability. In other words, a small perturbation of $y^{\eta}$ may correspond to large deviations from the corresponding solutions to the minimizer of (\ref{eq:1.1}) \cite{Engl,kravaris_identification_1985} . The aforementioned lack of stability can lead to the divergence of standard optimization schemes used to solve the least-squares problem (\ref{eq:1.1}) \cite{Nagy,Engl_bio}. In this paper we propose a computational approach for the solution of (\ref{eq:1.1}) by means of the regularizing Levenberg-Marquardt (LM) scheme developed by Hanke in \cite{Hanke,Hanke2}. The regularizing LM scheme belongs to the class of so-called iterative regularization techniques designed to compute stable approximation to inverse ill-posed problems like the one posed by the minimization of (\ref{eq:1.1}). In contrast to the other approaches where the problem is first regularized (e.g. by Tikhonov's method) and then optimized with a standard solver, with iterative regularization techniques, the aim is to regularize the problem within an algorithm that also provides an approximation to a minimizer of (\ref{eq:1.1}). In other words, iterative regularization schemes provide stable estimates that converge to a minimizer of $\Phi$ in the limit of small noise. A review of the analysis and applications of iterative regularization techniques can be found in \cite{Iterative}. For the regularizing LM scheme, the mathematical analysis of the convergence and regularizing properties is developed in \cite{Hanke,Hanke2}. This mathematical framework of the LM scheme motivates the implementation that we propose for history matching in petroleum reservoirs. The main objective of this paper is to numerically show that the proposed implementation of the regularizing LM scheme is a robust methodology for generating accurate estimates of geologic parameters given production data with small noise. Furthermore, we provide numerical comparisons of the performance of the proposed implementation with respect to one of the most standard approaches for deterministic history matching. While the regularizing LM scheme aims at solving a history matching problem posed differently from standard methods, there exist some technical similarities between our implementation and the standard approaches. We exploit those similarities to provide guidelines for a straightforward implementation of the proposed algorithm given routines and codes from standard optimization methods

The paper is organized as follows. Relevant literature is discussed in Section \ref{lr}. In Section \ref{ir} we introduce the application of the regularizing LM scheme of \cite{Hanke,Hanke2} to the history matching problem. The fundamental theoretical aspects of the regularizing LM scheme are discussed in subsection \ref{lm}. Computational aspects relevant to the implementation of the regularizing LM scheme are presented in subsection \ref{ca}. In subsection \ref{sa} we discuss one of the standard approaches for history matching based on the method used in \cite{Li,Tavakoli,svdRML,Oliver}. In subsection \ref{sim} we show fundamental similarities between the proposed implementation of the regularizing LM scheme and the aforementioned standard approach. In Section \ref{Numerics} we display numerical experiments to show the capabilities of the regularizing LM scheme for generating stable approximations to the proposed history matching approach. Our implementation of the regularizing LM scheme is applied to an incompressible oil-water reservoir model described in the Appendix. The regularizing properties of the LM scheme with respect to the noise level are studied in subsection \ref{noisesec}. The performance of the LM scheme with respect to relevant tunable parameters is illustrated in subsection \ref{taurho}. The efficacy of the regularizing LM scheme for different choices of variance in the prior covariance operator is investigated in  Section \ref{eq:ir1}. Finally, the proposed regularizing LM scheme is compared against the standard method of  \cite{Li,Tavakoli,svdRML,Oliver}. Conclusions and final remarks are provided in Section \ref{Conclu}.

\section{Literature Review}\label{lr}

For history matching applications posed in terms of (\ref{eq:1.1}), regularization has been typically addressed by reparameterizing the geologic properties with a small number of parameters (see \cite{OliverReview}, section 3.2 and references therein). A parameterization in terms of a finite dimensional (hence compact) set, ensures the well-posedness of (\ref{eq:1.1})  \cite{Isakov}. However, for history matching applications, only problems parameterized with very few (of the order of 10) parameters have been treated by minimizing a functional like (\ref{eq:1.1}). For reservoirs with highly heterogeneous geologic properties, thousands or even millions of parameters are required to fully resolve relevant geologic features. For those reservoirs, parameterizing the geologic properties with a small number of parameters may not be possible. In that case, minimizing (\ref{eq:1.1}) with standard optimization techniques may diverge due to the lack of stability described above. One of the most standard approaches for history matching that addresses the ill-posedness of the inverse problem is to minimize \cite{OliverReview}
\begin{eqnarray}\label{eq:sa}
J(u)\equiv  \frac{1}{2}\vert\vert \Gamma^{-1/2}(y^{\eta}-G(u))\vert\vert_{Y}^{2} +\frac{1}{2}\vert\vert C^{-1/2}(u-\overline{u})\vert\vert_{X}^2 
\end{eqnarray}
which can be thought as Tikhonov regularization of (\ref{eq:1.1}) \cite{Oliver}. Indeed, the second term in the right hand side of (\ref{eq:sa})
 \begin{eqnarray}\label{eq:1.3}
R(u)\equiv \frac{1}{2}\vert\vert C^{-1/2}(u-\overline{u})\vert\vert_{X}^2.
\end{eqnarray}
is a regularization term that alleviates the ill-posedness of the inverse problem (\ref{eq:1.1}).  Under some assumptions on $G$, the theory of Tikhonov regularization for nonlinear problems ensures that (\ref{eq:sa}) has a solution which is continuous with respect to $y^\eta$ \cite[Theorem 10.2]{Engl}. Therefore, the problem is well-posed and any standard optimization technique can be implemented for solving (\ref{eq:sa}). However, as the size of $R(u)$ decreases, the regularization properties of (\ref{eq:sa}) rely on the proper size of $R(u)$ relative to $\Phi(u)$ \cite[Theorem 10.3]{Engl}. Intuitively, if the regularization term $R(u)$ is ``too small'' (\textit{underestimated}) compare to $\Phi(u)$, the regularization provided by $R(u)$ may not suffice, resulting in the lack of \textit{stability}. On the other hand, if $R(u)$ is ``too large'' (\textit{overestimated}), the minimizer (\ref{eq:sa}) may produce estimates that lack \textit{fidelity} due to a potential poor data match. For a given fixed error covariance matrix $\Gamma$, then the relative size of $R(u)$ with respect to $\Phi(u)$ is determined by the covariance operator $C$. For reservoir applications the standard practice is to select $C$ based on geologic data which has no connection to the stability and fidelity issues of the inverse problem described above. Therefore, for some choice of $C$, the potential risk of instabilities and lack of fidelity in minimizing (\ref{eq:sa}) may arise. History matching applications where the minimization of (\ref{eq:sa}) with the Gauss-Newton led to instabilities have been reported in \cite{Li,Oliver}. In order to alleviate these instabilities, a standard Levenberg-Marquardt method has been proposed in \cite{Li,Tavakoli,svdRML,Oliver}. The well-posedness of the minimization of (\ref{eq:sa}) is a fundamental assumption for the application of those standard techniques. However, as we mentioned above, for some choices based on geological information of $C$ may result in an insufficient regularization of the term $R(u)$ which in turn, may lead to numerical instabilities in more general settings.

By minimizing (\ref{eq:1.1}) instead of (\ref{eq:sa}), the implementation of the regularizing LM scheme that we propose in this paper avoids the potential lack of fidelity of the standard approach previously discussed. While in the standard approach of minimizing (\ref{eq:sa}) the term $R(u)$ enforces the prior geological knowledge, in our implementation of the regularizing LM scheme, the prior mean is incorporated in the initial guess of the iterative algorithm and the geological constraint imposed by the prior covariance is enforced in the definition of the parameter space $X$. It is fundamental to emphasize that the regularizing LM scheme applied to the minimization of (\ref{eq:1.1}) is a regularization technique that aims at producing stable computational approximations to a minimizer of (\ref{eq:1.1}). Therefore, in contrast to the standard approach where the minimization of (\ref{eq:sa}) is assumed well-posed and so standard optimization techniques can be applied, the regularizing LM approach of Hanke postulates an algorithm that alleviates the ill-posedness in the minimization of (\ref{eq:1.1}) while computing an approximation that converges to a minimizer of (\ref{eq:1.1}) for small observational noise \cite{Hanke,Hanke2}.

Iterative regularization techniques such as the regularizing LM scheme, have been successfully used for the solution of a wide class of inverse problems in several disciplines. In particular, for the inversion of data in subsurface flow models, in \cite{IglesiasDawson2} the authors study a simplified version of the regularizing LM scheme (see discussion of Section \ref{ir}) for the inversion of pressure in single-phase Darcy flow. In \cite{Iglesias5} a truncated regularizing Newton-Conjugate gradient is implemented to invert combined surface deformation and pressure data in a coupled flow-geomechanics problem. 
In this context of data inversion, the present work aims at extending the treatment in \cite{IglesiasDawson2} by using a two-phase (oil-water) reservoir model which, in contrast to the model studied in \cite{IglesiasDawson2}, is nonlinear with respect to the state variables (pressures and saturations). The increase in nonlinearity of the forward operator imposes severe challenges on the regularization properties of the technique under consideration. However, the present work offers the numerical evidence that the theory of Hanke may be applied to the forward operator that arises from the prototypical oil-water reservoir model that we consider in our numerical experiments. Further investigations of the regularizing LM scheme and other iterative regularization technique may lead to the development of efficient tools for history matching applications.

\section{Iterative regularization for history matching}\label{ir}

In this section we present the application of the regularizing LM scheme of \cite{Hanke,Hanke2} for history matching by means of minimizing $\Phi$ defined in (\ref{eq:1.1}) where $G$ is an arbitrary reservoir model that captures the flow dynamics perfectly. For our experiments of Section \ref{Numerics}, we use a forward operator $G$ that we derived from a prototypical incompressible oil-water model (see Appendix). As we indicated in Section \ref{intro}, due to the ill-posedness of the minimization of (\ref{eq:1.1}), a regularization algorithm is required to compute stable solutions to the inverse problem. While a broad spectrum of iterative regularization techniques can be used, here we consider the regularizing LM technique because of the computational similarities with standard LM methods that are typically used for standard approaches in history matching that we described in Section \ref{sa}.
 
Assume we are provided measurements possibly corrupted by noise $y^{\eta}$ with the noise level denoted by $\eta$ and defined by 
\begin{eqnarray}\label{eq:noise}
\vert\vert \Gamma^{-1/2}(y^{\eta}- G(u^{\dagger}))\vert\vert_{Y}\leq \eta
\end{eqnarray}
where $u^{\dagger}$ denotes the true geologic properties of the reservoir. Note that if we knew the truth $u^{\dagger}$, in the absence of observational error (i.e. $\eta\to 0$), we would measure the model predictions of the truth  $y\equiv G(u^{\dagger})$. Therefore, $\eta$ in (\ref{eq:noise}) is an upper bound for the observational noise. In practice, $\eta$ can be defined from the measurement information also used for defining the measurement error covariance $\Gamma$.

\subsection{The regularizing Levenberg-Marquardt method}\label{lm}

The aim of the regularizing LM scheme is to compute stable approximations to a minimizer of (\ref{eq:1.1}). In other words, we want to compute $u_{\eta}$ such that $u_{\eta}\to u$ as $\eta\to 0$, where $u$ is a minimum of (\ref{eq:1.1}) in the limit $\eta\to 0$. The approximation $u^{\eta}$ is the limit of a finite sequence of estimates $\{u_{m}^{\eta}\}_{m=1}^{N}$ computed as we now describe. Given, the estimate $u_{m}^{\eta}$ at the $m$th iteration of the scheme, the aim is to construct an update $u_{m+1}^{\eta}=u_{m}^{\eta}+\Delta u_{m}^{\eta}$, where the increment $\Delta u _{m}^{\eta}$ is obtained by solving 
\begin{eqnarray}\label{eq:linear}
y^{\eta}-G(u_{m}^{\eta})=DG(u_{m}^{\eta})\Delta u_{m}^{\eta},
\end{eqnarray}
where $DG(u_m)$ is the Frechet derivative of $G$ at $u_{m}$. Note that (\ref{eq:linear}) is a linearized version of the equation $y^{\eta}=G(u)$ satisfied by a minimizer $u$ in the case that the minimum in (\ref{eq:1.1}) corresponds to $\Phi(u)=0$. Note that the truth increment defined by $\Delta u_{m}^{\dagger}=u^{\dagger}-u_{m}^{\eta}$ satisfies
\begin{eqnarray}\label{eq:linearV2}
G(u^{\dagger})-G(u_{m}^{\eta})-R(u_{m}^{\eta},u^{\dagger})=DG(u_{m}^{\eta})\Delta u_{m}^{\dagger},
\end{eqnarray}
where $R(u_{m}^{\eta},u^{\dagger})$ is the Taylor remainder of $G$ at $u^{\dagger}$ around $u_{m}^{\eta}$. As we mentioned in Section \ref{intro}, the ill-posedness in the minimization of (\ref{eq:1.1}) can be attributed to the compactness of the forward operator which is often encountered in PDE-constrained inverse problems \cite{Isakov}.  From the compactness of the Frechet derivative of a compact operator \cite{Colton}, it follows that the linear operator $DG(u_{m}^{\eta})$ is compact for each $m\in \mathbb{N}$. Therefore, the linear inverse problem (\ref{eq:linear}) also requires regularization. In the regularizing LM scheme of Hanke \cite{Hanke,Hanke2}, Tikonov regularization is applied to (\ref{eq:linear}) by computing 
\begin{eqnarray}\label{eq:4.0}
\Delta u_{m}^{\eta}(\alpha)=\textrm{argmin}_{w\in X}J_{LM}^{m}(w,\alpha)
\end{eqnarray}
 where 
\begin{eqnarray}\label{eq:4.1}
J_{LM}^{m}(w,\alpha)\equiv \frac{1}{2}\vert\vert \Gamma^{-1/2}(y^{\eta}-G(u_{m}^{\eta})-DG(u_{m}^{\eta})
w)\vert\vert_{Y}^{2}
+\frac{1}{2}\alpha\vert\vert C^{-1/2}w\vert\vert_{X}^2
\end{eqnarray}
The choice of $\alpha$ is fundamental to ensure the proper regularization of the inverse problem. Hanke proposes $\alpha$ such that
\begin{eqnarray}\label{eq:4.2}
\vert\vert \Gamma^{-1/2} (y^{\eta}-G(u_{m}^{\eta})-DG(u_{m}^{\eta})\Delta u_{m}^{\eta}(\alpha))\vert\vert_{Y}^{2}
\ge \rho^{2} \vert\vert \Gamma^{-1/2} (y^{\eta}-G(u_{m}^{\eta}))\vert\vert_{Y}^{2}
\end{eqnarray}
for some $\rho\in (0,1)$ (note that $\Delta u_{m}^{\eta}(\alpha)$ in (\ref{eq:4.2}) depends on $\alpha$ via (\ref{eq:4.0})-(\ref{eq:4.1})).

To gain further insight of the regularizing LM scheme as well as the selection of $\alpha$, let us define 
\begin{eqnarray}\label{eq:3.4}
d^{\eta,m}\equiv y^{\eta}-G(u_{m}^{\eta}),\qquad & d^{m}\equiv G(u^{\dagger})-G(u_{m}^{\eta})-R(u_{m}^{\eta},u^{\dagger}),\nonumber\\
DG(u_{m}^{\eta})\equiv g^{m},& \overline{u}^{m}\equiv u_{m}-\overline{u}.
\end{eqnarray} 
which applied to (\ref{eq:linear}) and (\ref{eq:linearV2}) yields
\begin{eqnarray}\label{eq:linearV3}
d^{\eta,m}=g^{m}\Delta u_{m}^{\eta},\qquad  d^{m}=g^{m}\Delta u_{m}^{\dagger}
\end{eqnarray}
From definitions (\ref{eq:3.4}), expressions (\ref{eq:4.1}) and (\ref{eq:4.2}) become
\begin{eqnarray}\label{eq:4.3}
J_{LM}^{m}(w,\alpha)\equiv \frac{1}{2}\vert\vert \Gamma^{-1/2}(d^{\eta,m}-g^{m}
\Delta u_{m}^{\eta}(\alpha))\vert\vert_{Y}^{2}
+\frac{1}{2}\alpha\vert\vert C^{-1/2}\Delta u_{m}^{\eta}(\alpha)\vert\vert_{X}^2
\end{eqnarray}
and
\begin{eqnarray}\label{eq:4.4}
\vert\vert \Gamma^{-1/2} (d^{\eta,m}-g^{m}\Delta u_{m}^{\eta}(\alpha))\vert\vert_{Y}^{2}
\ge \rho^{2} \vert\vert \Gamma^{-1/2} d^{\eta,m} \vert\vert_{Y}^{2}
\end{eqnarray}
respectively. Therefore, each iteration of the proposed scheme can be viewed as a Tikhonov regularization for the linear inverse problem of find $\Delta u_{m}^{\eta}$ given data $d^{\eta,m}$, where the latter is a noisy version of $d^{m}$. Note that, from (\ref{eq:3.4}) it follows that
\begin{eqnarray}\label{eq:linearV5}
\vert\vert   \Gamma^{-1/2}( y^{\eta}-y-R(u_{m}^{\eta},u^{\dagger}))\vert\vert_{Y}=\vert\vert \Gamma^{-1/2}(d^{m}-d^{\eta,m})\vert\vert_{Y}
\end{eqnarray}
The regularizing LM scheme assumes that it is possible to find $\rho\in (0,1)$ such that
\begin{eqnarray}\label{eq:linearV6}
\vert\vert   \Gamma^{-1/2}( y^{\eta}-y-R(u_{m}^{\eta},u^{\dagger}))\vert\vert_{Y}=\vert\vert \Gamma^{-1/2}(d^{m}-d^{\eta,m})\vert\vert_{Y}\leq \rho \vert\vert \Gamma^{-1/2}d^{\eta,m}\vert\vert_{Y}.
\end{eqnarray}
The inequality in the previous expression implies that the size of the error in the data $d^{m}$ must be smaller than the size of the observations $d^{\eta,m}$. It is certainly hopeless to invert data whose error is of the order of the size of the observations. The $\rho$ in (\ref{eq:linearV6}) is used in expression (\ref{eq:4.2}) for choosing the regularization parameter $\alpha$. Moreover, from (\ref{eq:linearV6}) it is easy to see that the selection of $\alpha$ according to (\ref{eq:4.4}) implies
\begin{eqnarray}\label{eq:4.4V2}
\vert\vert \Gamma^{-1/2} (d^{\eta,m}-g^{m}\Delta u_{m}^{\eta}(\alpha))\vert\vert_{Y}^{2}
\ge  \vert\vert \Gamma^{-1/2}(d^{m}-d^{\eta,m})\vert\vert_{Y}
\end{eqnarray}
which is the \textit{discrepancy principle} applied to the inverse problem $d^{\eta,m}=g^{m}\Delta u_{m}^{\eta}$. The discrepancy principle states that the estimate $\Delta u_{m}^{\eta}(\alpha)$ of the solution to the inverse problem (\ref{eq:linearV3}) cannot produce an output $g^{m}\Delta u_{m}^{\eta}(\alpha)$ whose associated error is better than the noise level. For a discussion of the discrepancy principle in the context of linear inverse problems the reader is referred to \cite{Groetsch}. Let us now denote  by $\alpha_{m}$ a solution to inequality (\ref{eq:4.4}). Then, the update of the regularizing LM scheme is defined by
\begin{eqnarray}\label{eq:4.6}
u_{m+1}^{\eta}\equiv u_{m}^{\eta}+\Delta u_{m}^{\eta}(\alpha_{m})=u_{m}^{\eta}+\textrm{argmin}_{w\in X}J_{LM}^{m}(w,\alpha^{m})
\end{eqnarray}
which provides a new estimate of the geologic properties. The existence of $\alpha_{m}$ is proven in \cite{Groetsch,Iterative} (see also discussion below). The minimizer of (\ref{eq:4.1}) with $\alpha_{m}$ given by (\ref{eq:4.2}) provides a regularized solution to the linear inverse problem (\ref{eq:linear}). Furthermore, the regularizing LM scheme is terminated provided the $(k+1)$th iteration produces an estimate $u_{k+1}^{\eta}$ such that
\begin{eqnarray}\label{discrepancy}
\vert\vert \Gamma^{-1/2} (y^{\eta}-G(u_{k+1}^{\eta}))\vert\vert_{Y}\leq \tau\eta\leq \vert\vert \Gamma^{-1/2} (y^{\eta}-G(u_{k}^{\eta}))\vert\vert_{Y}
\end{eqnarray}
for $\tau>1/\rho$. The resulting estimate $u^{\eta}\equiv u_{k+1}^{\eta}$ is the desired stable approximation to the inverse problem of minimizing (\ref{eq:1.1}). Expression (\ref{discrepancy}) is also an application of the discrepancy principle which, in this context, states that the data misfit obtained with approximation to the inverse problem $u^{\eta}$ should not be smaller than the noise level $\eta$. Intuitively, if $\rho\approx 1$  (with $\rho<1$) then $\tau$ can be chosen $\tau \approx 1$ which, in turn, may result in estimates that provide a good data fit. The regularizing LM scheme is now summarized below: 

\begin{algorithm}[Regularizing Levenberg-Marquardt Scheme]\label{al:LM}
Consider the initial estimate $u_{0}^{\eta}=\overline{u}$. Choose parameters $\rho<1$ and $\tau >1/\rho$. For each $m=1,\dots,k$,
\begin{itemize}
\item[(1)]\textbf{Forward simulation}. Given $u_{m}^{\eta}$ simulate the model response $G(u_{m}^{\eta})$.
\item[(2)] \textbf{Stopping rule (Discrepancy Principle)}. If (\ref{discrepancy}) holds then stop (i.e. $m=k+1$). Output: $u_{m}^{\eta}$.
\item[(3)] \textbf{Update}. Define $u_{m+1}^{\eta}$ according to (\ref{eq:4.6}) with $J_{LM}^{m}$ defined in (\ref{eq:4.1}) and $\alpha^{m}$ is chosen according to the (\ref{eq:4.2}).
\end{itemize}
\end{algorithm}

\begin{remark}\label{remark2}
In the regularizing LM scheme, prior knowledge $\overline{u}$ of the unknown is incorporated as the initial guess of the LM algorithm. In addition, the prior covariance $C$ is included in the definition of the parameter space, which formally, can be defined as the completion of the original space $X$ under the norm $\vert\vert C^{-1/2}\cdot \vert\vert_{X}$. This choice of the space is reflected in the second term of the right hand side of (\ref{eq:4.1})
\end{remark}
  
The application of the discrepancy principle for the selection of $\alpha$ in (\ref{eq:4.2}) as well as the termination of the algorithm (\ref{discrepancy}) are key aspects for the regularization properties of the regularizing LM scheme.  In particular, we recall the following result proven in \cite[Theorem 2.3]{Hanke}.
\begin{theorem}[Hanke \cite{Hanke}]\label{Hanke_Theorem}
Let $\rho\in(0,1)$ and $\tau>1/\rho$. Assume that $DG$ is locally bounded and that $G$ satisfies
\begin{eqnarray}\label{eq:4.7}
\vert\vert G(u)-G(\tilde{u})-DG(u)(u-\tilde{u})\vert\vert_{Y}\leq C\vert\vert u-\tilde{u}\vert\vert_{X}\vert\vert G(u)-G(\tilde{u})\vert\vert_{Y}
\end{eqnarray}
locally in $X$. If $u_{0}$ is sufficiently close to a solution $u^{\star}$ of $G(u^{\dagger})=G(u^{\star})$, then, the discrepancy principle (\ref{discrepancy}) terminates the LM algorithm with parameters $\alpha$ from (\ref{eq:4.2}) after a finite number of iterations $k(\eta)$. Moreover, the corresponding approximations $u_{k(\eta)}^{\eta}$ converge to a solution of $G(u^{\dagger})=G(u)$ as $\eta\to 0$.
\end{theorem}

\begin{remark}\label{rem3}
From Theorem \ref{Hanke_Theorem} we see that as $\eta\to 0$, then the solution $u^{\eta}$ computed with the regularizing LM scheme converges to $u$ that satisfies $G(u^{\dagger})=G(u)$. Therefore, since from (\ref{eq:noise}) $y^{\eta}\to G(u^{\dagger})$ as $\eta \to 0$, it then follows that $u$ satisfies $\Phi(u)=0$ and so $u^{\eta}$ converges to a minimizer of $\Phi$ in the limit of $\eta \to 0$.
\end{remark}

In \cite{IglesiasDawson2} we implemented a particular case of Algorithm \ref{al:LM} for the estimation of absolute permeability with a single-phase (linear) reservoir model. Instead of choosing $\alpha^{n}$ as in (\ref{eq:4.2}), in the scheme of \cite{IglesiasDawson2} the Tikhonov parameter was chosen constant $\alpha=1$. This selection was sufficient to prove convergence in the same sense of Theorem \ref{Hanke_Theorem}. For the work reported in this paper, we initially implemented the algorithm of \cite{IglesiasDawson2} for the estimation of absolute permeability with the reservoir model of the Appendix. However, the need for an adaptive selection of $\alpha$ arose due to the highly nonlinear structure of the present forward model. While the rigorous application of Theorem \ref{Hanke_Theorem} for the forward operator $G$ of the Appendix remains an open problem, our numerical results give evidence that confirms the regularizing properties predicted by Hanke's theory.

\subsection{Computational Implementation of the regularizing LM scheme}\label{ca}

In this section we discuss computational aspects of the regularizing LM scheme. Our main goal is to provide a reproducible computationally efficient algorithm for history matching. We first notice that, for $\alpha$ fixed, the Euler-Lagrange equation associated to the minimization of (\ref{eq:4.3}) yields
\begin{eqnarray}\label{eq:4.8}
\Delta u_{m}^{\eta}(\alpha)= \Big[ DG^{\ast}(u_{m}^{\eta})\Gamma^{-1}DG(u_{m}^{\eta})+\alpha C^{-1}\Big]^{-1}DG^{\ast}(u_{m}^{\eta})\Gamma^{-1}(y^{\eta}-G(u_{m}^{\eta})).
\end{eqnarray}
where $DG^{\ast}(u_m^{\eta})$ is the adjoint operator of $DG(u_m^{\eta})$. Expression (\ref{eq:4.8}) involves the inversion of the operator  $DG^{\ast}(u_{m})\Gamma^{-1}DG(u_{m})+\alpha C^{-1}$ in the space $X$. However, for the reservoir application under consideration, the dimension of the parameter space $X$ is typically much larger than the dimension of the observation space $Y$. Therefore, for computational efficiency we consider the equivalence between (\ref{eq:4.8}) and
\begin{eqnarray}\label{eq:4.9}
\Delta u_{m}^{\eta}(\alpha)=C \, DG^{*}(u_{m}^{\eta}) \Big[ DG(u_{m}^{\eta}) \, C \, DG^{*}(u_{m}^{\eta})+\alpha\Gamma\Big]^{-1}(y^{\eta}-G(u_{m}^{\eta})).
\end{eqnarray}
which in finite dimensions can be shown from the matrix lemmas of \cite[Section 7.4]{Oliver}. In the infinite-dimensional case, the equivalence between (\ref{eq:4.8})-(\ref{eq:4.9}) is only formal. Note that, assuming that the sensitivities $DG(u_{m}^{\eta})$ and $DG^{*}(u_{m}^{\eta})$ are available, then either (\ref{eq:4.8}) or (\ref{eq:4.9}) can be easily computed for any given $\alpha$. It is therefore clear that the computation of $\alpha$ in (\ref{eq:4.2}) represents the main new aspect of the proposed implementation. However, the computation of $\alpha$ is fairly simple as we describe below.

Let us define
\begin{eqnarray}\label{eq:4.10}
\kappa_{m}^{\eta}(\alpha)\equiv  \vert\vert \Gamma^{-1/2}(y^{\eta}-G(u_{m}^{\eta})-DG(u_{m}^{\eta})\Delta_{m}^{\eta}u(\alpha))\vert\vert_{X}^{2}
\end{eqnarray}
We substitute expression (\ref{eq:4.9}) in (\ref{eq:4.10}) and from simple computations it follows that 
\begin{eqnarray}\label{eq:4.11}
\kappa_{m}^{\eta}(\alpha)=\alpha^2 \vert\vert  \Gamma^{1/2}[DG(u_{m}^{\eta}) \, C\, DG^{\ast}(u_{m}^{\eta})+\alpha\Gamma]^{-1}[y^{\eta}-G(u_{m}^{\eta})]\vert\vert_{Y}^2
\end{eqnarray}
From this expression we find that $\kappa_{m}^{\eta}(\alpha)$ is a continuous increasing function of $\alpha$. Moreover, it can be shown \cite[Chapter 4]{Iterative} that
\begin{eqnarray}\label{eq:4.13}
\kappa_{m}^{\eta}(\alpha)\in \Bigg[\frac{\rho^2}{\gamma} \vert\vert y^{\eta}-G(u_{m}^{\eta})\vert\vert^2, \vert\vert y^{\eta}-G(u_{m}^{\eta})\vert\vert^2\Bigg]
\end{eqnarray}
for all $\alpha\in [0,\infty)$ and for some $\gamma>1$. Moreover, the right end of the interval above is given by
\begin{eqnarray}\label{eq:4.12}
\lim_{\alpha\to \infty}\kappa_{m}^{\eta}(\alpha)=\vert\vert y^{\eta}-G(u_{m}^{\eta})\vert\vert^2
\end{eqnarray}
Since $\kappa_{m}^{\eta}(\alpha)$ is continuously increasing, it follows from (\ref{eq:4.13}) that there exists $\alpha^{\star}\in [0,\infty)$ such that
\begin{eqnarray}\label{eq:4.14}
\frac{\rho^2}{\gamma} \vert\vert y^{\eta}-G(u_{m}^{\eta})\vert\vert^2\leq \rho^{2}\vert\vert y^{\eta}-G(u_{m}^{\eta})\vert\vert^2= \kappa_{m}^{\eta}(\alpha^{\star})\leq  \vert\vert y^{\eta}-G(u_{m}^{\eta})\vert\vert^2
\end{eqnarray}
Note that any $\alpha_{m}$ such that $\alpha^{\star}\leq \alpha_m$ will therefore satisfy $\kappa_{m}^{\eta}(\alpha^{\star})\leq \kappa_{m}^{\eta}(\alpha_{m})$ which, from (\ref{eq:4.14}) implies (\ref{eq:4.2}) as required. Computationally, we can determine such $\alpha_{m}$ by constructing $\alpha_{m}^{j}\to \infty $ as $j\to \infty$. Let us consider, for example, $\alpha_{m}^{j+1}=2^{j+1}\alpha_{m}^{j}$ where $\alpha_{m}^{0}>0$ is an initial guess for $\alpha_{m}$. We claim that there exists $J<\infty$ such that
\begin{eqnarray}\label{eq:4.16}
\rho^{2}\vert\vert y^{\eta}-G(u_{m}^{\eta})\vert\vert^2\leq \kappa_{m}^{\eta}(\alpha_{m}^{J})
\end{eqnarray}
If no such $J$ exists then  
\begin{eqnarray}\label{eq:4.16}
\kappa_{m}^{\eta}(\alpha_{m}^{j})< \rho^{2}\vert\vert y^{\eta}-G(u_{m}^{\eta})\vert\vert^2
\end{eqnarray}
for all $j\in \mathbb{N}$. In particular, for sufficient large $j$, from (\ref{eq:4.12}) we find
\begin{eqnarray}\label{eq:4.16}
\vert\vert y^{\eta}-G(u_{m}^{\eta})\vert\vert^2< \rho^{2}\vert\vert y^{\eta}-G(u_{m}^{\eta})\vert\vert^2
\end{eqnarray}
which contradicts the hypothesis of $\rho<1$.  We define $\alpha_{m}\equiv \alpha_{m}^{J}$ and the update of the regularizing LM scheme
\begin{eqnarray}\label{eq:4.17}
u_{m+1}^{\eta}=u_{m}^{\eta} +C \, DG^{*}(u_{m}^{\eta})[DG(u_{m}^{\eta}) \, C \, DG^{\ast}(u_{m}^{\eta})+\alpha_{m} \Gamma]^{-1}[y^{\eta}-G(u_{m}^{\eta})]
\end{eqnarray}
Note that the computation of $\alpha_{m}$ requires the evaluation of $\kappa_{m}^{\eta}(\alpha)$ which from (\ref{eq:4.11}) involves the inversion of $[DG(u_{m}^{\eta}) \ C \ DG^{\ast}(u_{m}^{\eta})+\alpha_{m} \Gamma]^{-1}$. However, $DG(u_{m}^{\eta})\ C \ DG^{\ast}(u_{m}^{\eta})$ has to be assembled only once per iteration of the scheme (see the update equation (\ref{eq:4.17})). The cost of inverting $[DG(u_{m}^{\eta})\ C\ DG^{\ast}(u_{m}^{\eta})+\alpha_{m} \Gamma]^{-1}$ for different $\alpha_{m}$'s is negligible for the application under consideration due the small dimensionality of the observation space. Therefore, the cost of computing $\alpha_{m}$ that satisfies (\ref{eq:4.2}) is negligible compared to the cost of evaluating $G(u_{m}^{\eta})$ and assembling  $DG(u_{m}^{\eta})\ C \ DG^{\ast}(u_{m}^{\eta})$ which both, in turn, constitute the main computational cost per iteration of the proposed implementation of the regularizing LM scheme. With the aforementioned considerations we propose a computationally efficient implementation of the regularization LM scheme.
 
\begin{algorithm}[Regularizing LM Scheme (implementable version)]\label{al:LM2}
Let $u_{0}\in X$ be an initial guess. Choose parameters $\rho\in (0,1)$ and $\tau>1/\rho$. For each $n=1,\dots$,
\begin{itemize}
\item[(1)]\textbf{Solution to the forward model}. Given $u_{m}^{\eta}$ evaluate the forward operator $G(u_{m}^{\eta})$.
\item[(2)] \textbf{Stopping rule (Discrepancy Principle)}. If
\begin{eqnarray}\label{discrepancy2}
\vert\vert \Gamma^{1/2}(y^{\eta}-G(u_m^{\eta}))\vert\vert_{Y}\leq \tau\eta
\end{eqnarray}
stop. Output: $u_{m}^{\eta}$.
\item[(3)] Compute the sensitivity matrices $DG(u_{m}^{\eta})$, its adjoint operator  $DG(u_{m}^{\eta})^{\ast}$ and assemble matrix $DG(u_{m}^{\eta}) \, C \, DG^{\ast}(u_{m}^{\eta})$. Let $\alpha_{m}^{0}>0$ and $\alpha_{m}^{j+1}=2^{j}\alpha_{m}^{j+1}$. Let $J$ be such that 
\begin{eqnarray}\label{eq:4.18}
\rho^{2}\vert\vert y^{\eta}-G(u_{m}^{\eta})\vert\vert^2\leq \kappa_{m}^{\eta}(\alpha_{m}^{J})\nonumber\\
\equiv \alpha^2 \vert\vert  \Gamma^{1/2}[DG(u_{m}^{\eta}) \, C \, DG^{\ast}(u_{m}^{\eta})+\alpha_{m}^{J}\Gamma]^{-1}[y^{\eta}-G(u_{m}^{\eta})]\vert\vert_{Y}^2 
\end{eqnarray}
\textbf{Update}. Define 
\begin{eqnarray}\label{eq:4.19}
u_{m+1}^{\eta}=u_{m}^{\eta} +C \, DG^{*}(u_{m}^{\eta})[DG(u_{m}^{\eta}) \, C \, DG^{\ast}(u_{m}^{\eta})+\alpha_{m}^{J} \Gamma]^{-1}[y^{\eta}-G(u_{m}^{\eta})]
\end{eqnarray}
\end{itemize}
\end{algorithm}

\subsection{The standard approach for history matching}\label{sa}

As we indicated in Section \ref{intro}, one of the most standard approaches for history matching consist of minimizing (\ref{eq:sa}). For analogy with our proposed implementation for solving (\ref{eq:1.1}), we based the following discussion on the the application of the Levenberg-Marquardt algorithm used in \cite{Li,Tavakoli,svdRML,Oliver} for the minimization of (\ref{eq:sa}) in the standard approach. The aforementioned method consist of computing the sequence $u_{m+1}=u_{m}+\Delta u$ where the step $\Delta u$ satisfies
\begin{eqnarray}\label{eq:3.7}
\Big[ DG^{\ast}(u_{m})\Gamma^{-1}DG(u_{m})+C^{-1}+\lambda_{m} C^{-1}\Big]\Delta u= DG^{\ast}(u^{n})\Gamma^{-1}[y^{\eta}-G(u_{m}))-C^{-1}(u_{m}-\overline{u})]\nonumber\\
\end{eqnarray}
for some $\lambda_{m}>0$. The proposed update (\ref{eq:3.7}) is a scaled version of the standard LM algorithm for the solution of well-posed optimization problems \cite[Chapter 10]{Nocedal}.  For the history matching applications of \cite{Li,Tavakoli,svdRML,Oliver}, the suggested selection of $\lambda_{m}$ is the following. The initial $\lambda_{0}$ is chosen between $\sqrt{J(u_{0})/N_{d}}$ and $J(u_{0})/N_{d}$ where $N_{d}$ is the dimension of the observation space. For $m\ge 0$, $\lambda_{m+1}$ is chosen according to
\begin{eqnarray}\label{eq:3.9B}
\lambda_{m+1}=\left\{\begin{array}{cc}
\lambda_{m}/10 & \textrm{if}~~J(u_{m+1})< J(u_{m})\\
10\lambda_{m} & \textrm{if}~~J(u_{m+1})\ge J(u_{m})\end{array}\right.
\end{eqnarray}
In addition, the stopping criteria for the LM technique of \cite{Li,Tavakoli,svdRML,Oliver} is based on the following two stopping criteria
\begin{eqnarray}\label{eq:3.9C}
\frac{\vert J(u_{m+1})-J(u_{m})\vert }{J(u_{m+1})}\leq \epsilon_{0},\\ \frac{\vert\vert u_{m+1}-u_{m}\vert\vert_{X} }{ \vert\vert u_{m+1} \vert\vert_{X} }\leq \epsilon_{1}\label{eq:3.9D}
\end{eqnarray}
In order to understand the standard LM approach for minimizing (\ref{eq:sa}), note that (\ref{eq:3.7}) can be derived from the Euler-Lagrange equations for the minimization of 
\begin{eqnarray}\label{eq:3.8}
Q^{m}(\Delta u)\equiv  \frac{1}{2}\vert\vert \Gamma^{-1/2}(y^{\eta}-G(u_{m})-DG(u_{m})\Delta u)\vert\vert_{Y}^{2} \nonumber\\+\frac{1}{2}\vert\vert C^{-1/2}(\Delta u-(u_{m}-\overline{u})) \vert\vert_{X}^2+\frac{1}{2}\lambda_{m}\vert\vert C^{-1/2}\Delta u\vert\vert_{X}^2
\end{eqnarray}
In other words, $\Delta u_{m}=\arg\min_{v\in X} Q^{m}(v)$, and so each iteration step of the LM method of \cite{Li,Tavakoli,svdRML,Oliver} is the solution of a least-squares Tikhonov-type problem on the linearized inverse problem. Note that the choice $\lambda_{m}=0$ in (\ref{eq:3.8}) suppresses the extra regularization term in the right hand side of (\ref{eq:3.8}). Indeed, the initial motivation of the LM scheme used in \cite{Li,Tavakoli,svdRML,Oliver} was to alleviate the lack of stability of the Gauss-Newton (GN) method of \cite{Li} and \cite[section 8.4.2]{Oliver} which corresponds to $\lambda_{m}=0$ in (\ref{eq:3.7}) (\ref{eq:3.8}). The LM scheme of \cite{Li,Tavakoli,svdRML,Oliver} for the minimization of (\ref{eq:sa}) is an efficient strategy provided that the minimization of $J$ is a well-posed problem. However, as we indicated before, the regularization term $R(u)$ in (\ref{eq:sa}) may be insufficient for some choices of $C$. For some choices of $C$, in the following section we present numerical experiments demonstrating that the selection of $\lambda$ in (\ref{eq:3.9B}) and the stopping criteria of (\ref{eq:3.9C})-(\ref{eq:3.9D}) may lead to both lack of stability and fidelity in the computation of estimates of geologic parameters.

\subsection{Computational similarities between the standard and the proposed approach}\label{sim}

We emphasize that the regularizing LM scheme presented in Section \ref{ir} is designed to compute stable approximation to a minimizer of $\Phi$ defined in (\ref{eq:1.1}). In contrast, the standard approach discussed in the preceding section is based on minimizing (\ref{eq:sa}) by means of a standard optimization algorithm. Therefore, the two approaches aim at solving two substantially different problems. Nonetheless, there are some computational similarities between the aforementioned approaches as we now discuss. Let us recall that each iteration step for minimizing (\ref{eq:sa}) in the standard approach is given by (\ref{eq:3.7}) which, from the lemmas in \cite[Section 7.4]{Oliver}, is equivalent to
\begin{eqnarray}\label{eq:3.9}
\Delta u=
C \, DG^{*}(u_{m}) \Big[ DG(u_{m}) \, C \, DG^{*}(u_{m})+(1+\lambda_{m})\Gamma\Big]^{-1}\Big[y^{\eta}-G(u_{m})\nonumber\\+\frac{1}{1+\lambda_{m}} DG(u_{m})(u_{m}-\overline{u})\Big]
+\frac{1}{1+\lambda_{m}}(u_{m}-\overline{u}).
\end{eqnarray}
On the other hand, the $m$th step computed with the regularizing LM scheme for approximating the minimizer of (\ref{eq:1.1}) is given by 
\begin{eqnarray}\label{eq:4.21}
\Delta u_{m}^{\eta}=C \, DG^{*}(u_{m}^{\eta})[DG(u_{m}^{\eta}) \, C\, DG^{\ast}(u_{m}^{\eta})+\alpha_{m} \Gamma]^{-1}[y^{\eta}-G(u_{m}^{\eta})]
\end{eqnarray}
The substantial similarities between expression (\ref{eq:3.9}) and (\ref{eq:4.21}) are evident although they converge to different functions. Notice that, at the discretization level, the computation of the matrices $C \, DG^{*}(u)$ and $DG(u)\, C\, DG^{\ast}(u)$ as well as the evaluation of $G(u)$ are needed for both approaches (obviously evaluated at different $u$'s). In addition, note that the terms $1/(1+\lambda_{m}) DG(u_{m})(u_{m}-\overline{u})$ and $1/(1+\lambda_{m})DG(u_{m})(u_{m}-\overline{u})\Big]$ are not required in (\ref{eq:4.21}). It is therefore clear that the main routines and codes used for computing (\ref{eq:3.9}) in the standard approach can be used for implementing the regularizing LM scheme step (\ref{eq:4.21}). In fact, a routine that assembles $C \, DG^{*}(u)$ and $DG(u) \, C\, DG^{\ast}(u)$ as well as the routine that evaluates $G(u)$ are sufficient for a straightforward implementation of the regularizing LM scheme Algorithm \ref{al:LM2}. The aforementioned computational similarities open the possibility to study the history matching problem in the sense presented in this paper by using available implementations for the standard approach. Moreover, provided that the same implementation for $C \, DG^{*}(u)$, $DG(u)\, C\, DG^{\ast}(u)$ and $G(u)$ are used for both approaches, from the previous discussion it follows that the two approaches have the same computational cost per iteration. For the results presented in the subsequent section, the operator $DG$ and $DG^{\star}$ are computed as described in Section 9.7 of \cite{Oliver}.

\section{Numerical Results }\label{Numerics}

In this section we present numerical experiments to show the capabilities of the regularizing LM scheme for estimating the log-permeability in the oil-water model of the Appendix. 

\subsection{Experimental setting}\label{expset}

We consider a synthetic experiment where the reservoir domain is $\Omega=[0,L]\times [0,L]$ and the prior knowledge of the subsurface is given in terms of a prior $\overline{u}=500~\textrm{md}$ (constant in $\Omega$) and a covariance operator 
\begin{eqnarray}\label{eq:3.10}
C =\kappa^{-1} C_{0}
\end{eqnarray}
where $C_{0}$ is a spherical covariance function \cite{Geos}
\begin{eqnarray}\label{eq:3.11}
C_{0}(z_1,z_2) =\left\{\begin{array}{cc}
1-\frac{3}{2}\frac{\vert\vert M_{\theta}(z_{1} -z_{2})\vert\vert }{a}+\frac{1}{2}\frac{\vert\vert M_{\theta}(z_{1} -z_{2})\vert\vert^3}{a^3}&\textrm{if}~~\vert\vert M_{\theta}(z_{1} -z_{2})\vert\vert <a\\
0 & \textrm{if}~~\vert\vert M_{\theta}(z_{1} -z_{2})\vert\vert \ge a\end{array}\right.
\end{eqnarray}
with $z_{i}=(x_{i},y_{i})$. In the previous expression, $M_{\theta}$ is a rotation matrix along the direction of maximum continuity with range denoted by $a$. Covariance functions like (\ref{eq:3.11}) are common in modeling geologic properties of reservoirs \cite{Geos}. The tunable parameter $\kappa$ in (\ref{eq:3.11}) will enable us to study the performance of the proposed approach with respect to different choices of the prior covariance parameterized in terms of $\kappa$.

\begin{figure}
\includegraphics[scale=0.35]{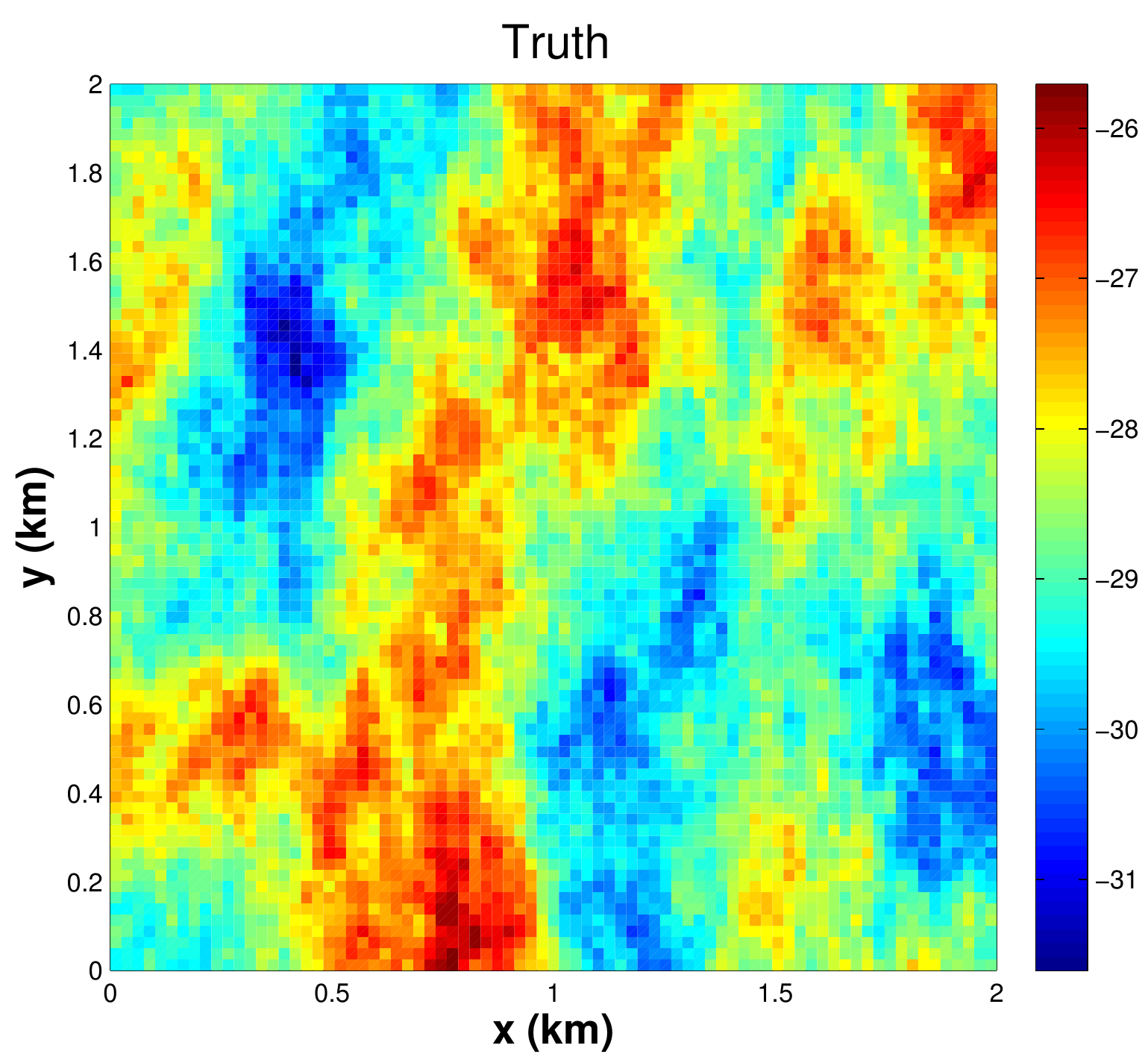}
\includegraphics[scale=0.35]{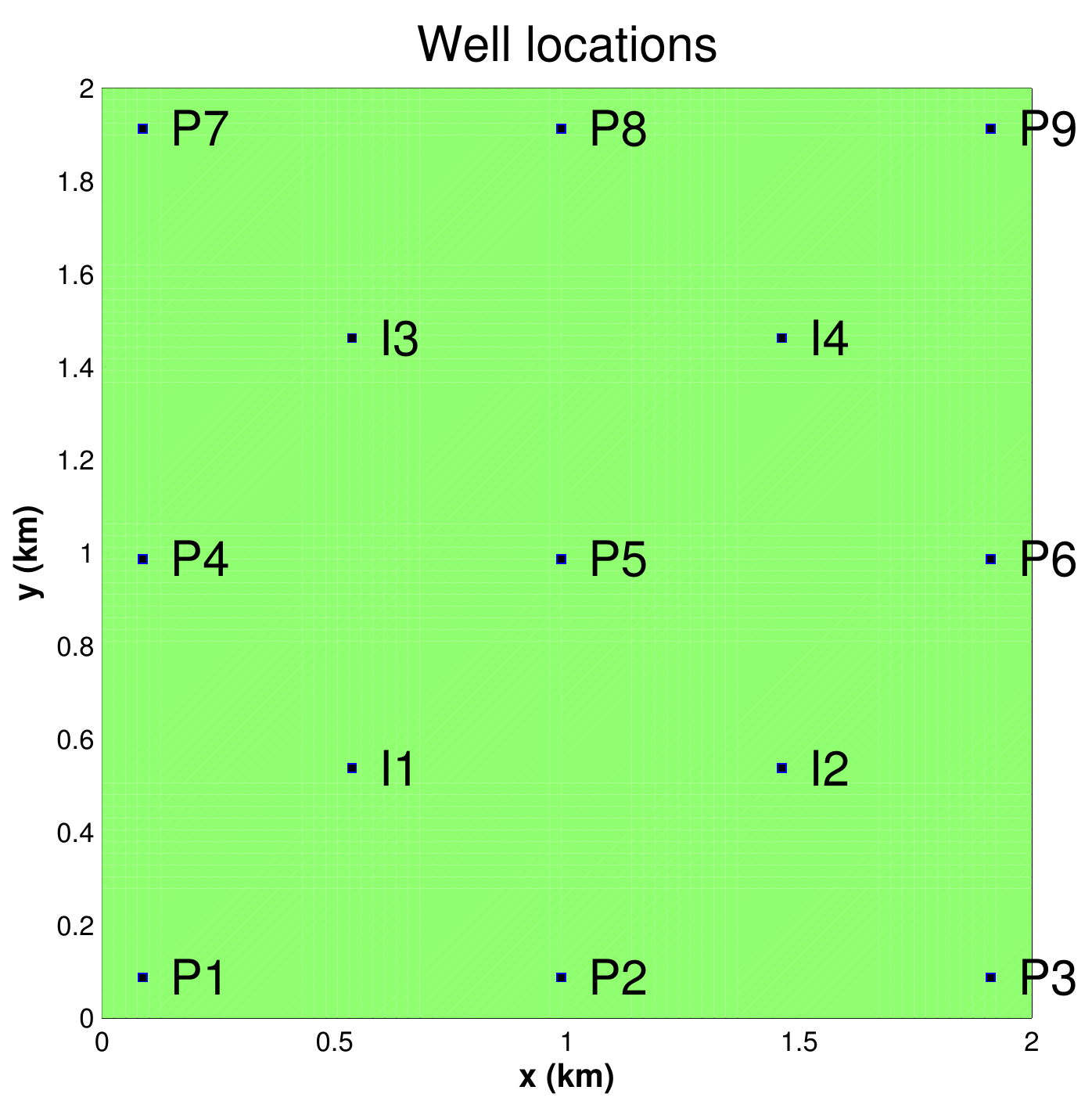}
\caption{Left: True log-permeability [$\log{\textrm{m}^2}$]. Right: Well configuration.}  
\label{Figure1}
\end{figure}

We consider $\kappa=1$ in (\ref{eq:3.10}) to be the ``correct'' covariance in the sense that the true (or reference) log-permeability is a Gaussian field with mean $\overline{u}$ and covariance $C=C_{0}$. In other words, $\kappa=1$ corresponds to the best case scenario where with our prior knowledge is consistent with the truth. In Figure \ref{Figure1} (left) we display the true permeability $u^{\dagger}$, sampled from the aforementioned distribution. We now consider a water flood described with the model presented in the Appendix. Nine production wells $P_{1},\dots,P_{9}$ and four injection wells $I_{1},\dots, I_{4}$ are considered in the configuration displayed in Figure \ref{Figure1} (right). Relevant data of the reservoir model is displayed in Table \ref{Table1}. We use the true log-permeability field of Figure \ref{Figure1} (left) to generate synthetic data as we now describe. First, the PDE system (\ref{eq:2.7})-(\ref{eq:2.7B}) is solved for $u=u^{\dagger}$, the resulting pressures and saturations are used in the expression for the measurement functional (\ref{eq:2.11})-(\ref{eq:2.14}) to find $G(u^{\dagger})$. Finally, synthetic data is generated by adding Gaussian random noise $\xi\sim N(0,\Gamma)$. More precisely, we define $y^{\eta}\equiv G(u^{\dagger})+\xi$. We consider a diagonal error measurement covariance $\Gamma$ with diagonal elements denoted by $\sigma_{i}^2$. The values of $\sigma_{i}$ associated to measurements of bottom hole pressure consist of some percentage (defined below) of the nominal value of the corresponding measured variable. In order to avoid zero values for the $\sigma_{i}$'s associated to measurements of water rates, for either water and oil rate measurements, the corresponding $\sigma_{i}$ is a percentage of the nominal value of the total flow rate (which is the well constraint). The aforementioned percentage is the same for both measurements of pressure and flow rates. The noise level is defined by
\begin{eqnarray}
\eta \equiv \vert \vert \Gamma^{-1/2}(  y^{\eta}-G(u^{\dagger})) \vert\vert_{Y}
\end{eqnarray}

\begin{table}
\caption{Reservoir model description}
\label{Table1}       
\begin{tabular}{lc|lcc}
\hline\noalign{\smallskip}
Variable  & Value&Variable  & Value\\
\noalign{\smallskip}\hline\noalign{\smallskip}
L  [$\textrm{m}^{3}$] & $2\times 10^{3}$ &$a_{w}$ & $0.3$\\
$c$  [$\textrm{Pa}^{-1}$] &  0.0 & $a_{o}$   & $0.9$  \\
$\nu_{o}$  [$\textrm{Pa s}$] & $10^{-2}$  &$^{\rm b}$ $P_{bh}^{l}$ [\textrm{Pa}]  & $2.7\times 10^{7}$ \\
$T$  [$\textrm{years}$] & 5 &$^{\rm b}$ $q_{w}^{l}$  [$\textrm{m}^{3}/\textrm{day}$]  & $2.6\times 10^3$  \\
$^{\rm a}$ $p_{0}$ [\textrm{Pa}] &$2.5\times 10^{7}$ &$s_{iw}$  & 0.2 \\
$^{\rm a}$ $s_{0}$   & 0.2 &$s_{ro}$    & 0.2 \\
$\nu_{w}$  [$\textrm{Pa s}$] & $5\times 10^{-4}$  \\
\noalign{\smallskip}\hline
\end{tabular}

$^{\rm a}$ Constant in $\Omega$. $^{\rm b}$ Constant in $[0,T]$.
\end{table}

\subsection{Performance of the LM scheme with respect to the observational noise level}\label{noisesec}

In this subsection we investigate the accuracy of the estimate $u^{\eta}$ obtained with the regularizing LM scheme as a function of the noise level $\eta$ with $\eta\to 0$. According to Theorem \ref{Hanke_Theorem}, $u^{\eta}$ converges to a minimizer of (\ref{eq:1.1}) as $\eta\to 0$ (see Remark \ref{rem3}). Although this converged solution may not necessarily be the truth $u^{\dagger}$ (due to possible non-uniqueness), the solutions may arguably reflect the main spatial features of the truth. We therefore consider the accuracy of the estimates in terms of their relative error with respect to the true log-permeability $u^{\dagger}$. 

Five sets of synthetic data associated to different noise levels $\{\eta_{j}\}_{j=1}^{5}$ are generated with the procedure previously described. For each set, a different percentage of the nominal value of the measured values is selected. The resulting sets of synthetic data $\{y^{\eta_{j}}\}_{j=1}^{5}$ provide noise levels (defined by (\ref{eq:noise})) that correspond to some fractions of the norm of the corresponding measurements $  \vert \vert \Gamma^{-1/2}  y^{\eta_{j}} \vert\vert_{Y} $, $j\in\{1,\dots,5\}$. More precisely, we have 
\begin{eqnarray}\label{eq:noise3}
\eta_{j} \equiv  f_{j}\vert \vert \Gamma^{-1/2}  y^{\eta_{j}} \vert\vert 
\end{eqnarray}
with $f_{1}=5\times 10^{-2}$, $f_{2}=10^{-2}$, $f_{3}=5\times 10^{-3}$, $f_{4}=10^{-3}$ and $f_{5}=5\times 10^{-4}$. 

\begin{figure}
\includegraphics[scale=0.35]{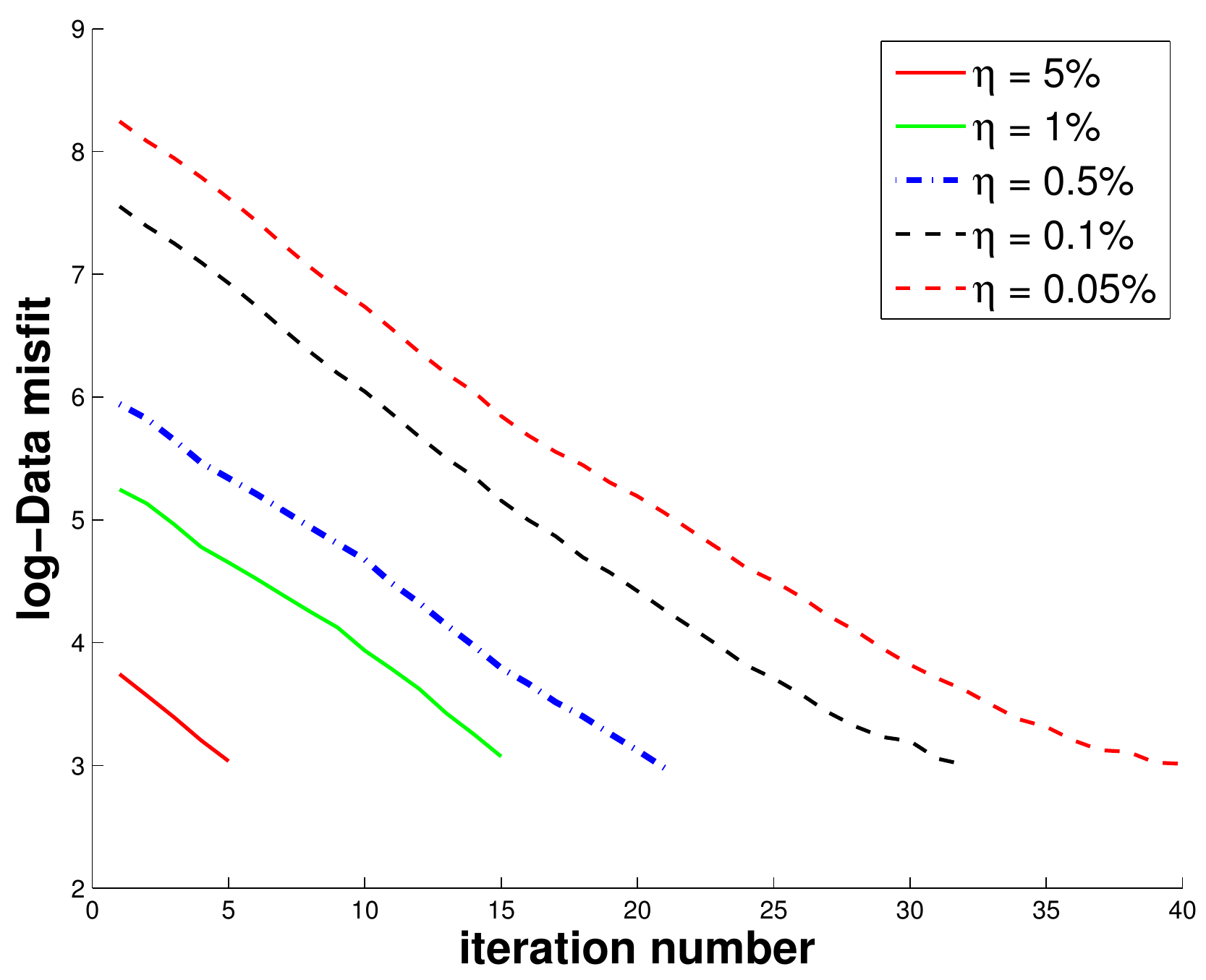}
\includegraphics[scale=0.35]{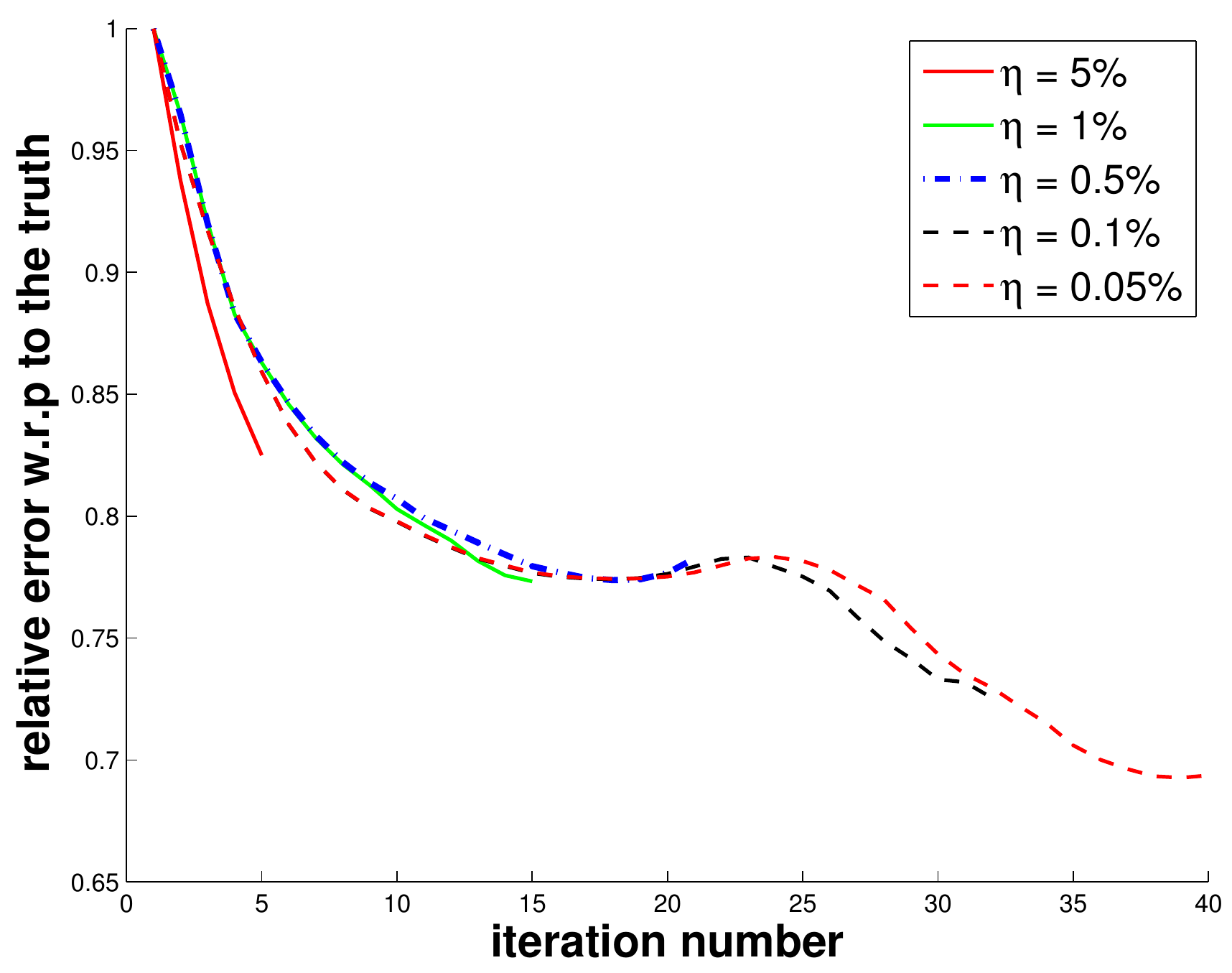}
\caption{Performance of the regularizing LM scheme with respect to the noise level . Right: data misfit. Left: relative error with respect to the truth}  
\label{Figure8}
\end{figure}

\begin{figure}
\includegraphics[scale=0.25]{./True}
\includegraphics[scale=0.25]{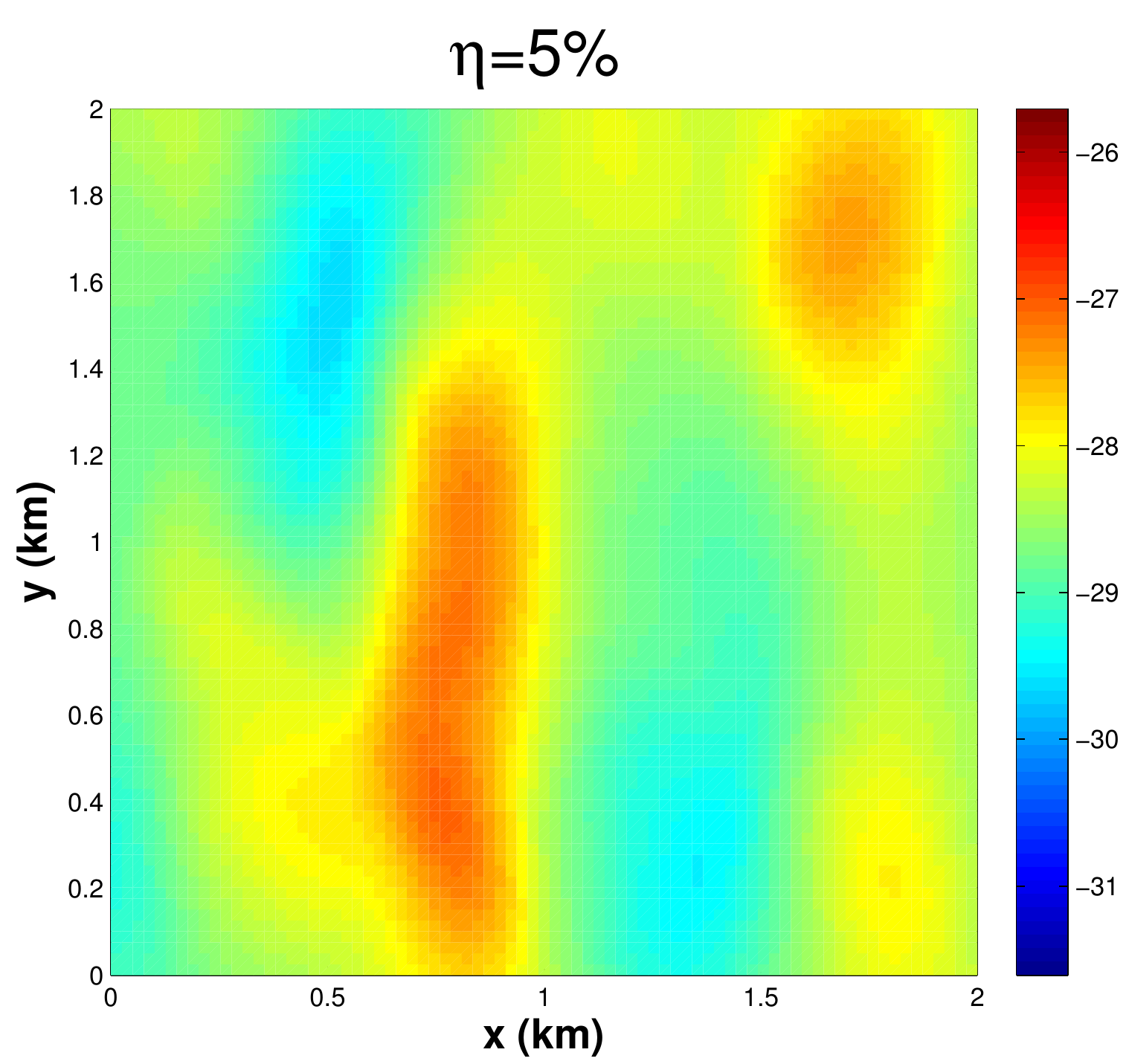}
\includegraphics[scale=0.25]{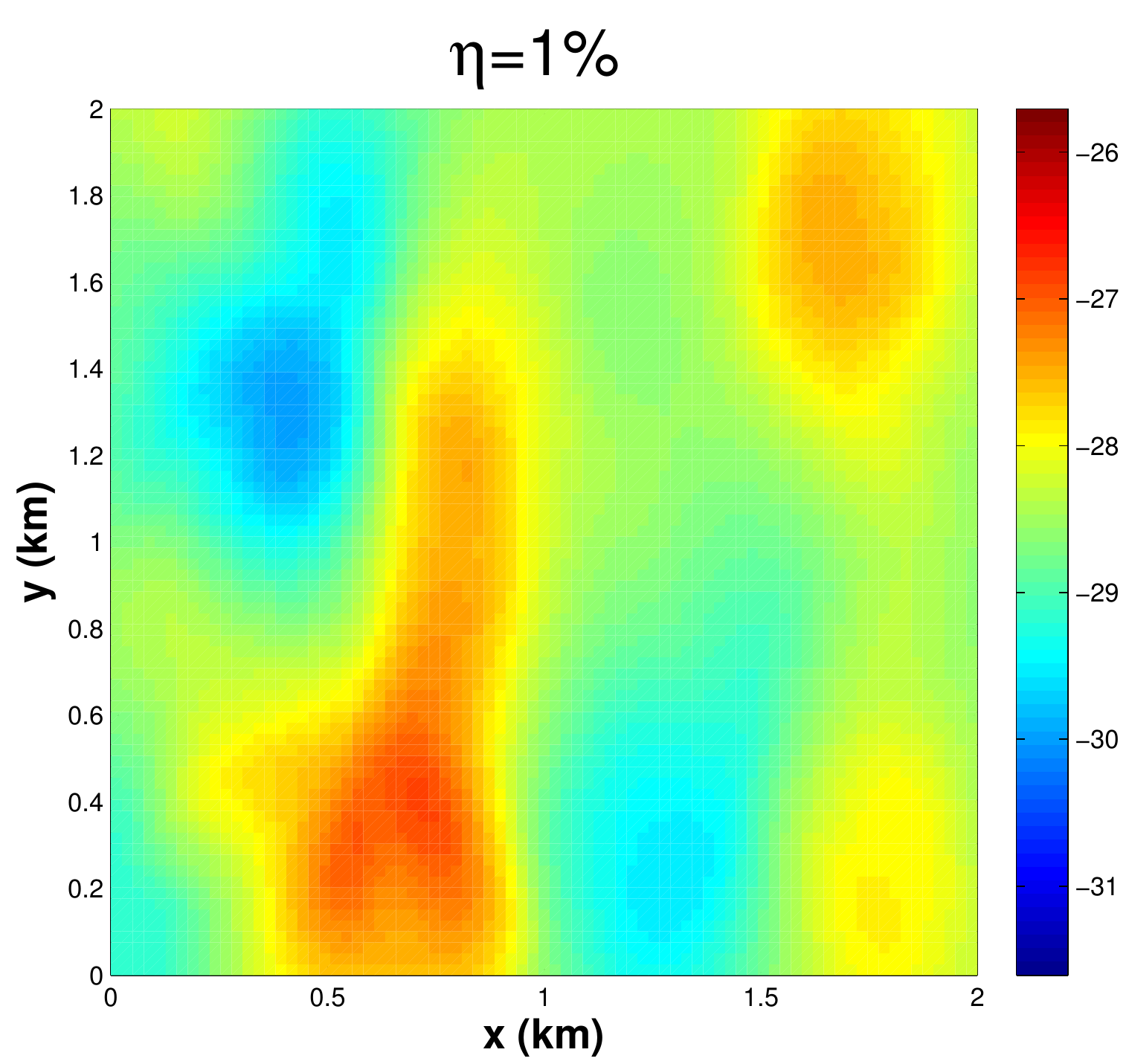}\\
\includegraphics[scale=0.25]{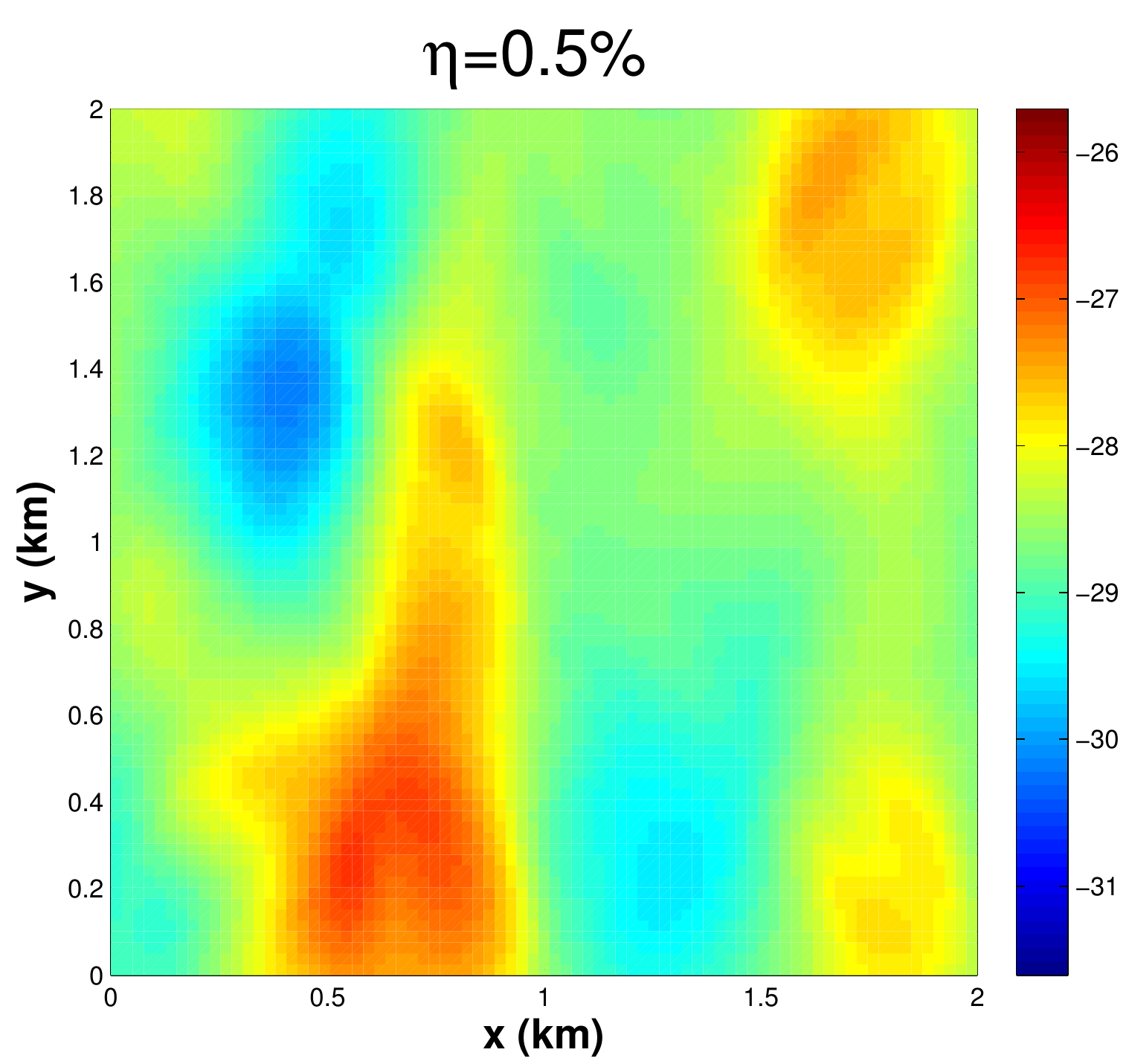}
\includegraphics[scale=0.25]{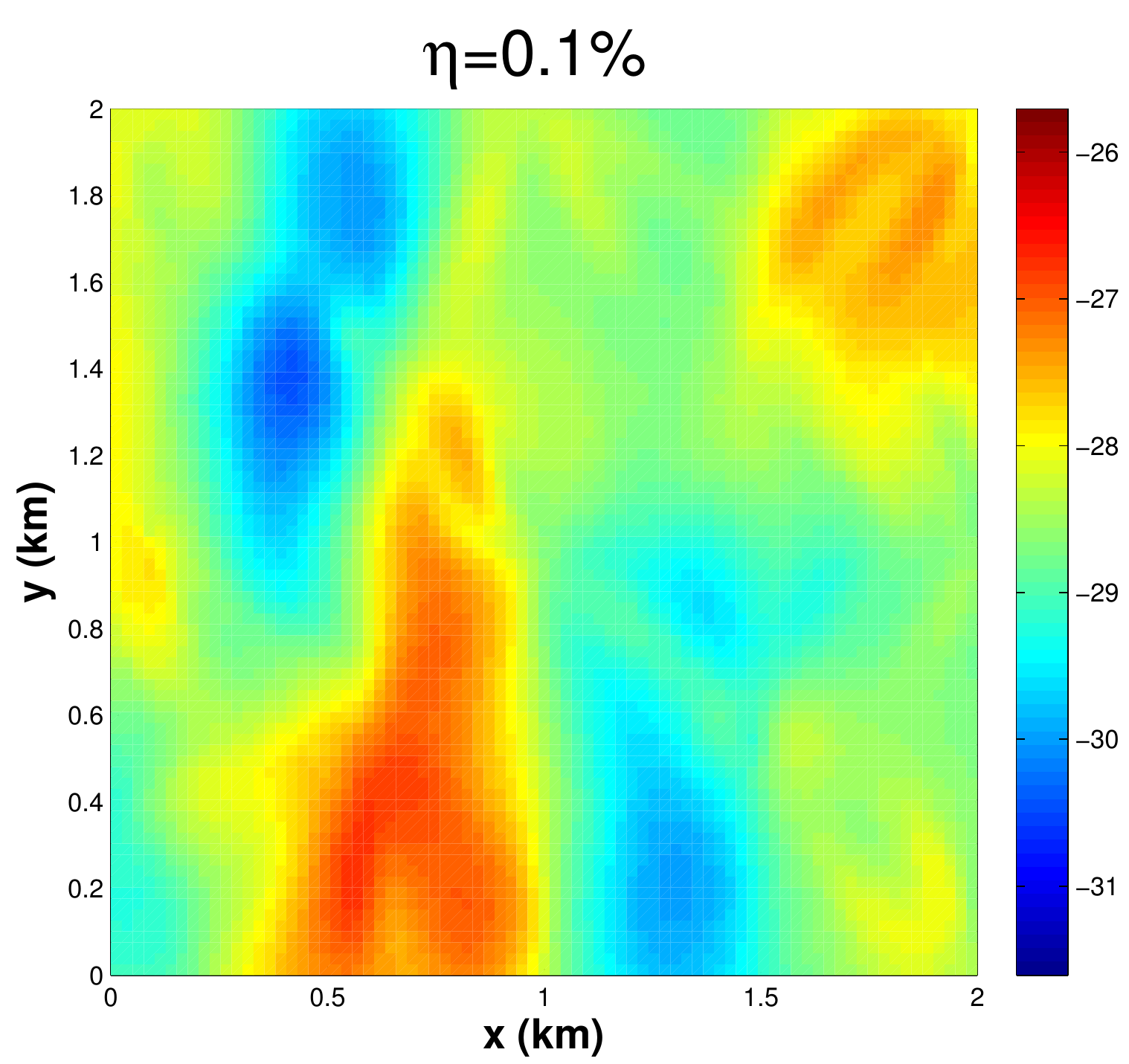}
\includegraphics[scale=0.25]{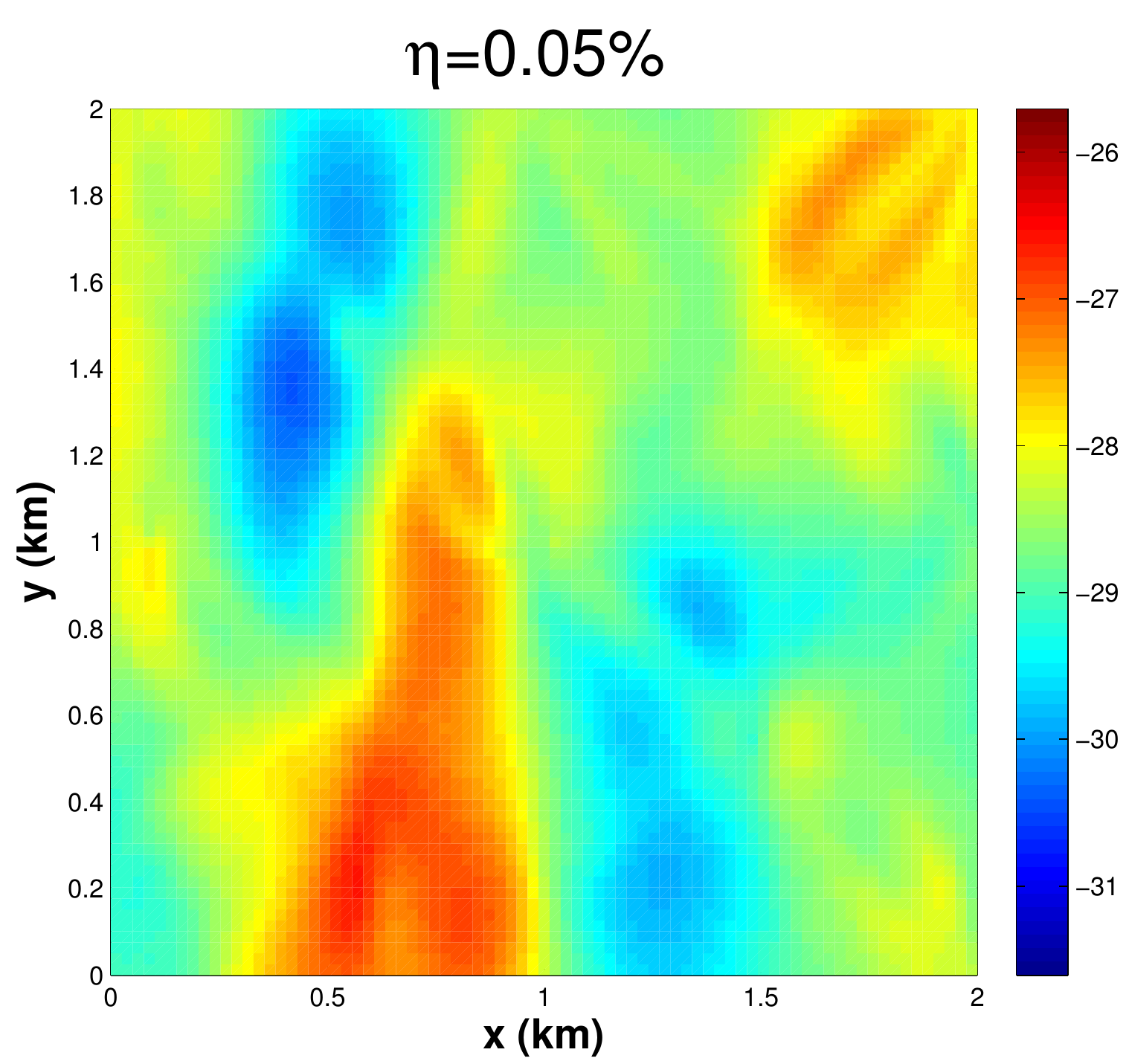}
\caption{Log-permeability estimates obtained with the regularizing LM scheme for different noise levels  [$(\log{\textrm{m}^2})$]}   \label{Figure9}
\end{figure}

For this set of experiments, the parameters for the regularizing LM scheme are selected as $\tau=1.2$ and $\rho=0.83$. Further choices of $\rho$ and $\tau$ are investigated in subsection \ref{taurho}. In addition, we consider $\kappa=1$ in (\ref{eq:3.10}).The performance of the regularizing LM scheme for each of the five sets of synthetic data corresponding to different noise levels is presented in Figure \ref{Figure8}. The data misfit $\vert \vert \Gamma^{-1/2}(  y^{\eta}-G(u^{\eta})) \vert\vert_{Y}$ is displayed in Figure \ref{Figure8} (left) and the relative error with respect to the truth $\vert\vert u_{m+1}-u_{m}\vert\vert_{X} /\vert\vert u_{m+1} \vert\vert_{X} $ is shown in Figure \ref{Figure8} (right). The stability in the computation of the numerical solutions is reflected in the decrease of the relative error with respect to the truth. Note that, as the noise level decreases, the accuracy with respect to the relative error increases. The dependence of the accuracy on the noise level can be visually appreciated from the log-permeability estimates presented in Figure \ref{Figure9}. For smaller noise in the observations, the regularizing LM scheme seems to provide stable and accurate estimates of the geologic properties.

By construction, the weight $\Gamma^{-1/2}$ (in the data misfit) depends inversely on the error $y^{\eta}-G(u^{\dagger})$. Therefore, even though the five experiments have the same initial guess $u_{0}=\overline{u}$, the initial value of the data misfit is larger for smaller noise levels. Furthermore, since the error that defines the noise level (\ref{eq:noise}) is also weighted by $\Gamma^{-1/2}$, the actual value in (\ref{eq:noise3}) is approximately similar for all the experiments. The difference, however, is in the corresponding fraction $f_{j}$ of the norm in (\ref{eq:noise3}) which is used in the label of Figure \ref{Figure8}.

\subsection{Parameters $\tau$ and $\rho$}\label{taurho}

The parameters $\rho\in (0,1)$ and $\tau>1/\rho$ are the tunable parameters in the regularizing LM scheme. In this section we present the numerical performance of the LM scheme for different choices of these parameters. We recall from Section \ref{ir} that if $\rho\approx 1$ ($\rho<1$), we then may choose $\tau\approx 1$ ($\tau>1/\rho$). Then the regularizing LM scheme terminates when the estimate $u_{m}^{\eta}$ produces a data misfit $\vert\vert \Gamma^{-1/2}(y^{\eta}-G(u_{m}^{\eta})\vert\vert_{Y}\approx \eta$. Small values of $\rho$ imply larger values of $\tau$ which may lead to estimates that provide a poor fit to the production data. It is therefore important to study the potential lack of accuracy due to the choices of $\rho$ and $\tau$.

\begin{figure}
\includegraphics[scale=0.35]{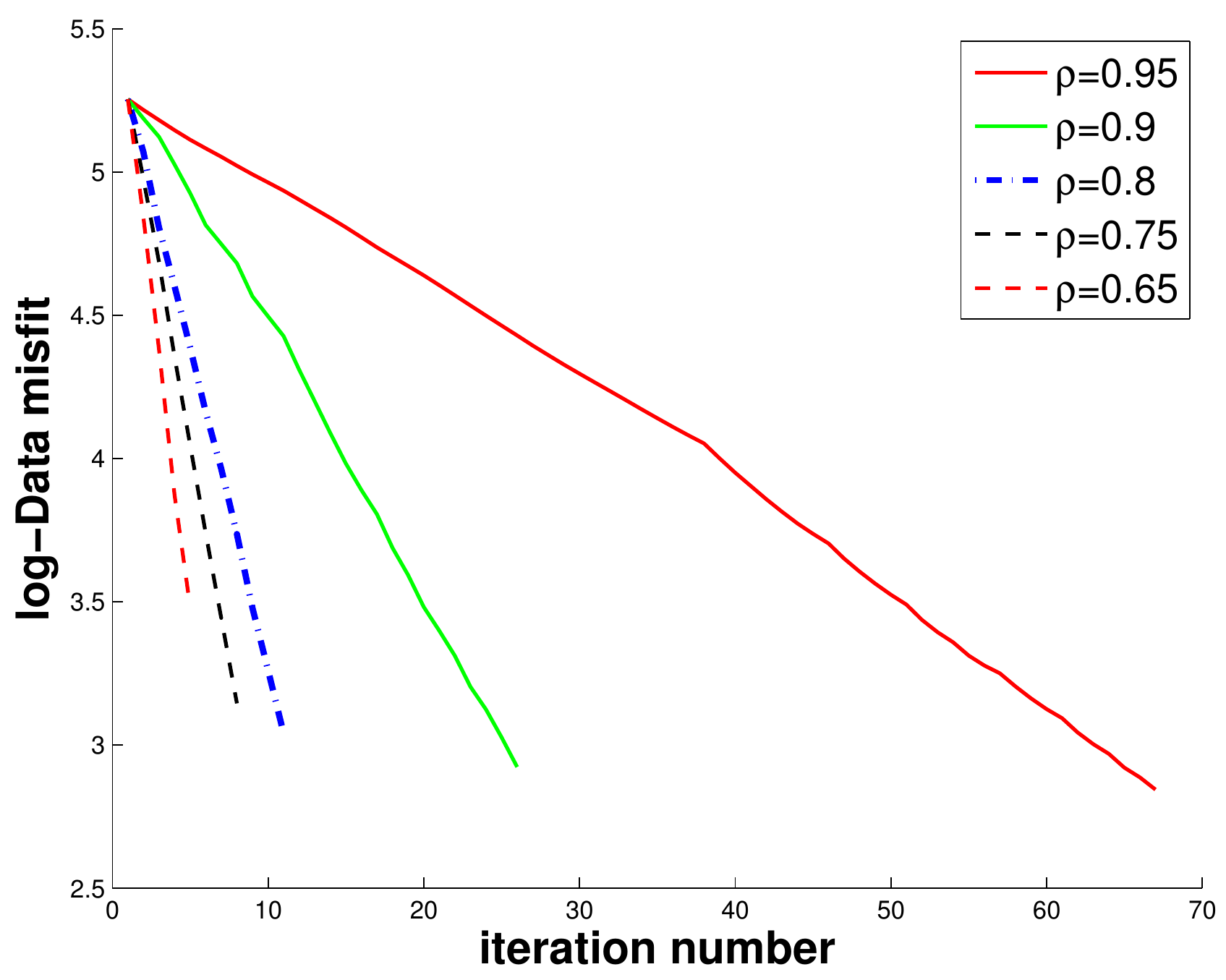}
\includegraphics[scale=0.35]{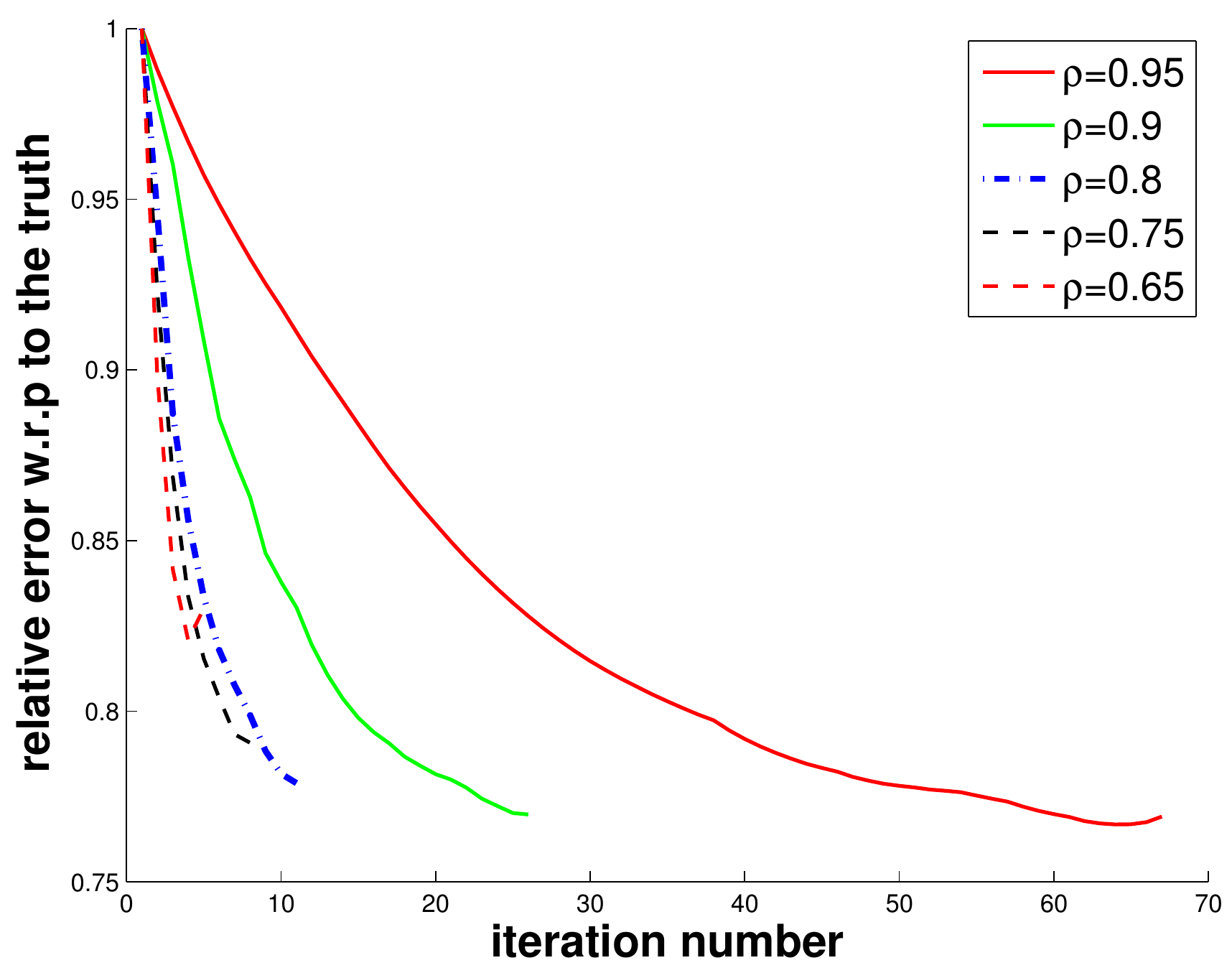}
\caption{Performance of the regularizing LM scheme with respect to the parameter $\rho$. Right: data misfit. Left: relative error with respect to the truth}  
\label{Figure10}
\end{figure}

\begin{figure}
\includegraphics[scale=0.25]{./True}
\includegraphics[scale=0.25]{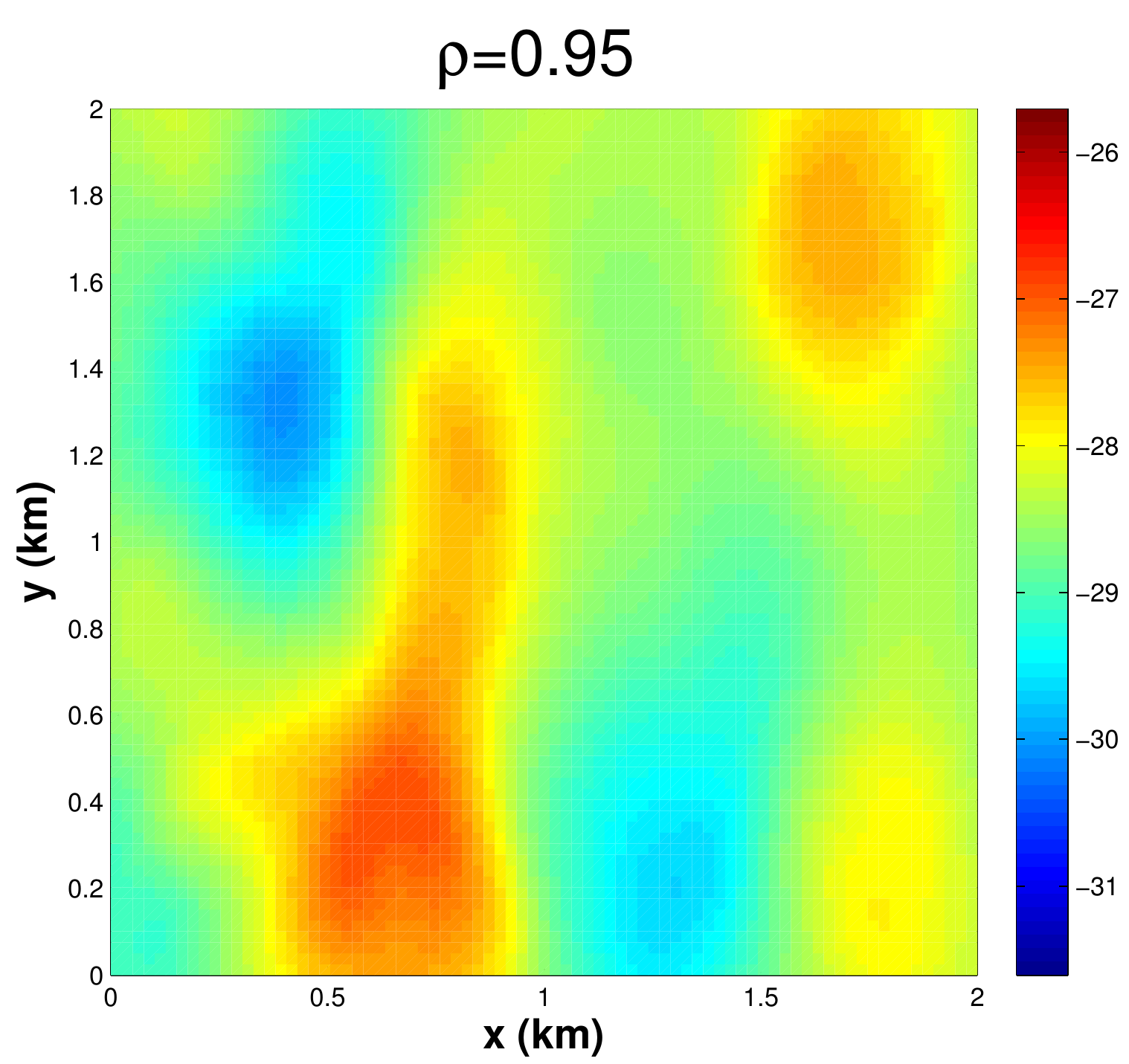}
\includegraphics[scale=0.25]{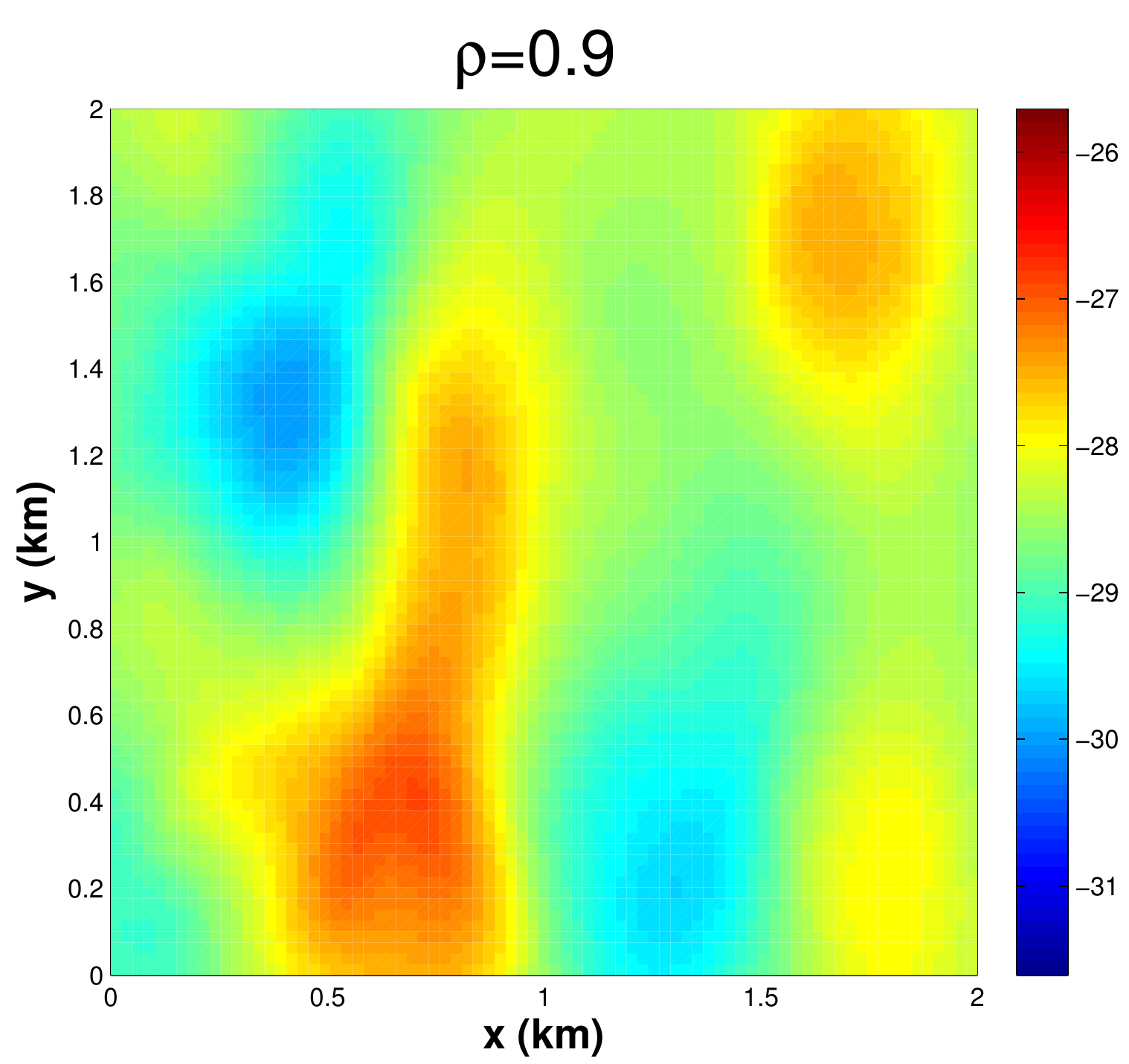}\\
\includegraphics[scale=0.25]{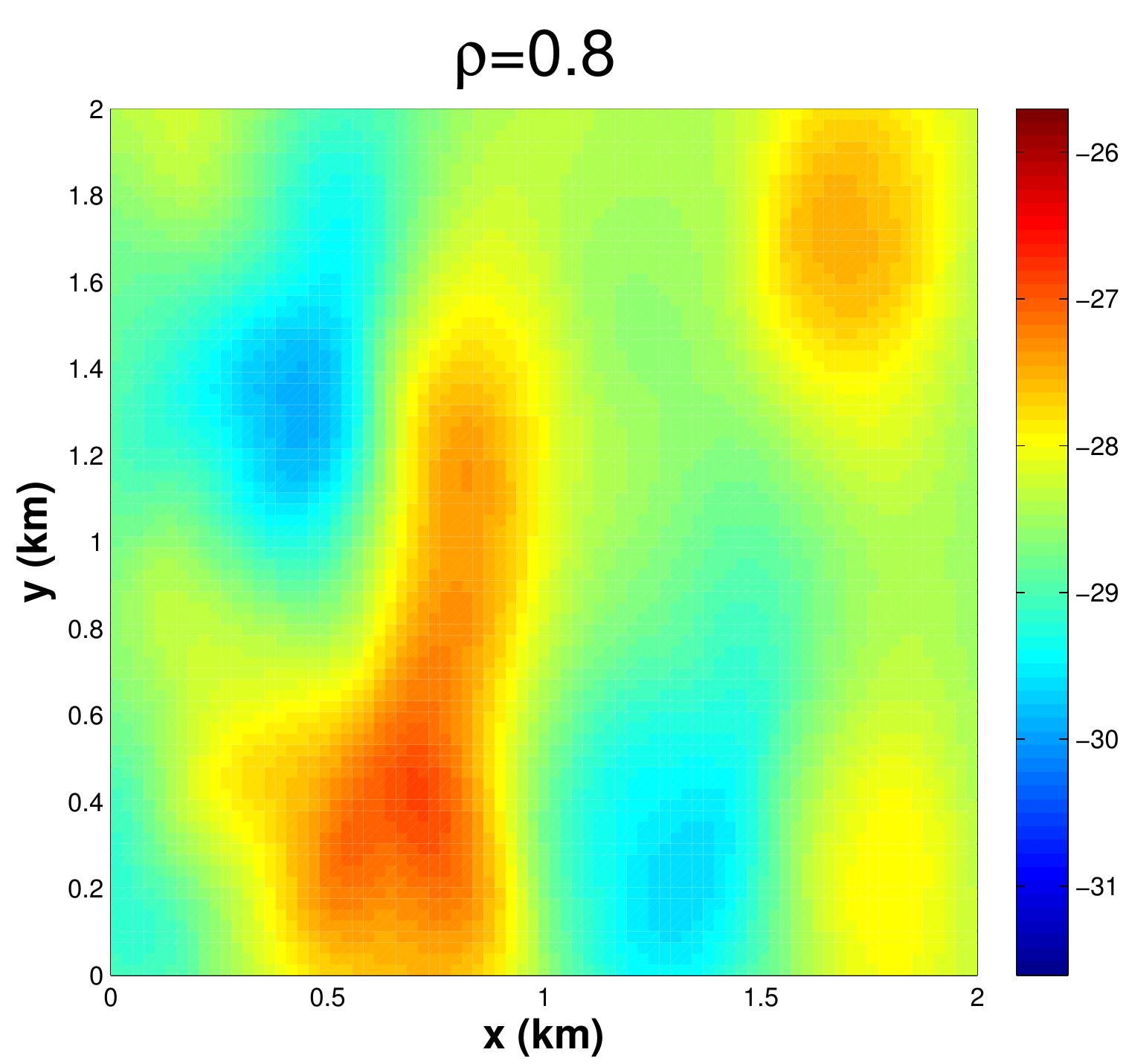}
\includegraphics[scale=0.25]{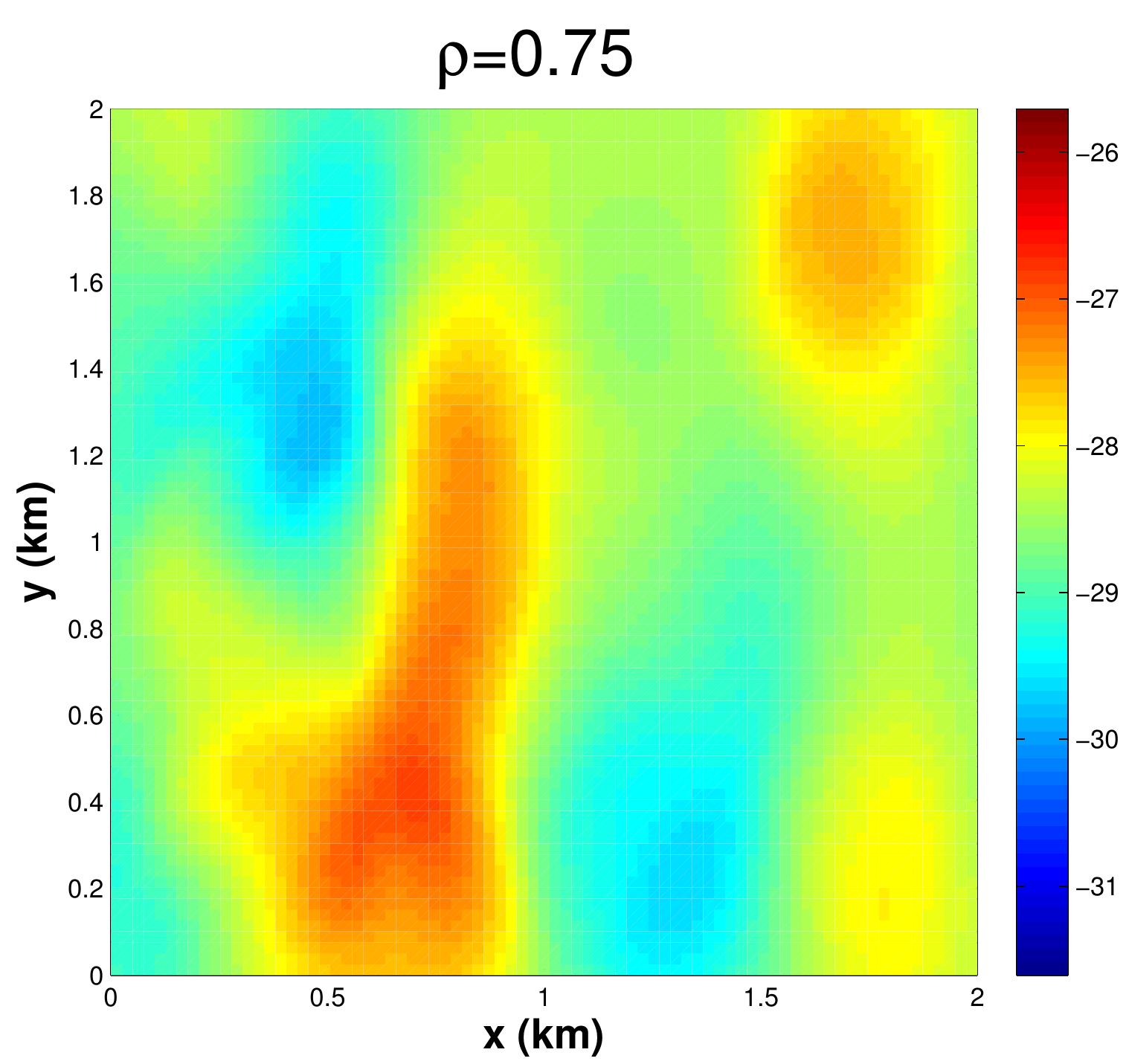}
\includegraphics[scale=0.25]{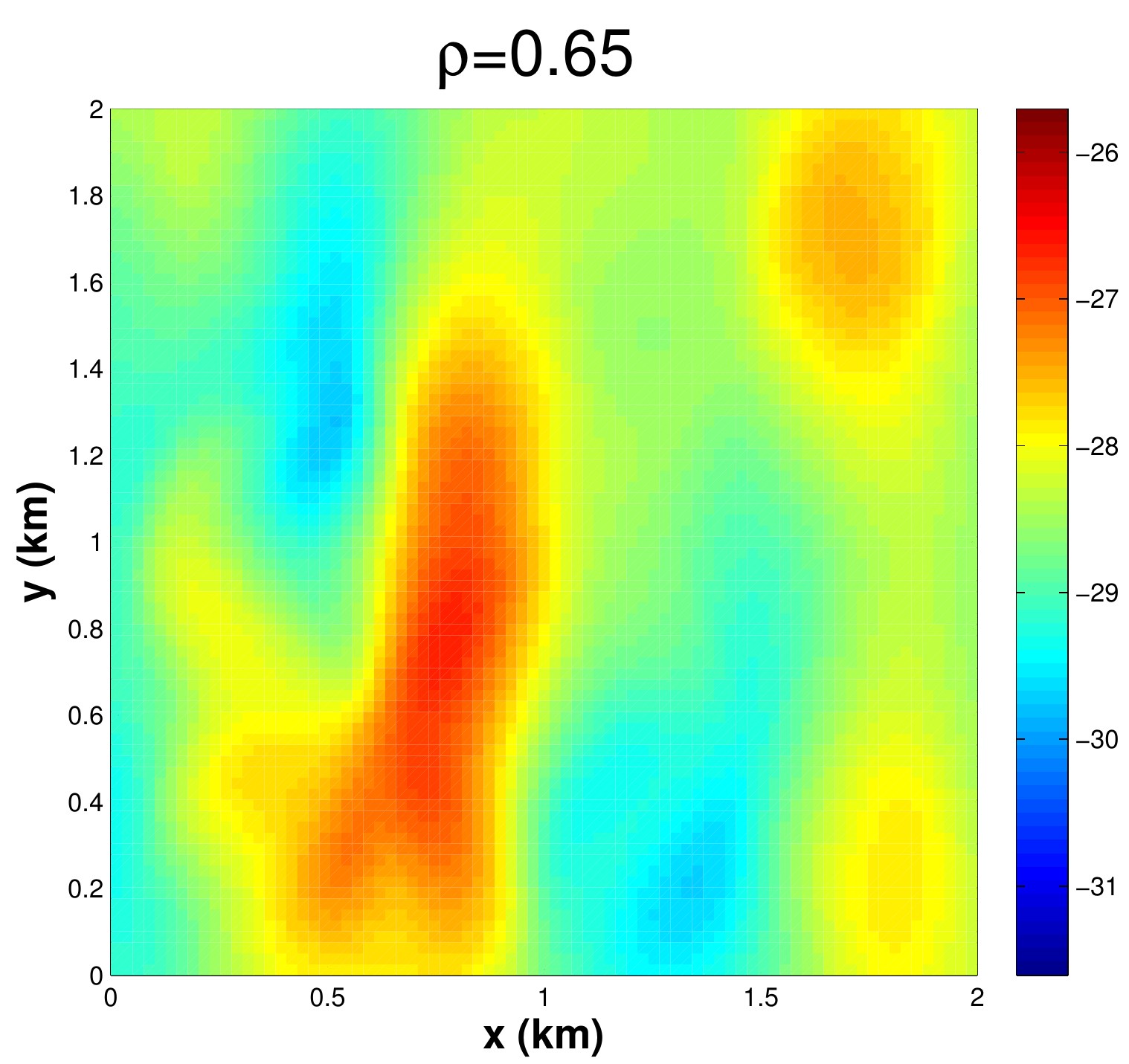}
\caption{Log-permeability estimates obtained with the regularizing LM scheme for different parameter $\rho$.  [$(\log{\textrm{m}^2})$]}  
\label{Figure11}
\end{figure}

We consider the same experimental setting as before for only one fixed set of synthetic data with $1\%$ of observational noise level. We consider several choices of $\rho$, with the corresponding $\tau$ defined by $\tau=1/(\rho-10^{-3})$. The performance of the LM scheme for these choices of parameters is presented in Figure \ref{Figure10} and Figure \ref{Figure11}. Although reasonable estimates were obtained for all these choices of $\rho$, from Figure  \ref{Figure10} we observe that more accurate estimates, in terms of the relative error with respect to the truth, are obtained when $\rho$ is indeed close to one. However, it is important to remark that an increase in the computational cost is associated with the improved accuracy for $\rho\approx 1$. These numerical experiments suggest that optimal choices in terms of computational efficiency and accuracy are obtained for $\rho\in [0.8, 0.9]$.

\subsection{Performance of the LM scheme with respect to the prior covariance. }\label{eq:ir1}

Recall that in the previous experiments we have chosen $\kappa=1$ in (\ref{eq:3.10}) which corresponds to the best-case scenario where the true log-permeability is consistent with the prior knowledge. In this subsection we investigate the performance of the regularizing LM scheme with respect to different choices of the prior covariance. In particular, we consider the case where the prior covariance is parameterized by (\ref{eq:3.10}), and we are interested in the performance of the LM scheme with respect to $\kappa$. The values of the parameters $\tau$ and $\rho$ in the regularizing LM scheme are as described in the subsection \ref{noisesec} and the synthetic data as the same as in subsection \ref{taurho}. 

\begin{figure}
\includegraphics[scale=0.35]{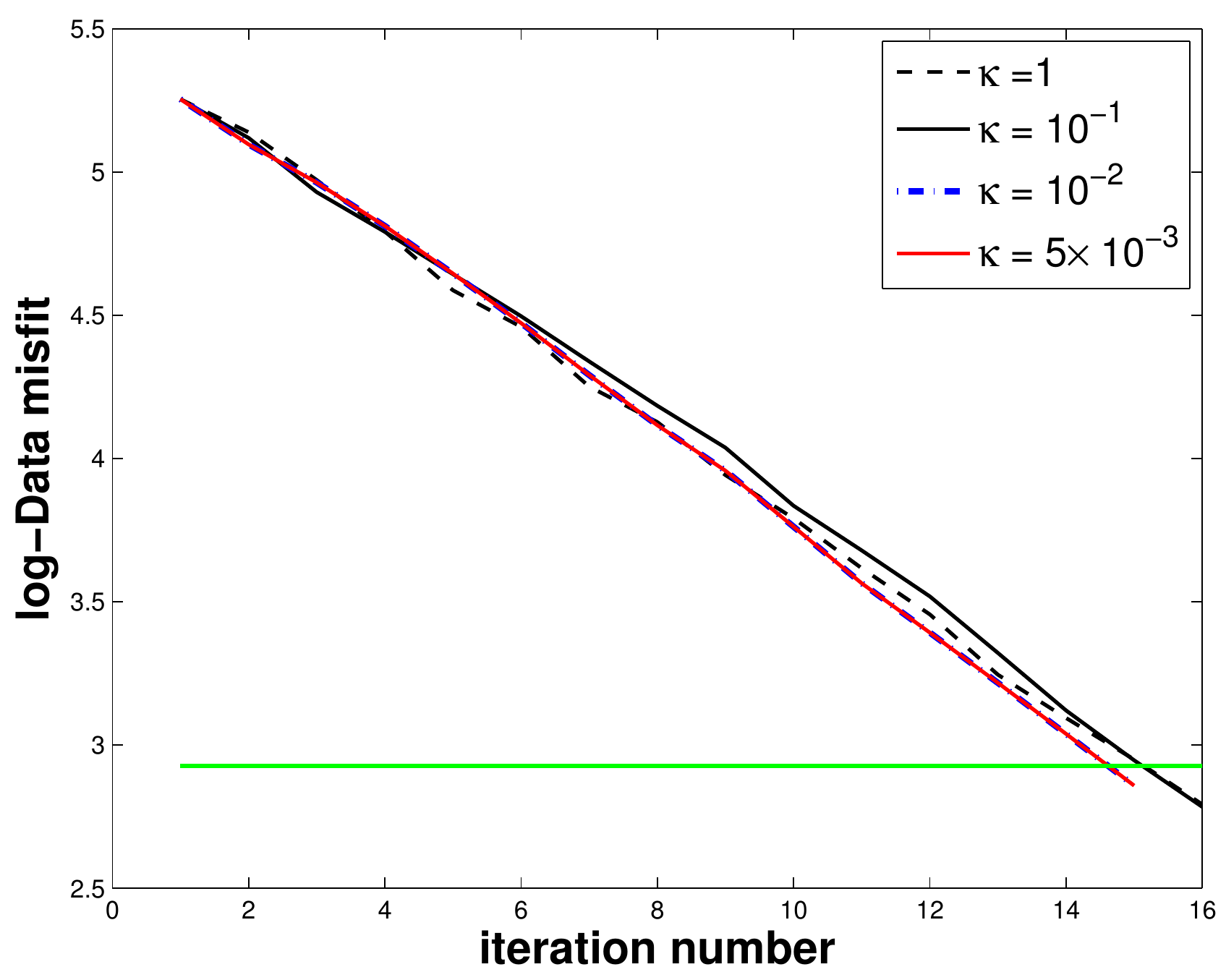}
\includegraphics[scale=0.35]{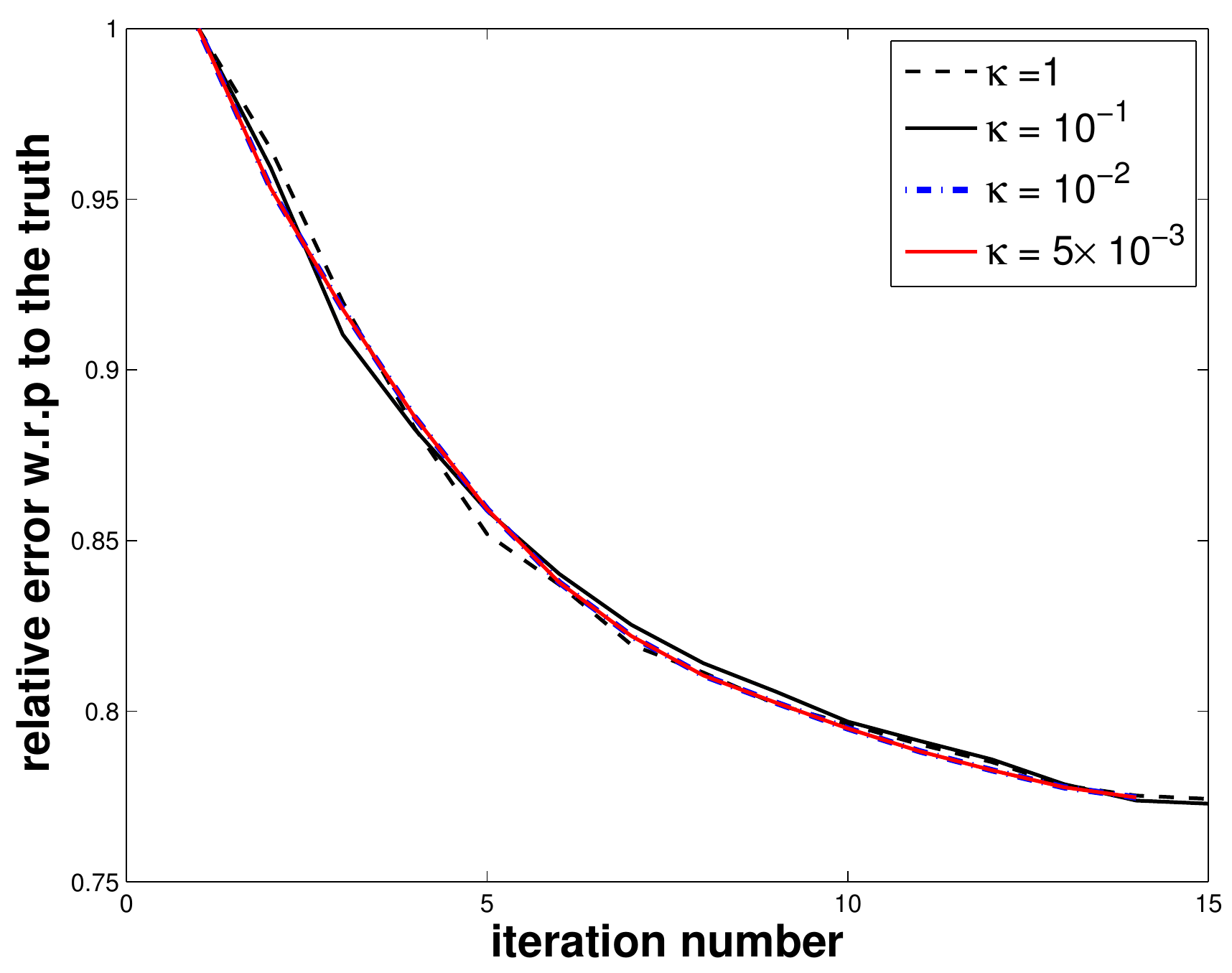}
\caption{Regularizing LM scheme for history matching. Right: data misfit. Left: relative error with respect to the truth}  
\label{Figure5}
\end{figure}

\begin{figure}
\includegraphics[scale=0.35]{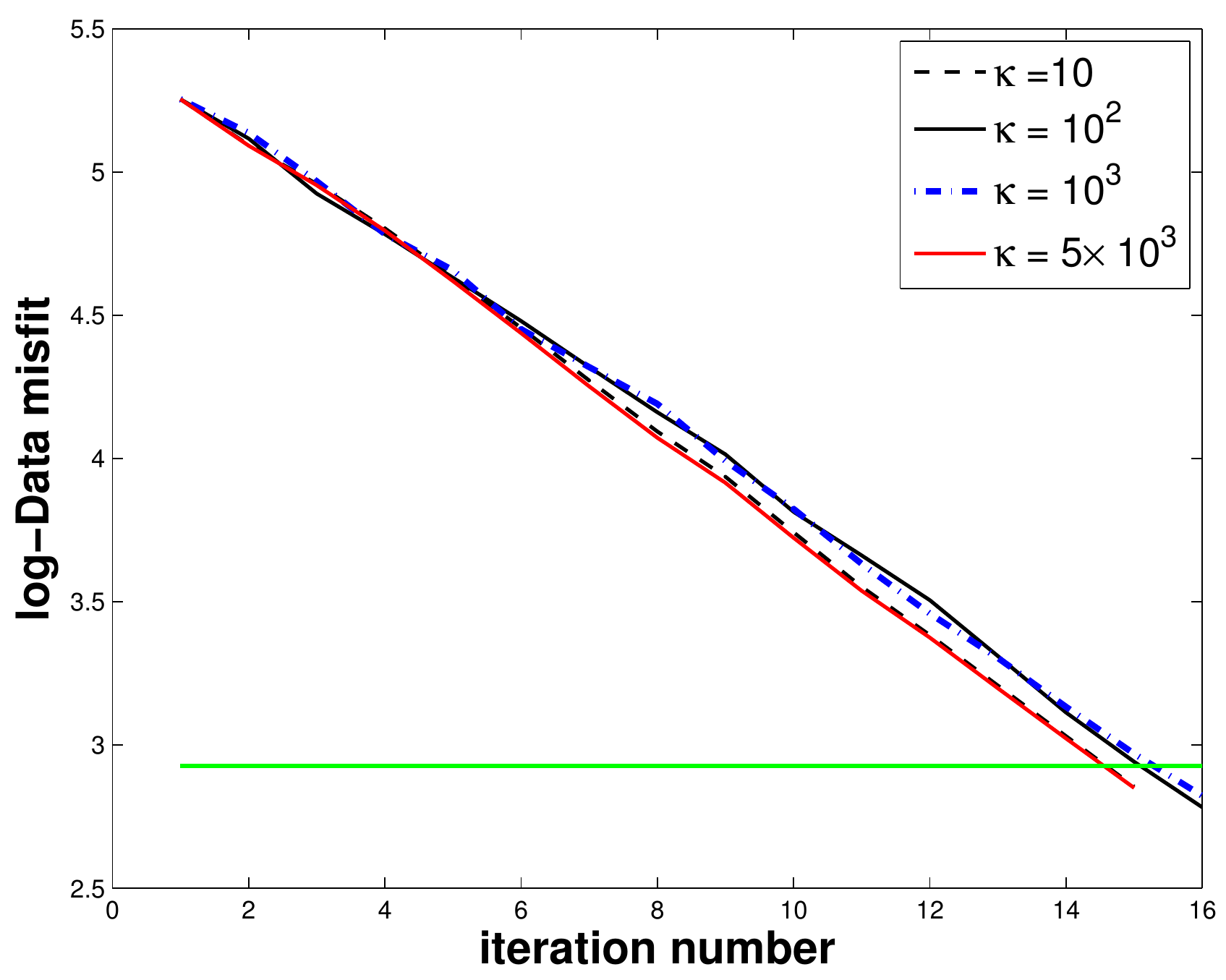}
\includegraphics[scale=0.35]{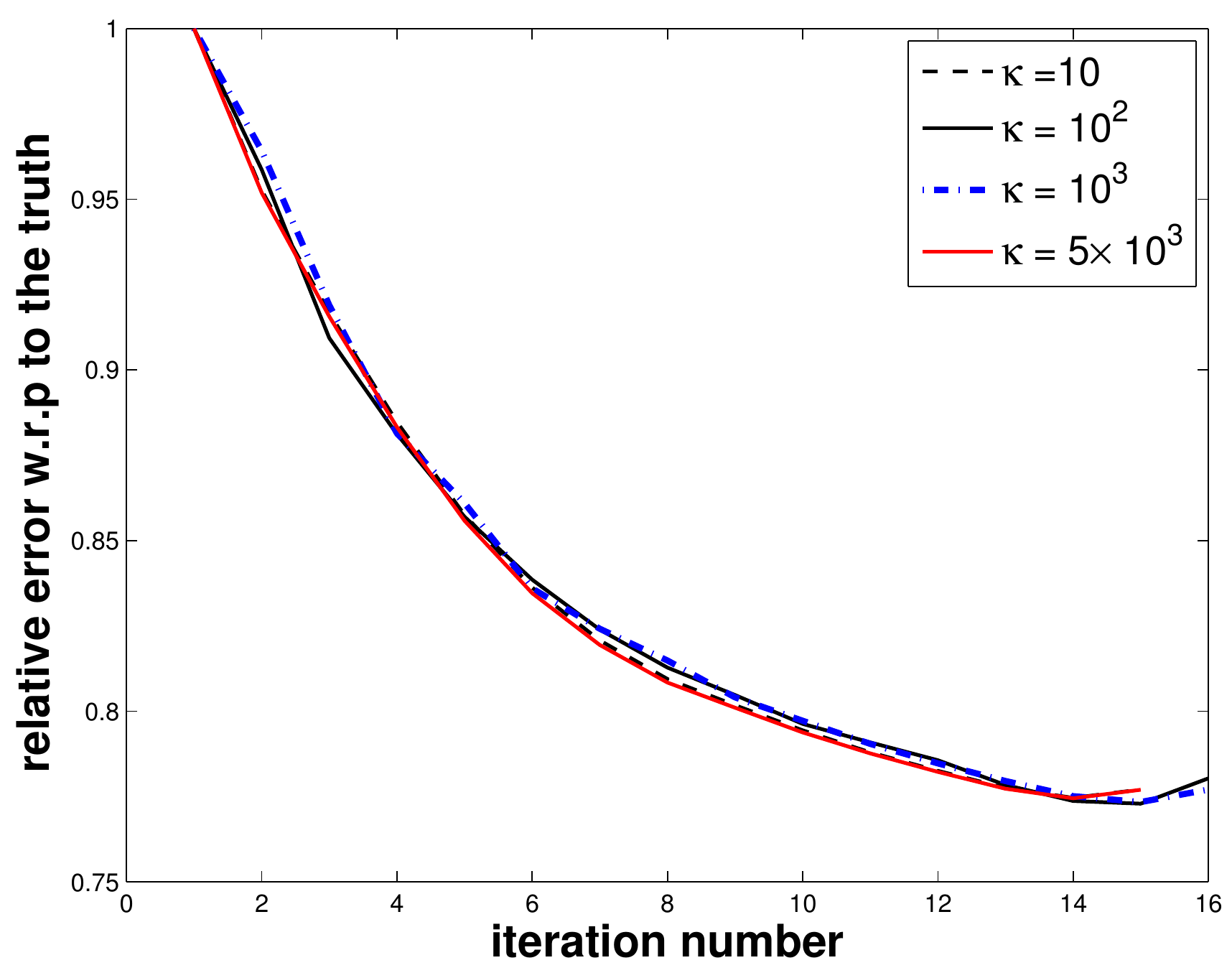}
\caption{Regularizing LM scheme for history matching. Right: data misfit. Left: relative error with respect to the truth}  
\label{Figure6}
\end{figure}

\begin{figure}
\begin{center}
\includegraphics[scale=0.25]{./True}
\includegraphics[scale=0.25]{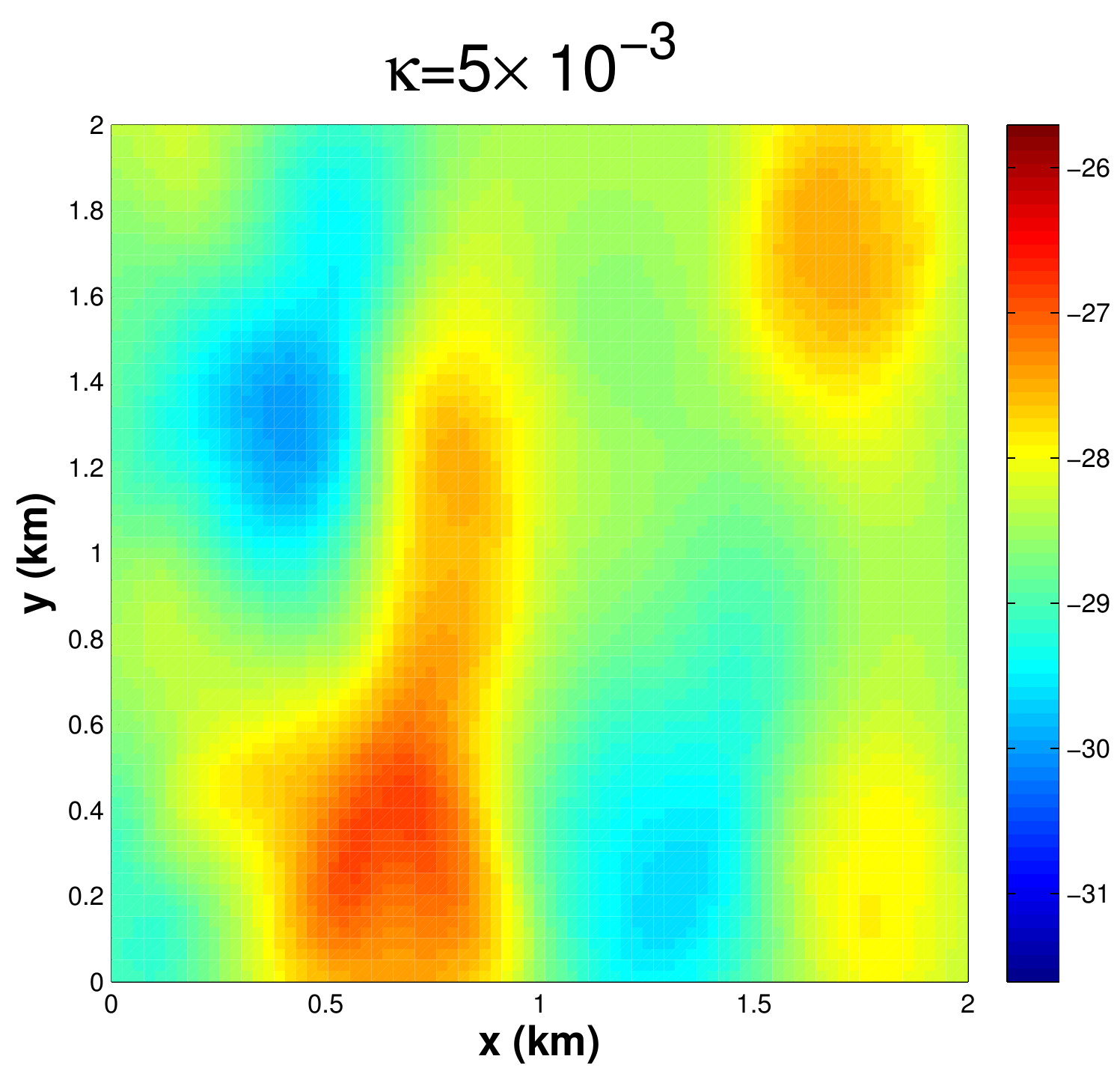}
\includegraphics[scale=0.25]{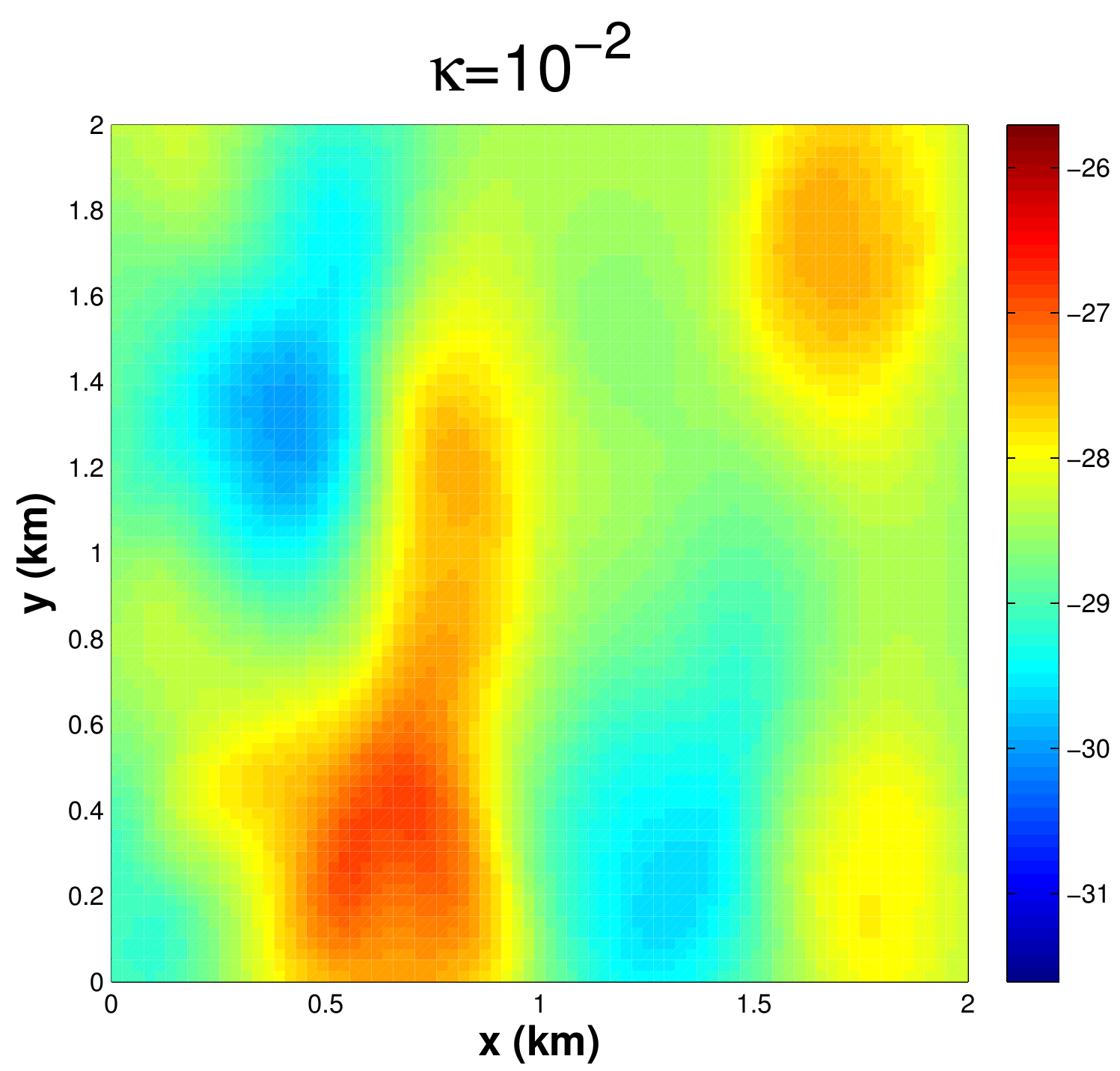}\\
\includegraphics[scale=0.25]{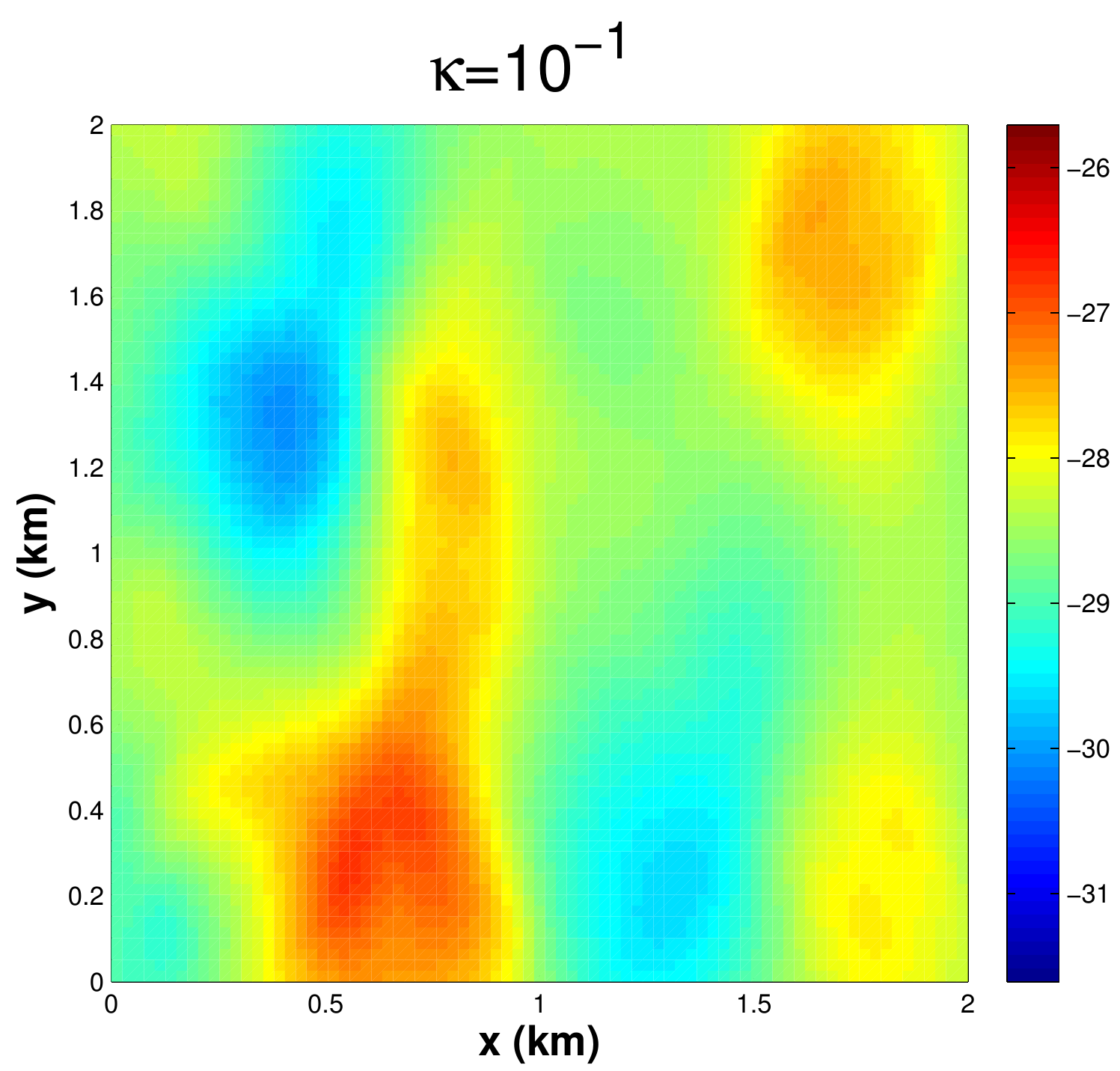}
\includegraphics[scale=0.25]{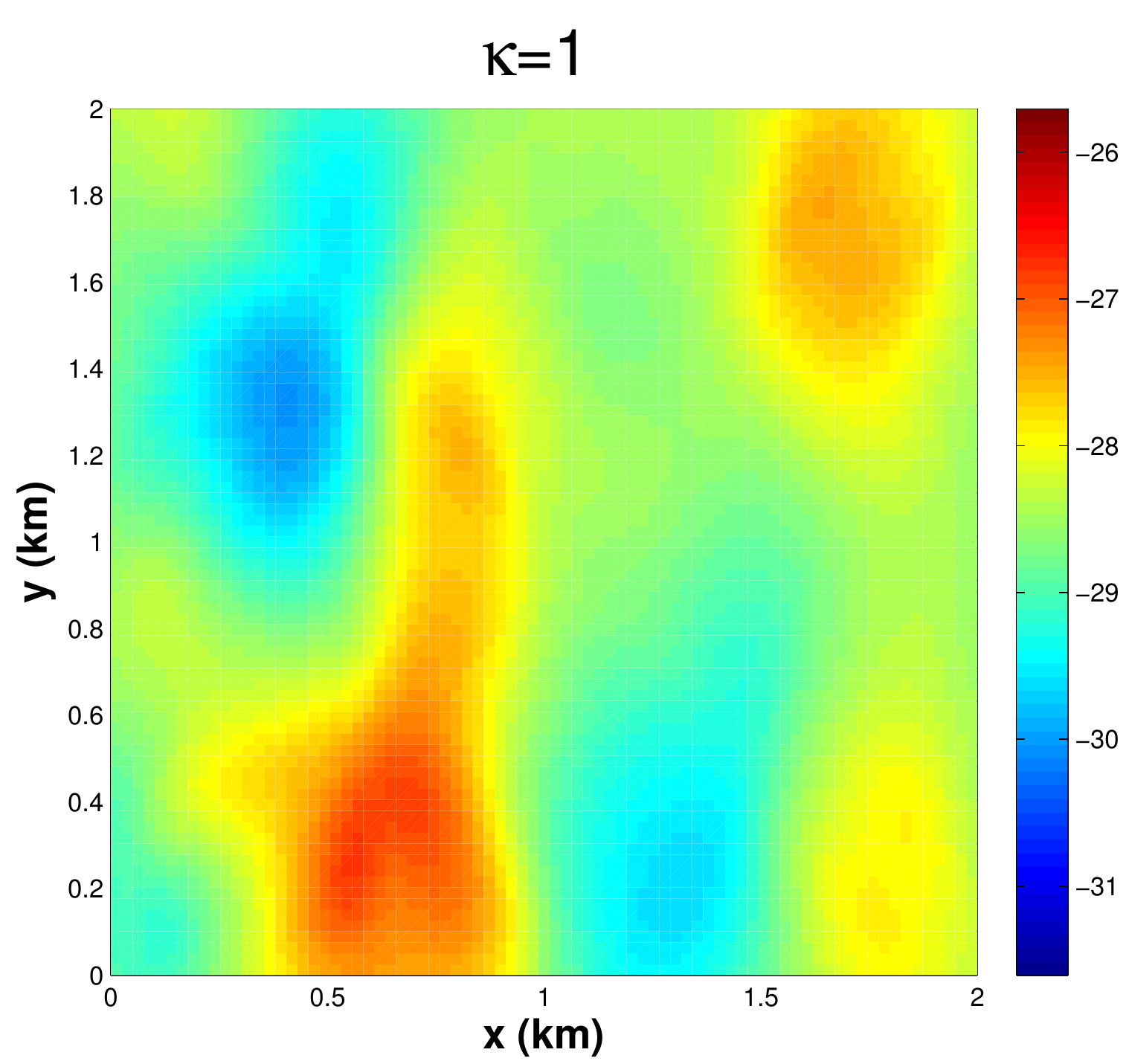}
\includegraphics[scale=0.25]{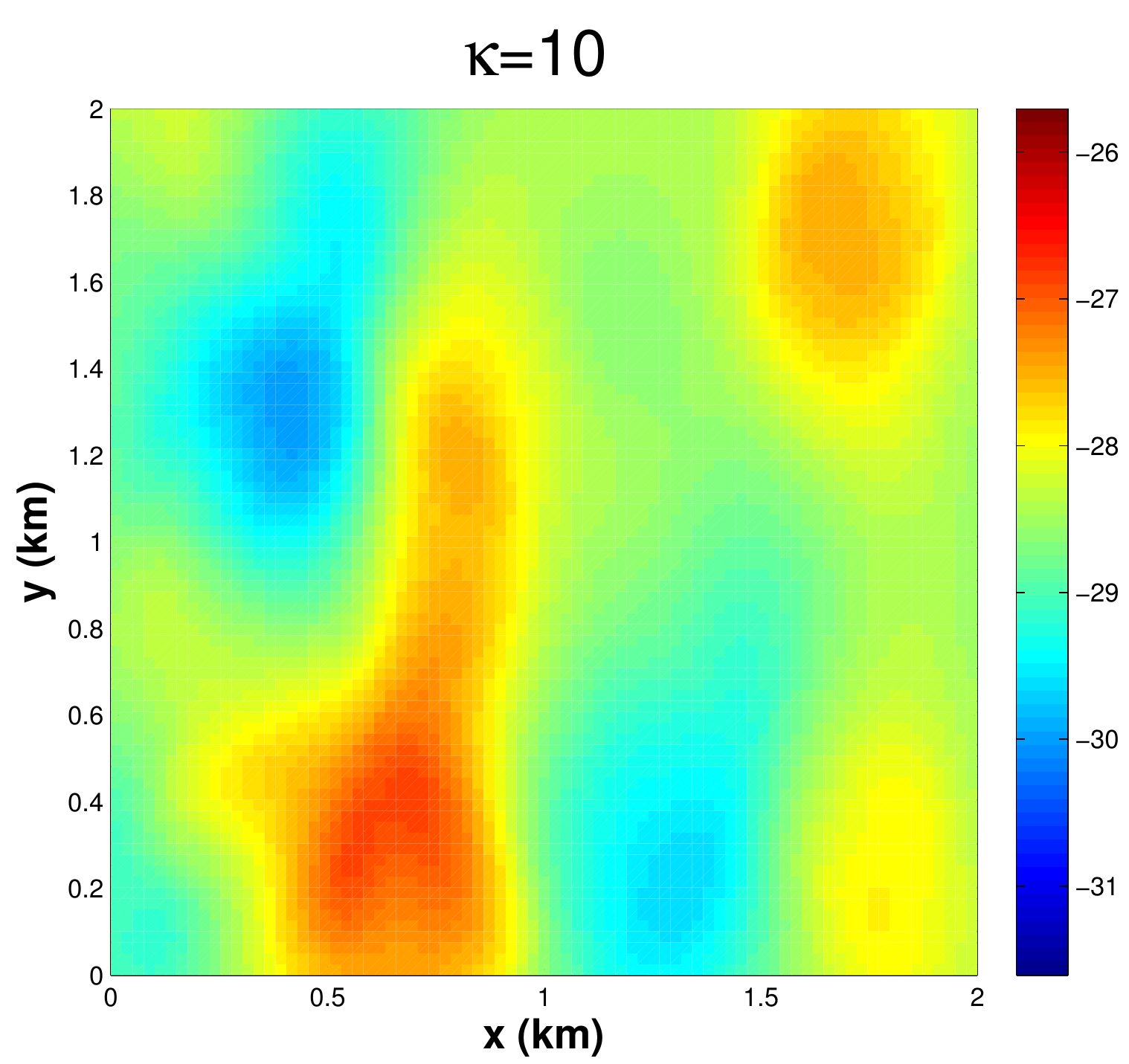}\\
\includegraphics[scale=0.25]{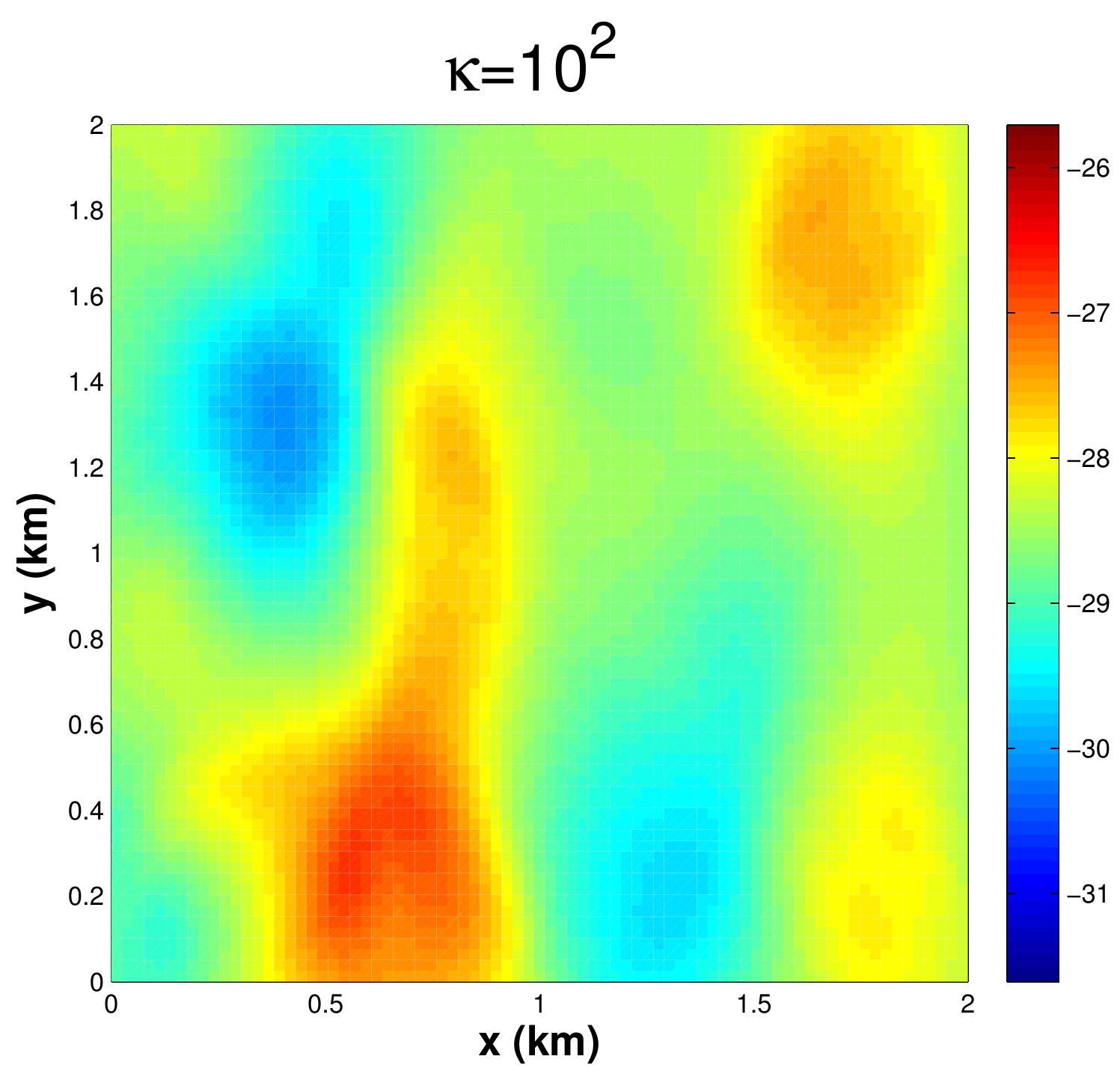}
\includegraphics[scale=0.25]{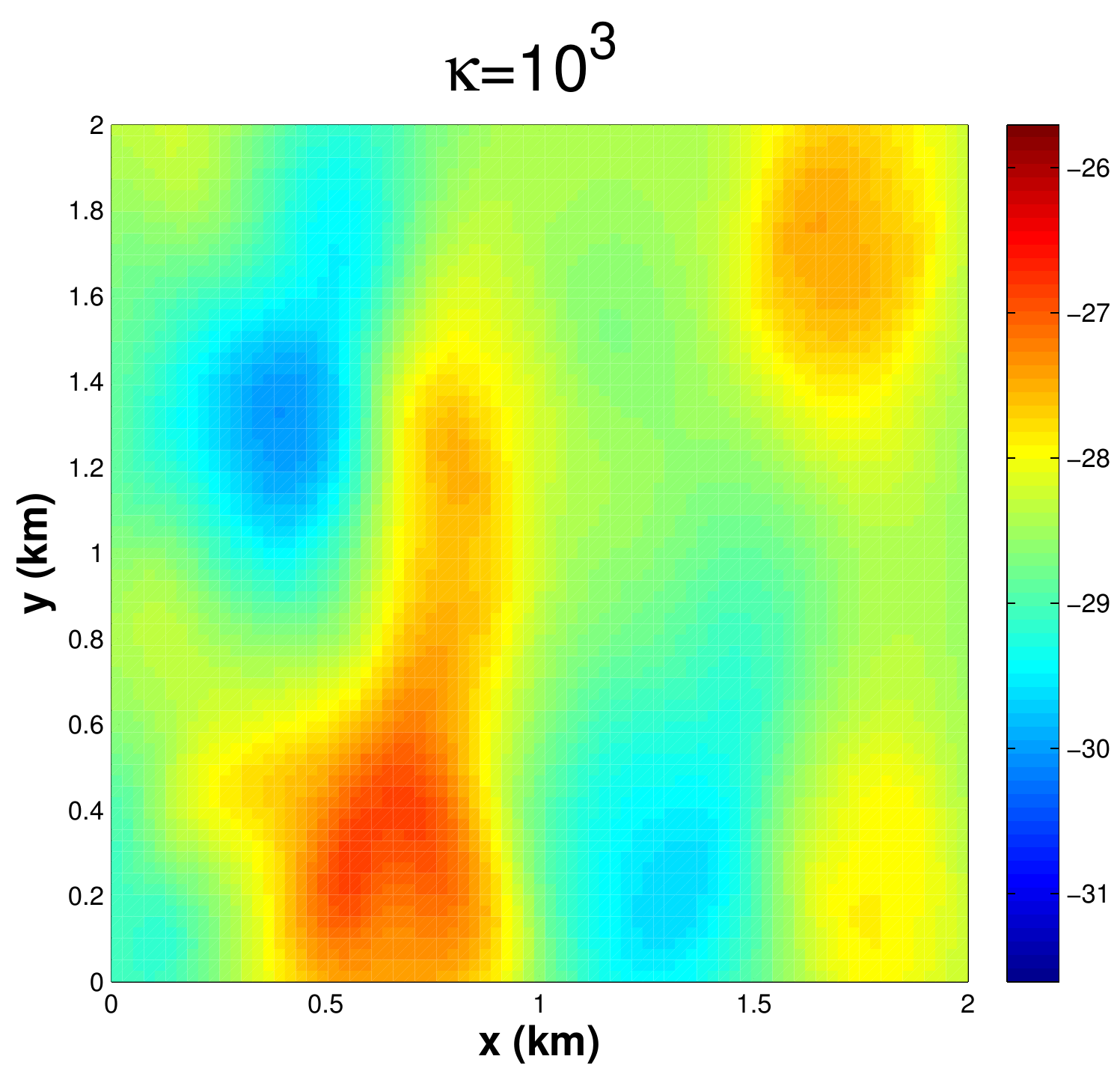}
\includegraphics[scale=0.25]{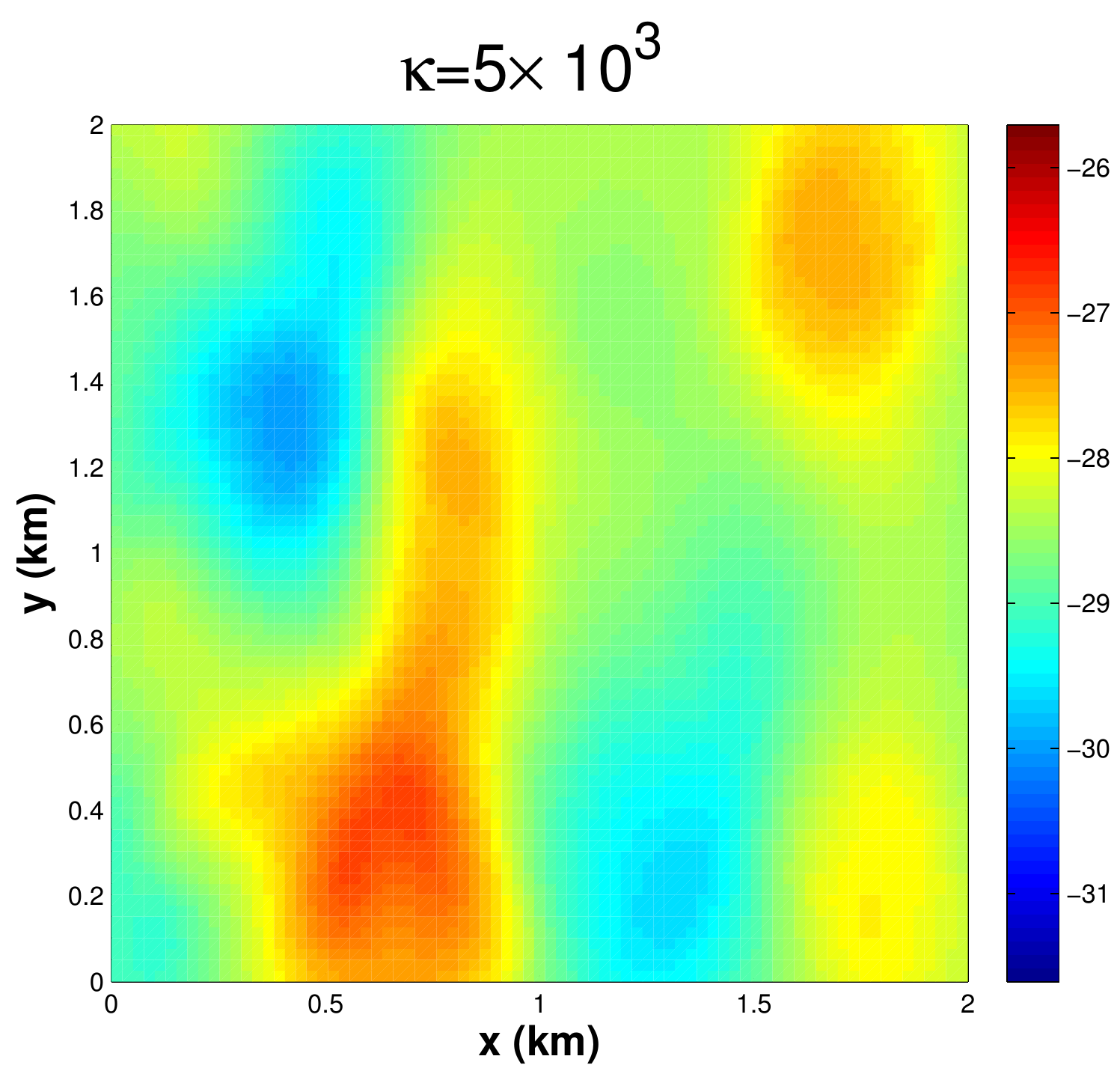}
\caption{Log-permeability estimates obtained with the regularizing LM scheme for different $\kappa$'s in (ref{eq:4.19})  [$(\log{\textrm{m}^2})$]}   \label{Figure7}
\end{center}
\end{figure}

\begin{figure}
\includegraphics[scale=0.225]{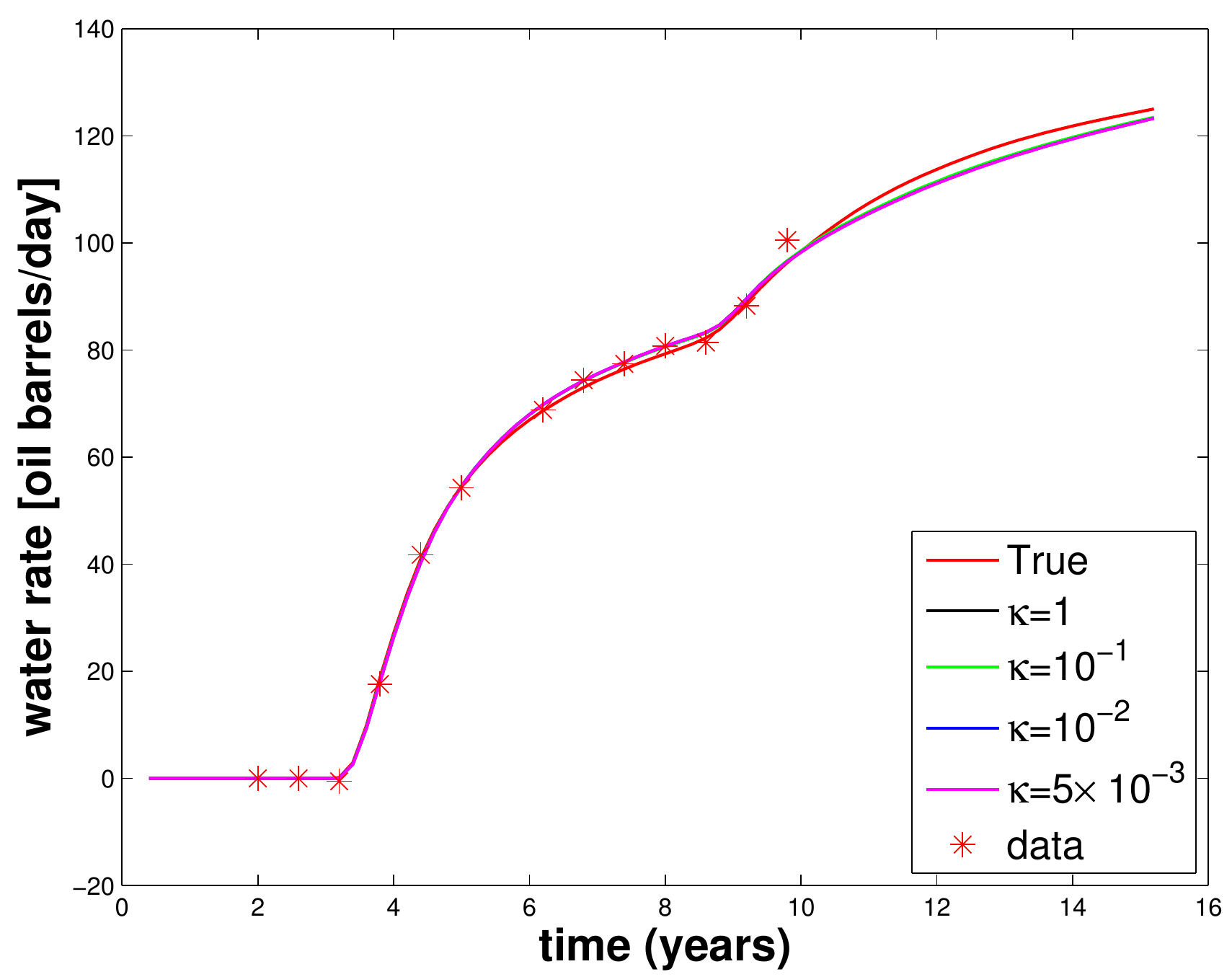}
\includegraphics[scale=0.225]{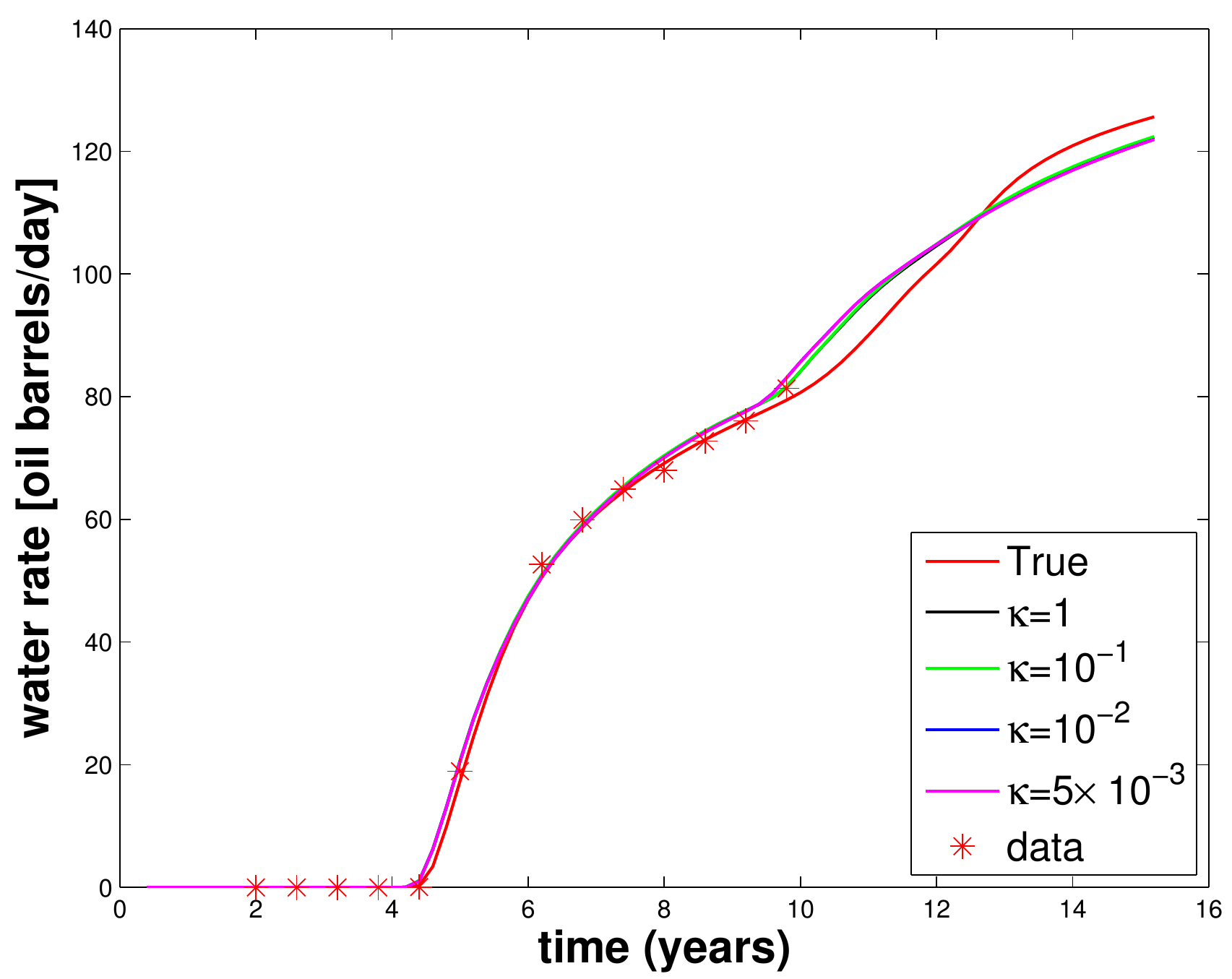}
\includegraphics[scale=0.225]{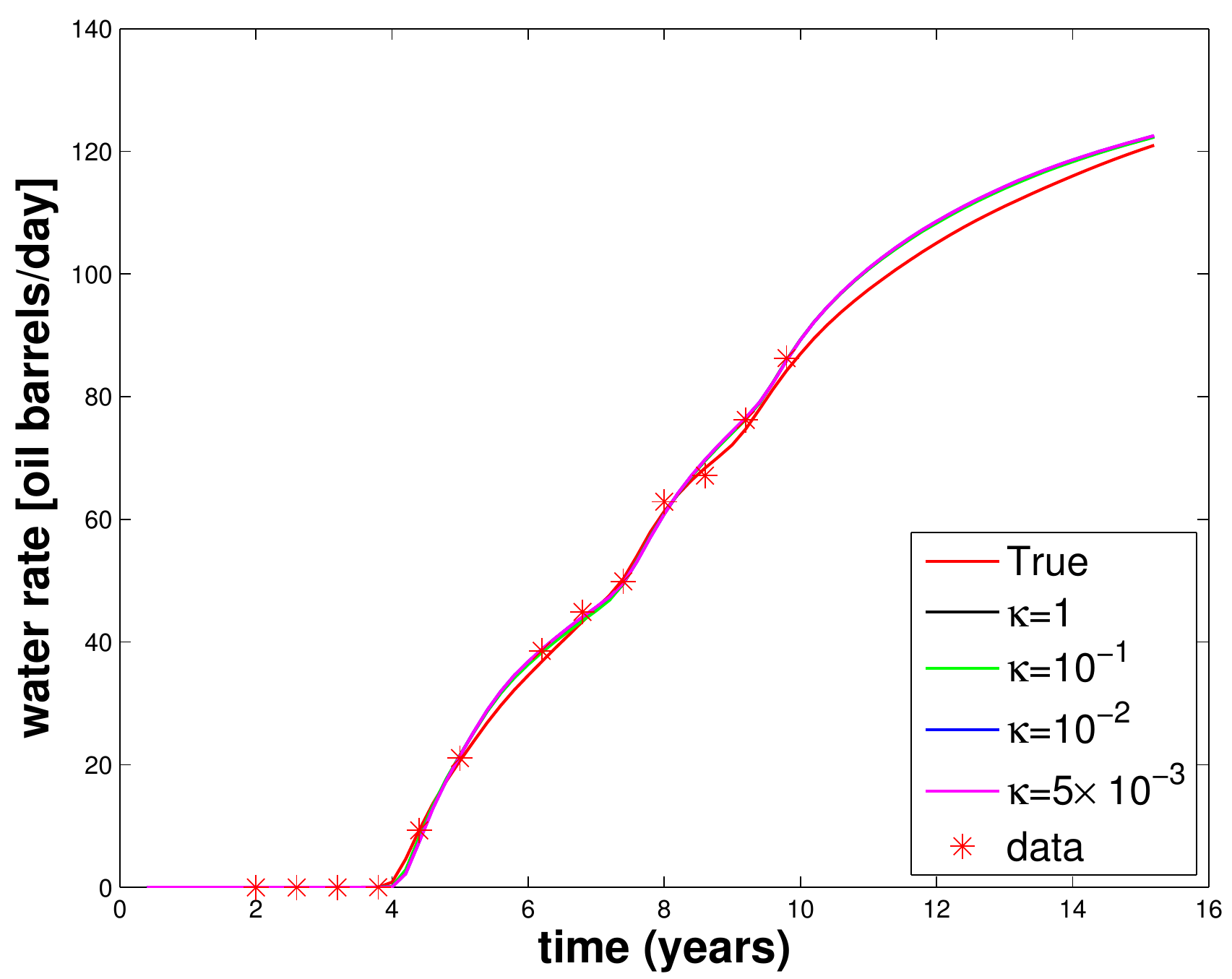}\\
\includegraphics[scale=0.225]{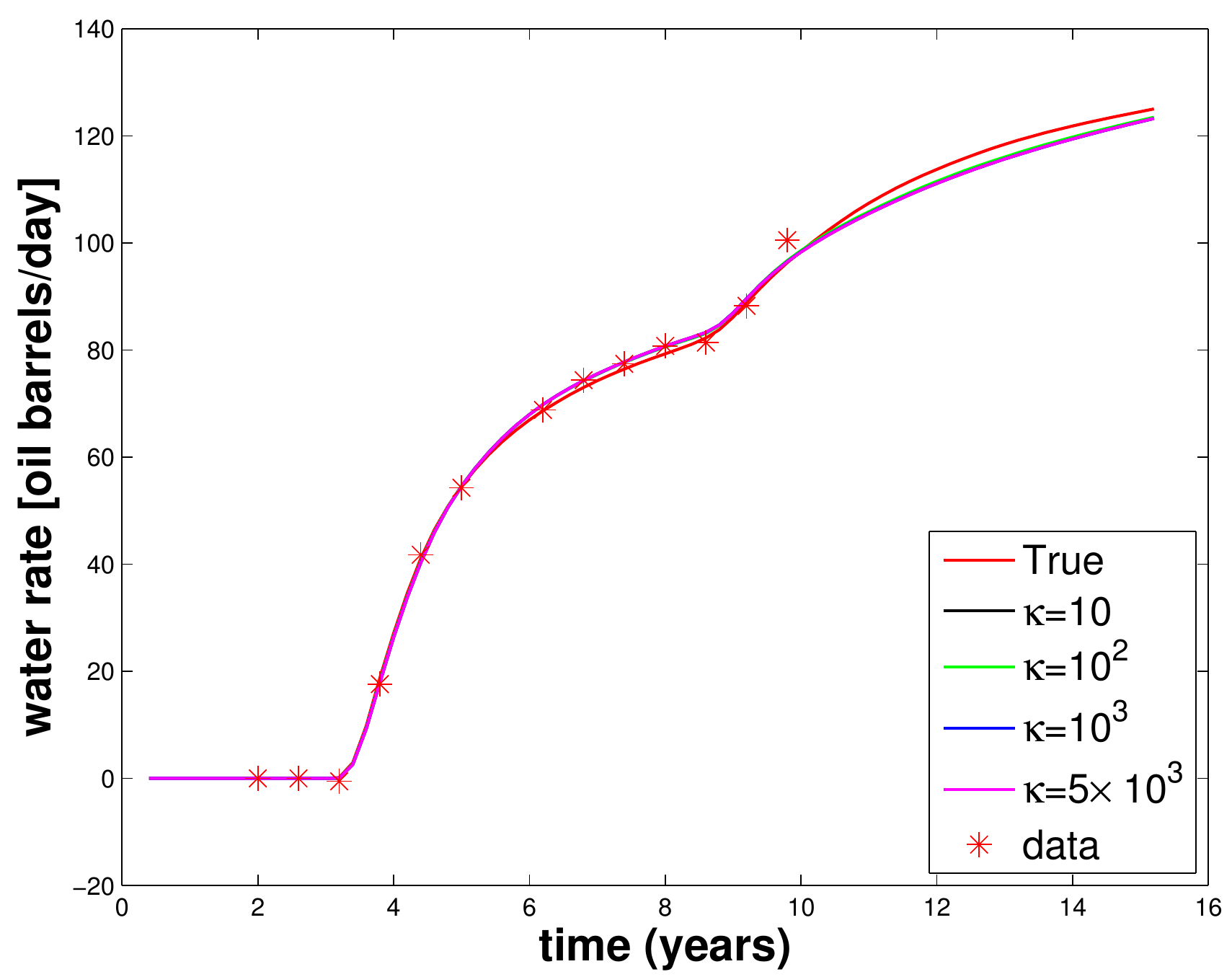}
\includegraphics[scale=0.225]{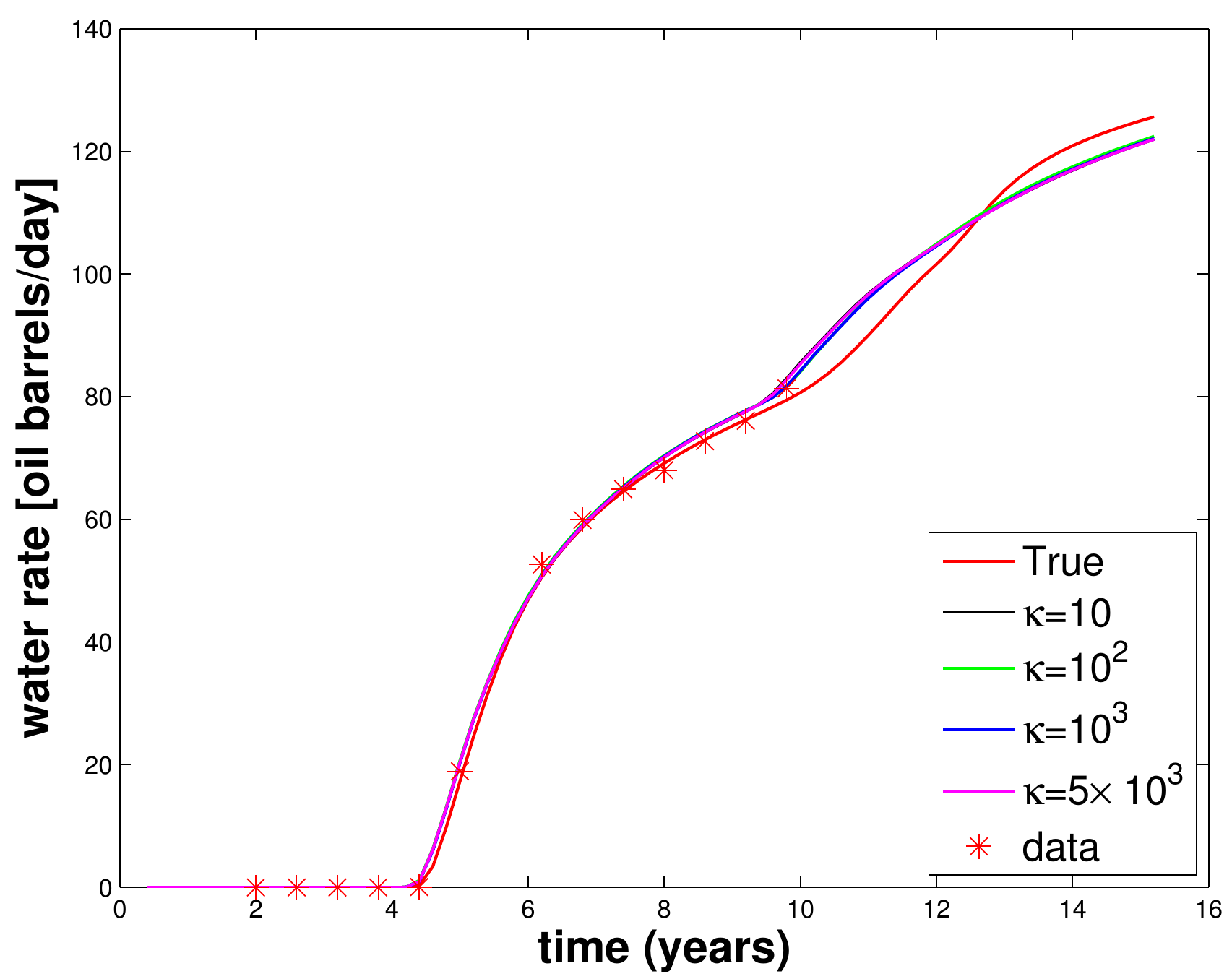}
\includegraphics[scale=0.225]{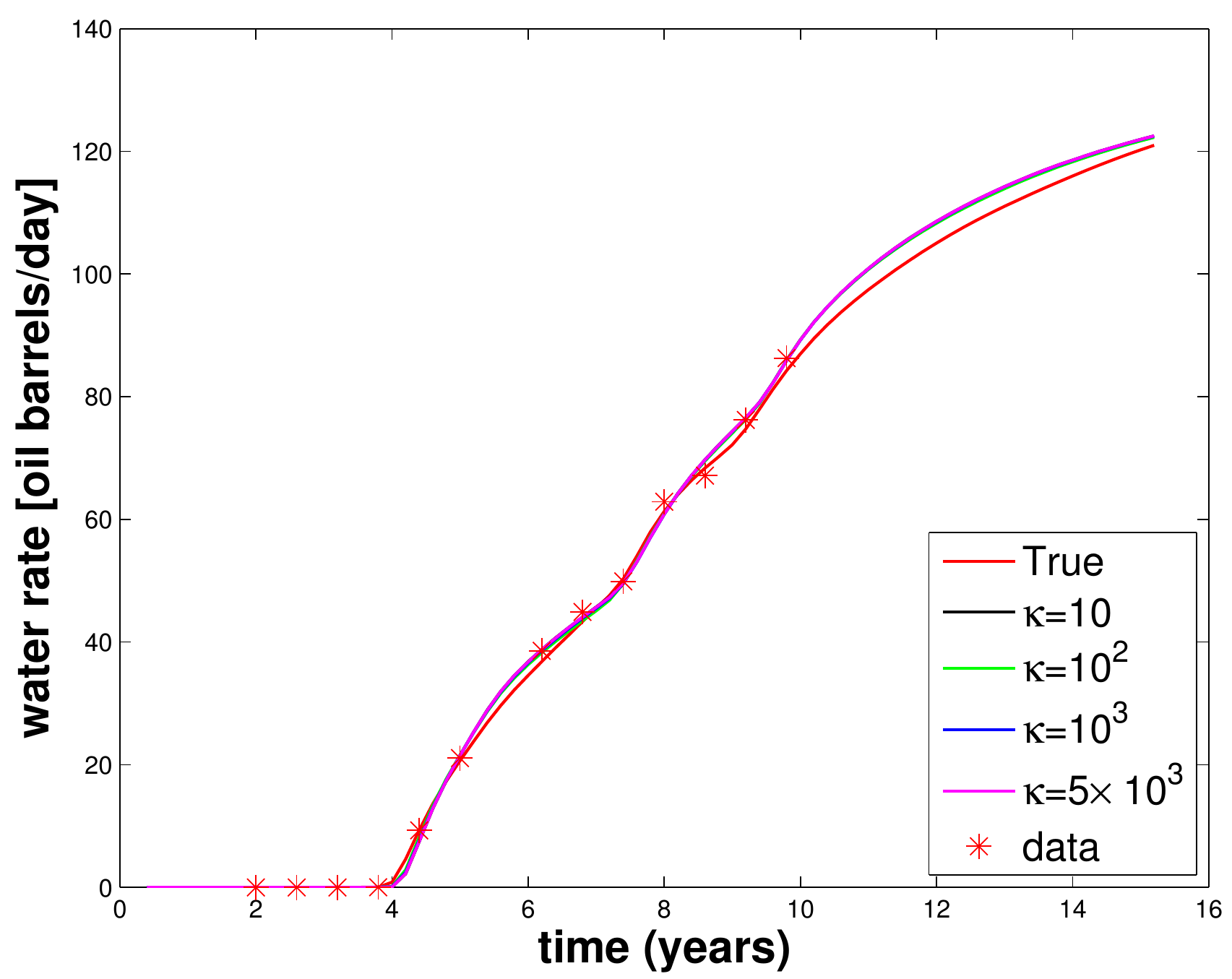}
\caption{Water rates [bbl/day]. From left to right: Wells $P_{2}$, $P_{4}$ and $P_{5}$. Top: Experiments for $\kappa\leq 1$. Bottom: Experiments for $\kappa\ge 10$ } 
\label{Figure7B}
\end{figure}

\begin{figure}
\includegraphics[scale=0.22]{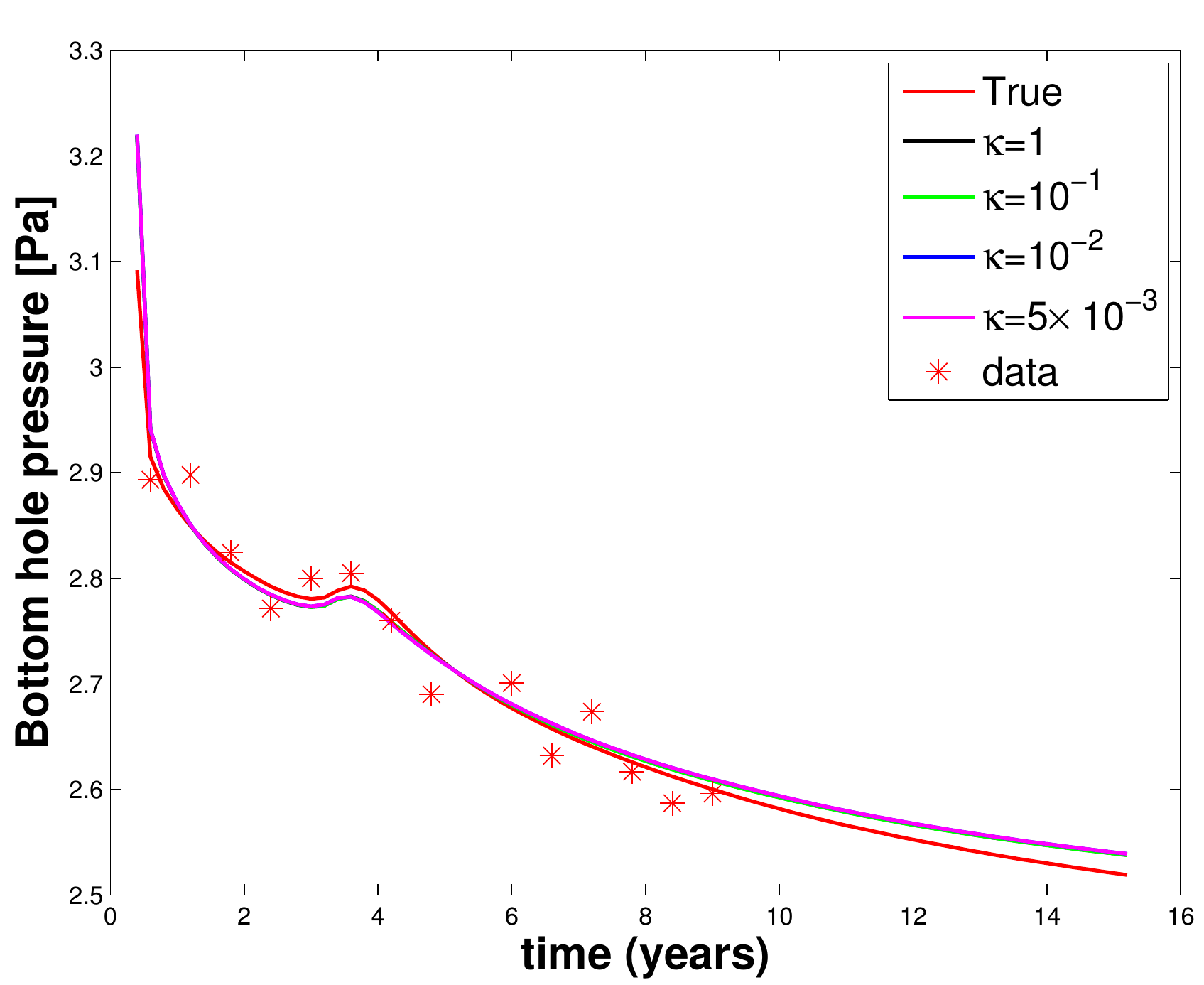}
\includegraphics[scale=0.22]{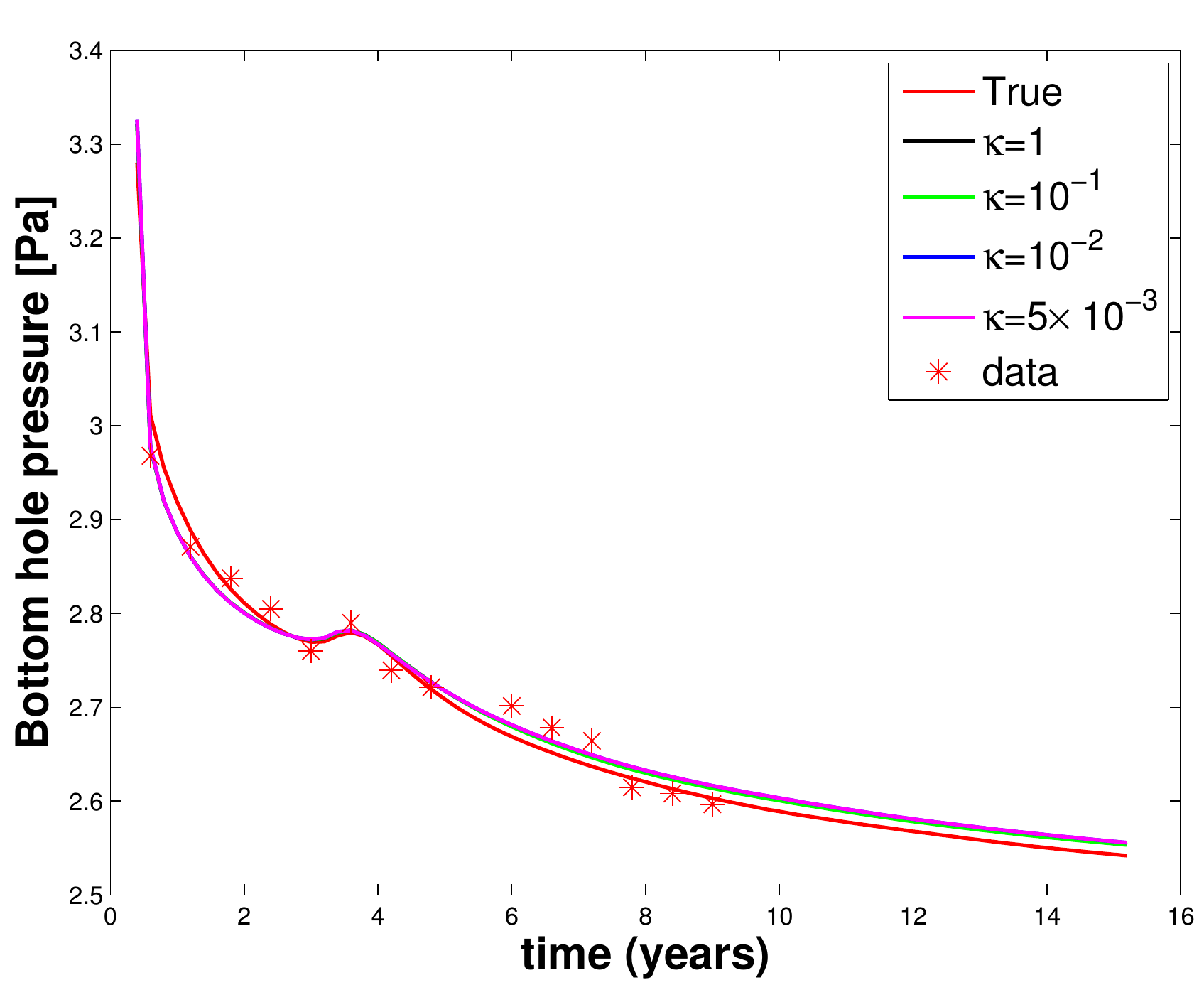}
\includegraphics[scale=0.22]{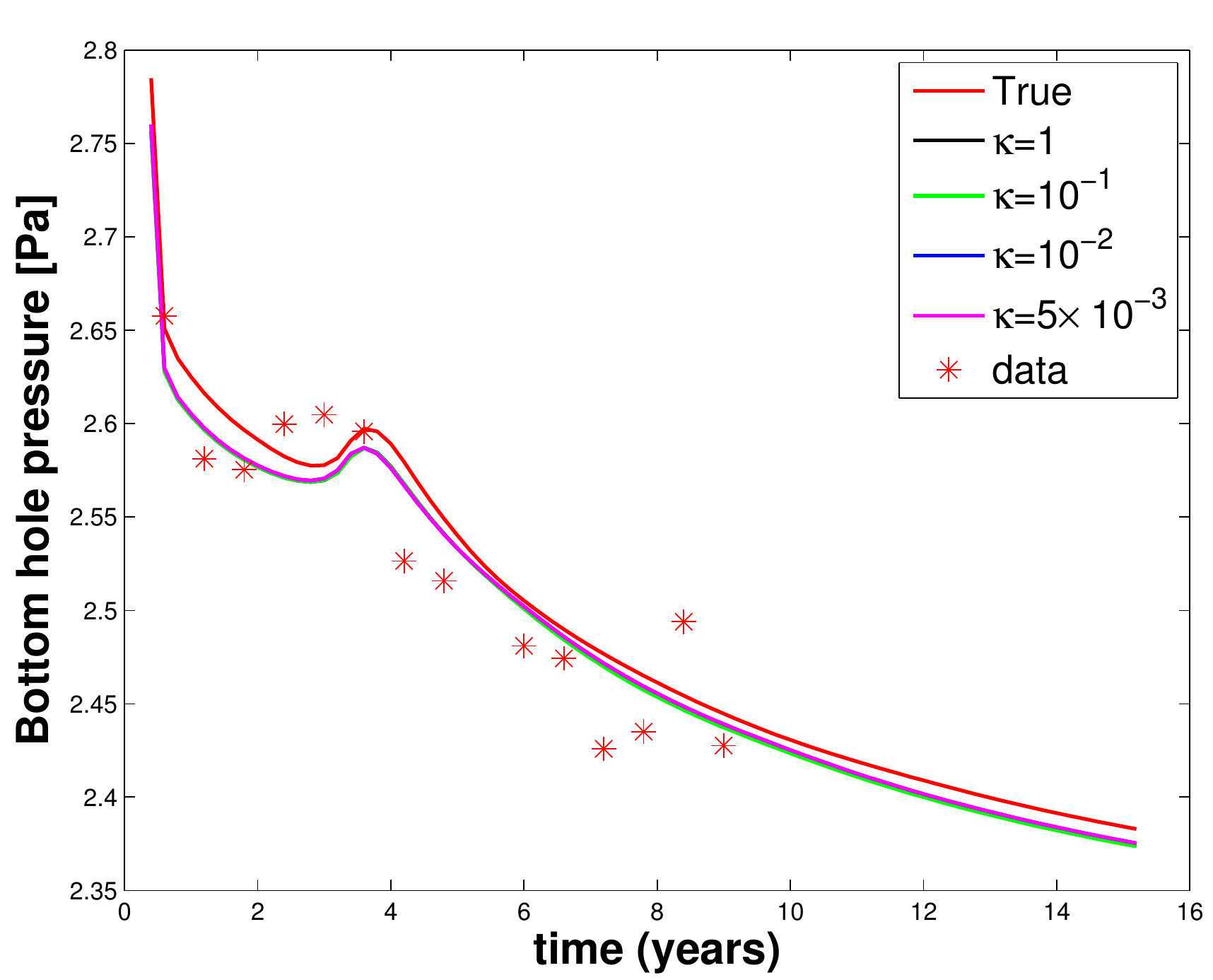}\\
\includegraphics[scale=0.22]{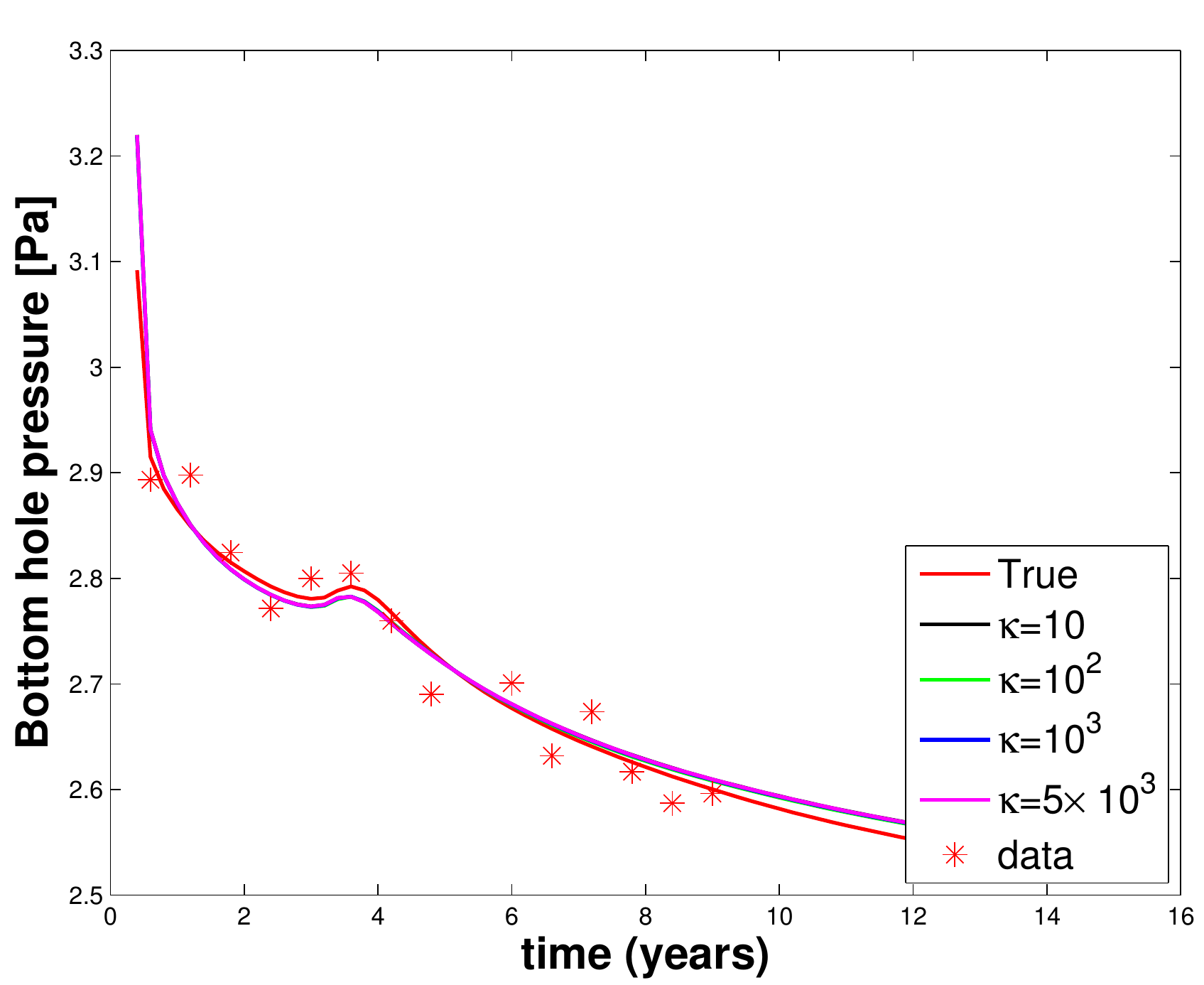}
\includegraphics[scale=0.22]{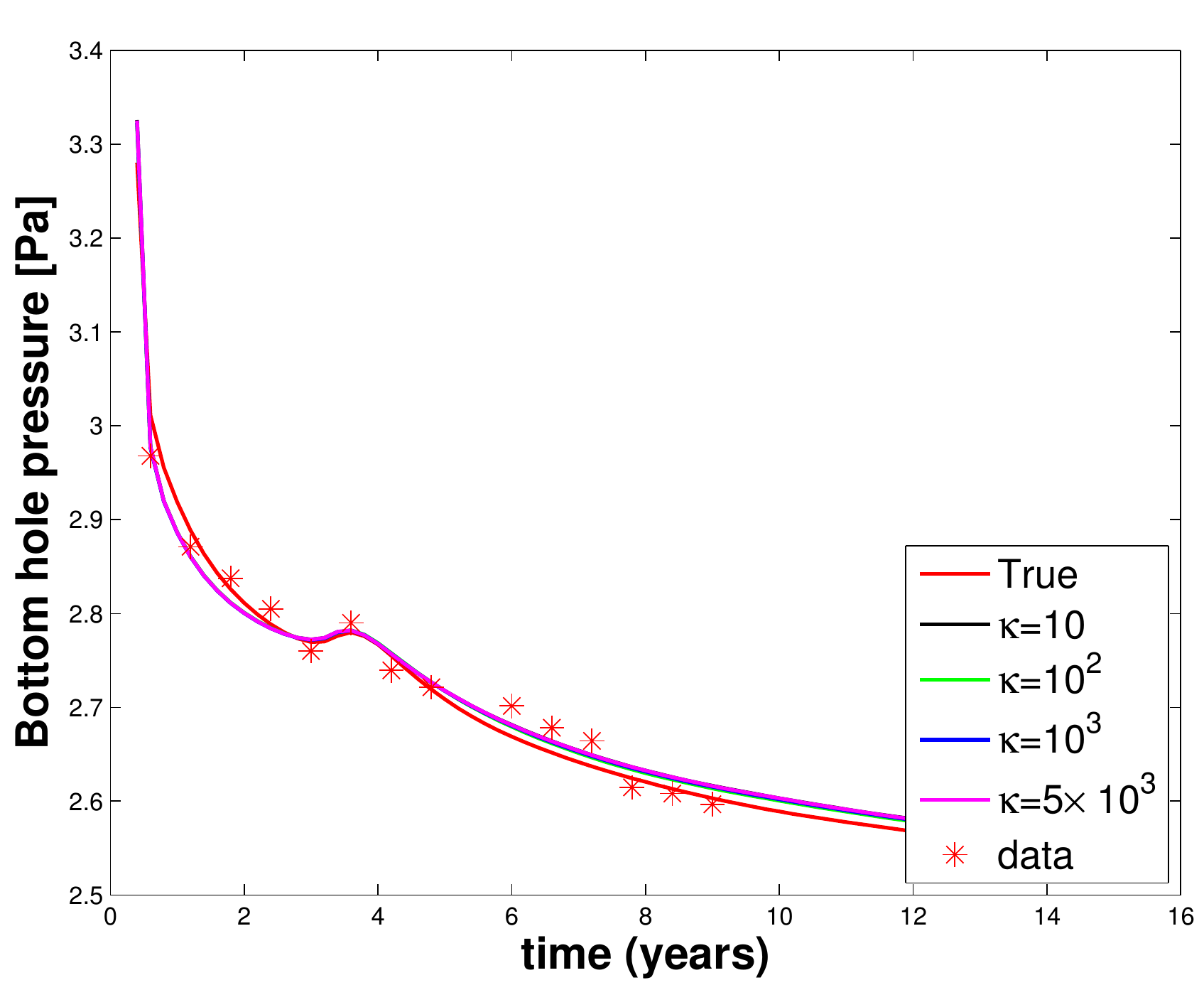}
\includegraphics[scale=0.22]{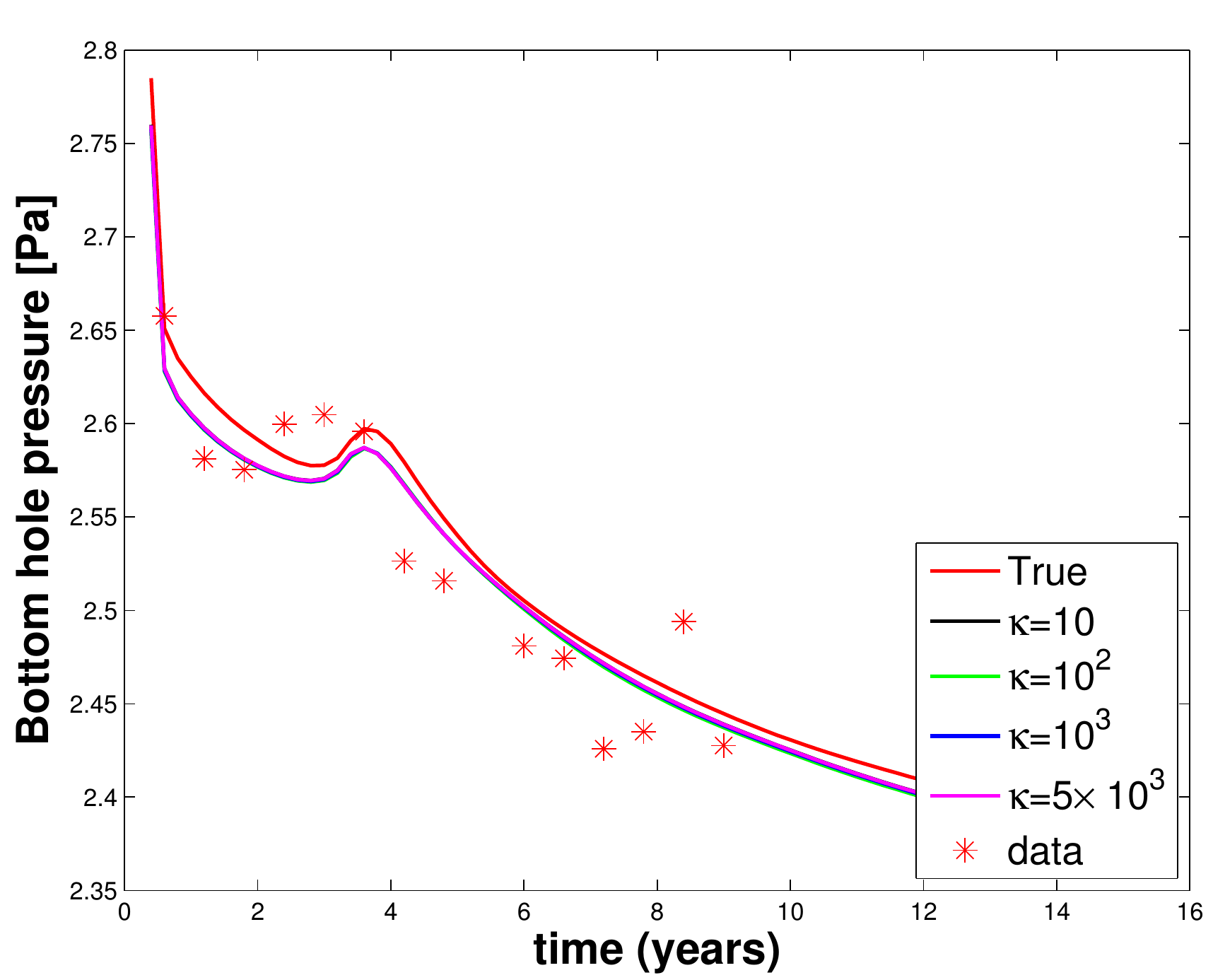}
\caption{Bottom hole pressure [Pa]. From left to right: Wells $I_{2}$, $I_{3}$ and $I_{4}$. Top: Experiments for $\kappa\leq 1$. Bottom: Experiments for $\kappa\ge 10$ } 
\label{Figure7C}
\end{figure}

In Figure \ref{Figure5} (left) we present the data misfit associated to the regularizing LM scheme for some choices of $\kappa$ with $\kappa\leq 1$ in the prior covariance (\ref{eq:3.10}). The horizontal line indicates the value of $\tau\eta$ used in the stoping criterion. We recall that the regularizing LM scheme is stopped after the estimate produces a data misfit below the aforementioned value. In Figure \ref{Figure5} (right) we display the relative error of the estimates with respect to the truth. In Figure \ref{Figure6} we present the performance of the regularizing LM scheme for $\kappa\ge 10$ in (\ref{eq:3.10}). The log-permeability estimates for all $\kappa$'s are displayed in Figure \ref{Figure7}. It is clear that the regularizing LM scheme produce similar estimates regardless the value of $\kappa$ in the covariance expression (\ref{eq:3.10}). For the present experiments, the estimates reach the stopping criterion after approximately 15 iterations. In Figure \ref{Figure7B} and Figure \ref{Figure7C} we display the model predictions obtained by simulating the water flood with the estimates of log-permeability produced with the regularizing LM scheme. As we expect from the similarities in all the estimates (Figure \ref{Figure7}) obtained for all $\kappa$'s considered here, the associated model predictions are all almost identical.

\subsection{Comparison with the standard approach}\label{compa}

\begin{figure}
\includegraphics[scale=0.35]{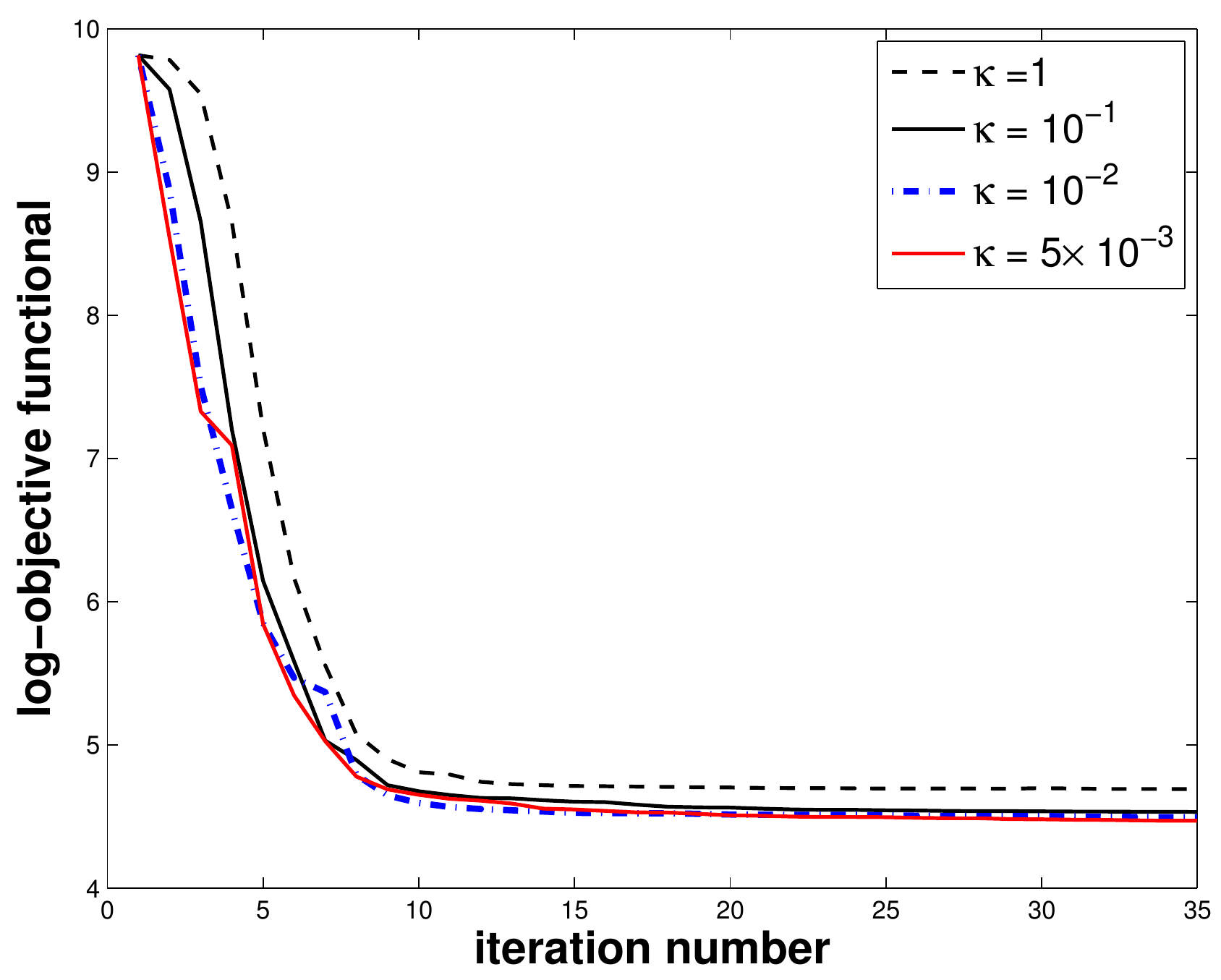}
\includegraphics[scale=0.35]{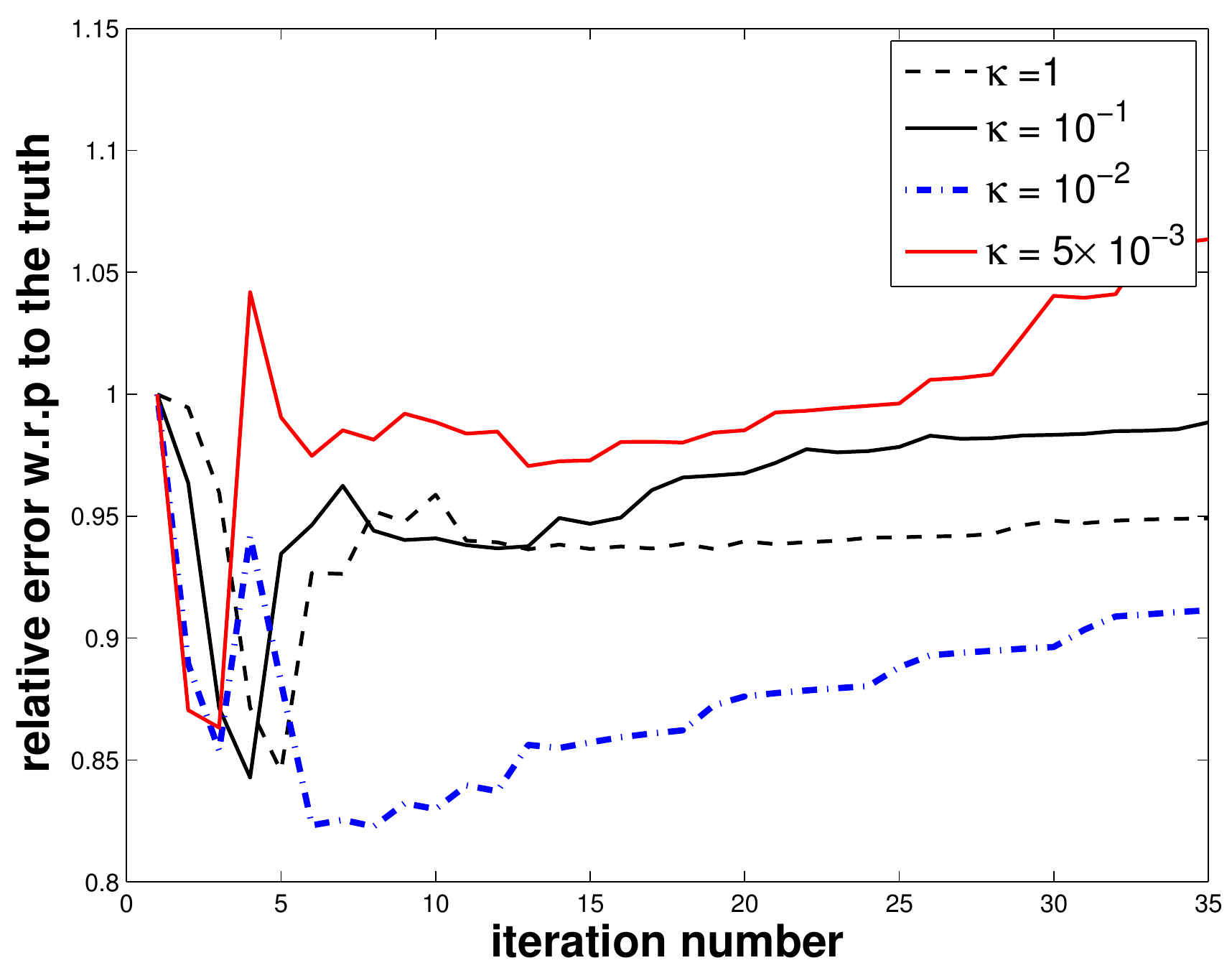}
\caption{Performance of the standard approach for history matching ($\kappa\leq 1$). Left: Objective functional (\ref{eq:sa}). Right: Relative error (left hand side (\ref{eq:3.9D}))}  

\label{Figure2}
\end{figure}

We consider the same set of synthetic data $y^{\eta}$, measurement error covariance $\Gamma$, prior mean $\overline{u}$ and $C$ (for the same $\kappa$'s) used in the experiment of subsection \ref{eq:ir1}. In this case, however, we find estimates of the log-permeability by means of the standard approach of \cite{Li,Tavakoli,svdRML,Oliver} described in subsection \ref{sa}. Note that with the choice of $C$ given by (\ref{eq:3.10})-(\ref{eq:3.11}), the objective functional that is minimized in the standard approach (\ref{eq:sa}) becomes
\begin{eqnarray}\label{eq:saB}
J(u)\equiv  \frac{1}{2}\vert\vert \Gamma^{-1/2}(y^{\eta}-G(u))\vert\vert_{Y}^{2} +\kappa\frac{1}{2}\vert\vert C_{0}^{-1/2}(u-\overline{u})\vert\vert_{X}^2 
\end{eqnarray}
 Therefore, $\kappa$ in (\ref{eq:3.10}) controls the relative size the the prior term with respect to the data misfit. In Figure \ref{Figure2} we report the performance of the experiments for $\kappa\leq 1$. The right panel of Figure \ref{Figure2} shows the relative error with respect to the truth of the estimate log-permeability field. In Figure \ref{Figure2} (left) we present the associated log-objective functional (\ref{eq:saB}). As the number of iteration increases, the method produces estimates that decreases the objective functional $J$. However, due to the lack of stability in the computation of (\ref{eq:3.9}), the error with respect to the truth increases after a certain number of iterations. Note that, even when the estimate is computed with ($\kappa=1$) the same covariance used for the generation of the truth, the corresponding error starts increasing after 5 iterations of the method. Additionally, Figure \ref{Figure2} reveals the potential failure of the stopping criteria (\ref{eq:3.9C})-(\ref{eq:3.9D}) in the standard approach of \cite{Li,Tavakoli,svdRML,Oliver}. More precisely, due to the ill-posedness of the inverse problem, a decrease of the objective functional (\ref{eq:sa}) may not be associated with a controlled change in the corresponding estimate $u$. Therefore, the choice of $\lambda$ based on (\ref{eq:3.9B}) may still lead to large values of the estimate for which (\ref{eq:3.9D}) may not be satisfied.

\begin{figure}
\includegraphics[scale=0.35]{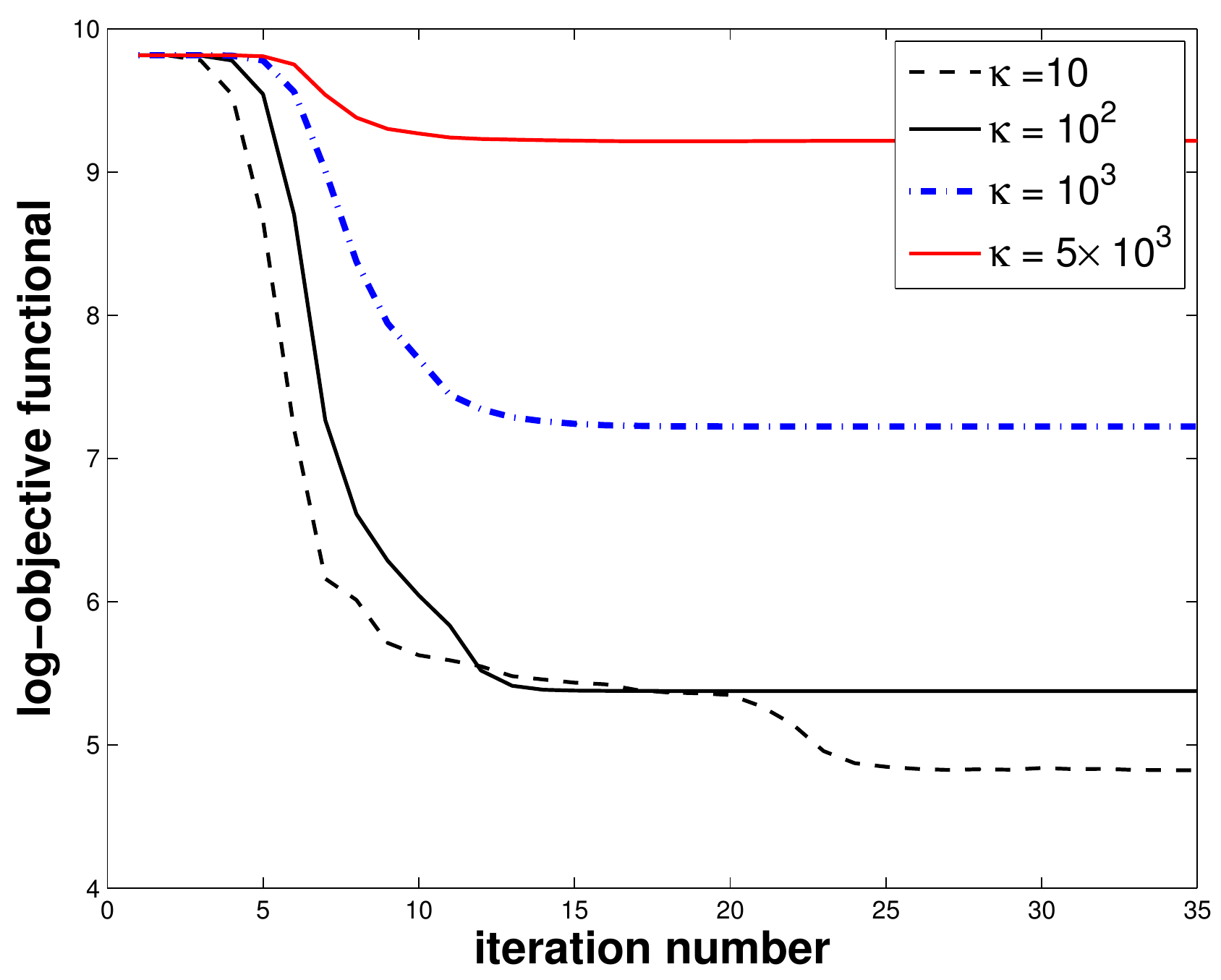}
\includegraphics[scale=0.35]{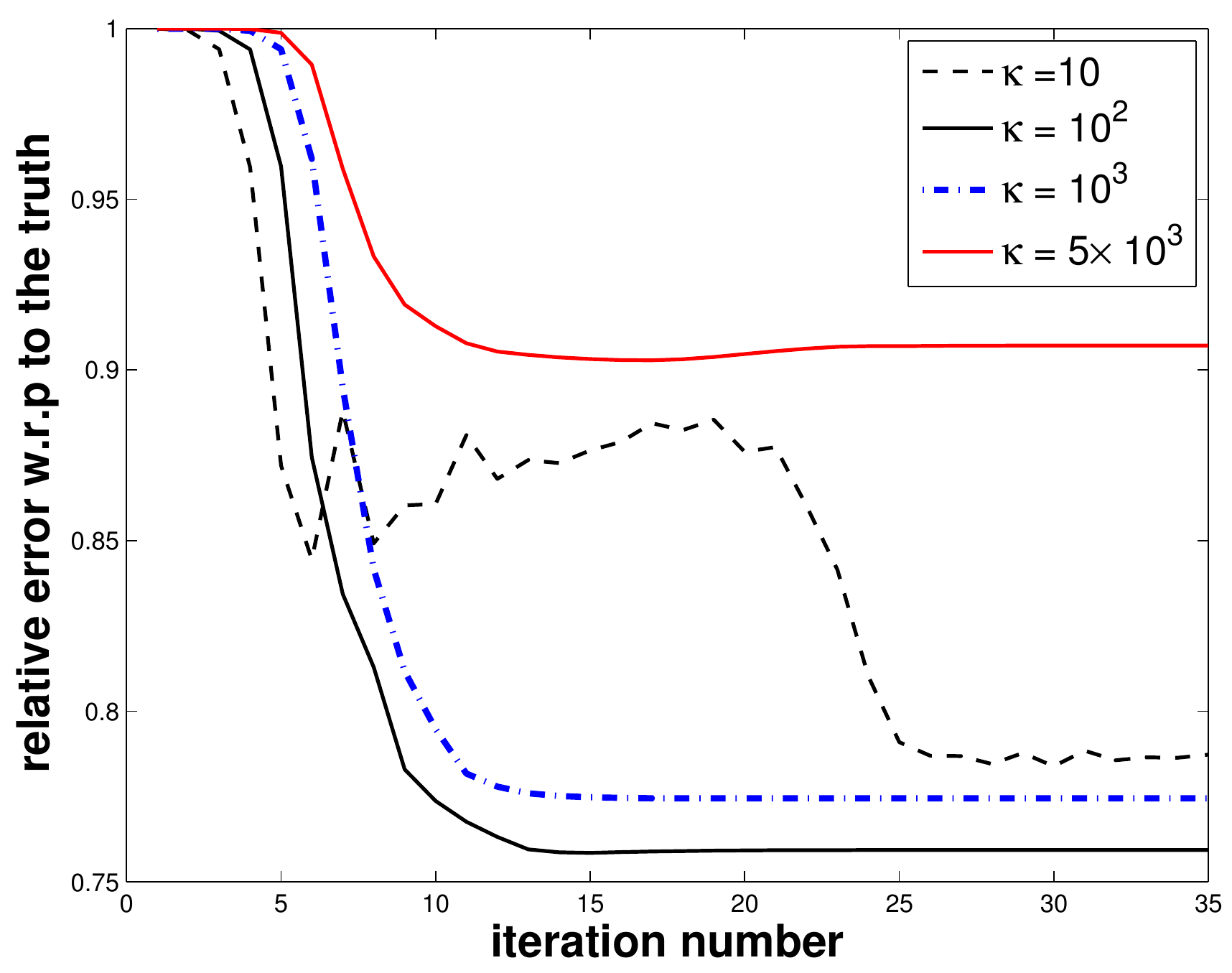}
\caption{Performance of the standard approach for history matching ($\kappa\ge 10$). Left: Objective functional (\ref{eq:sa}). Right: Relative error (left hand side (\ref{eq:3.9D})) }    \label{Figure3}
\end{figure}

\begin{figure}
\begin{center}
\includegraphics[scale=0.25]{./True}
\includegraphics[scale=0.25]{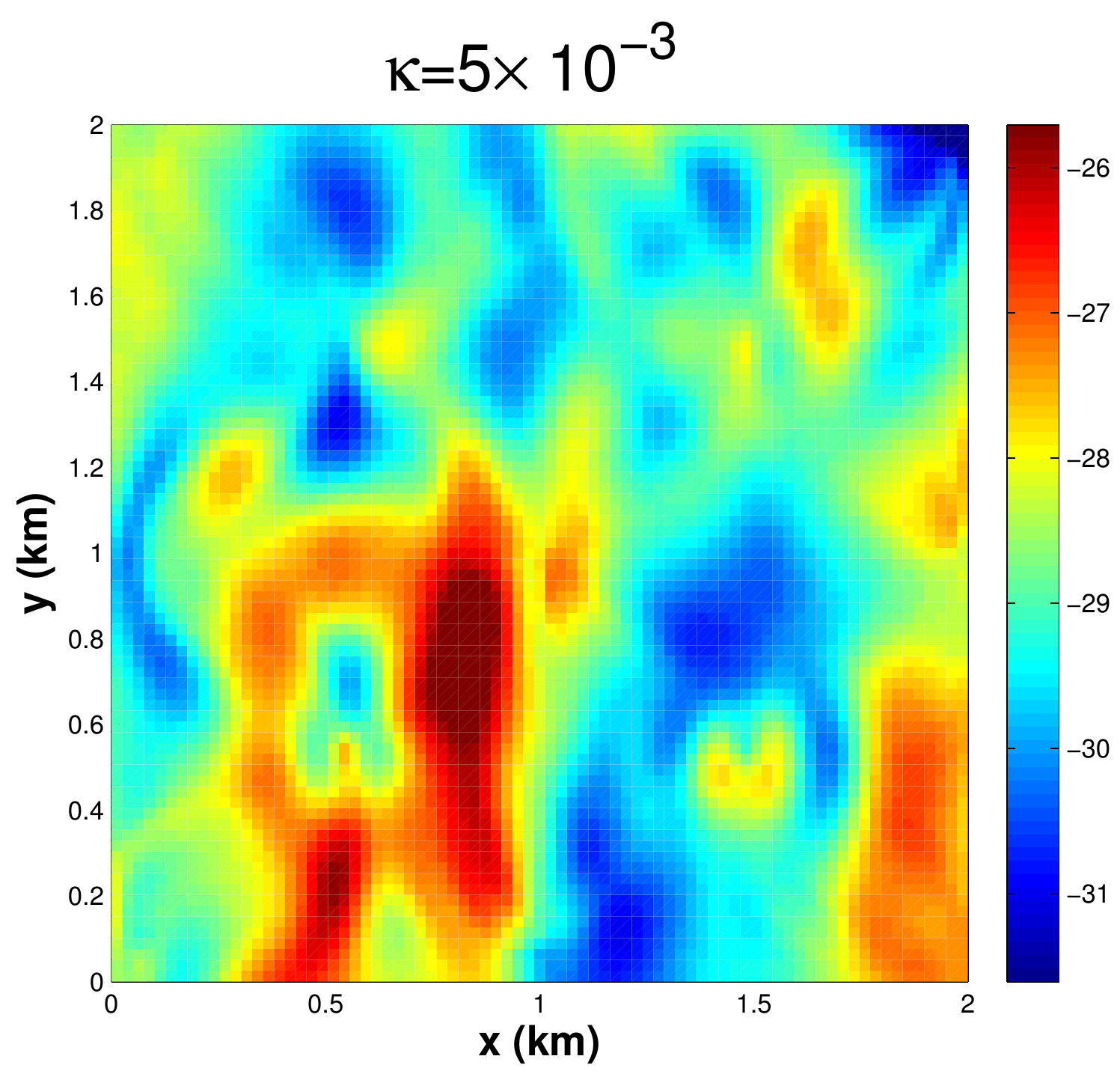}
\includegraphics[scale=0.25]{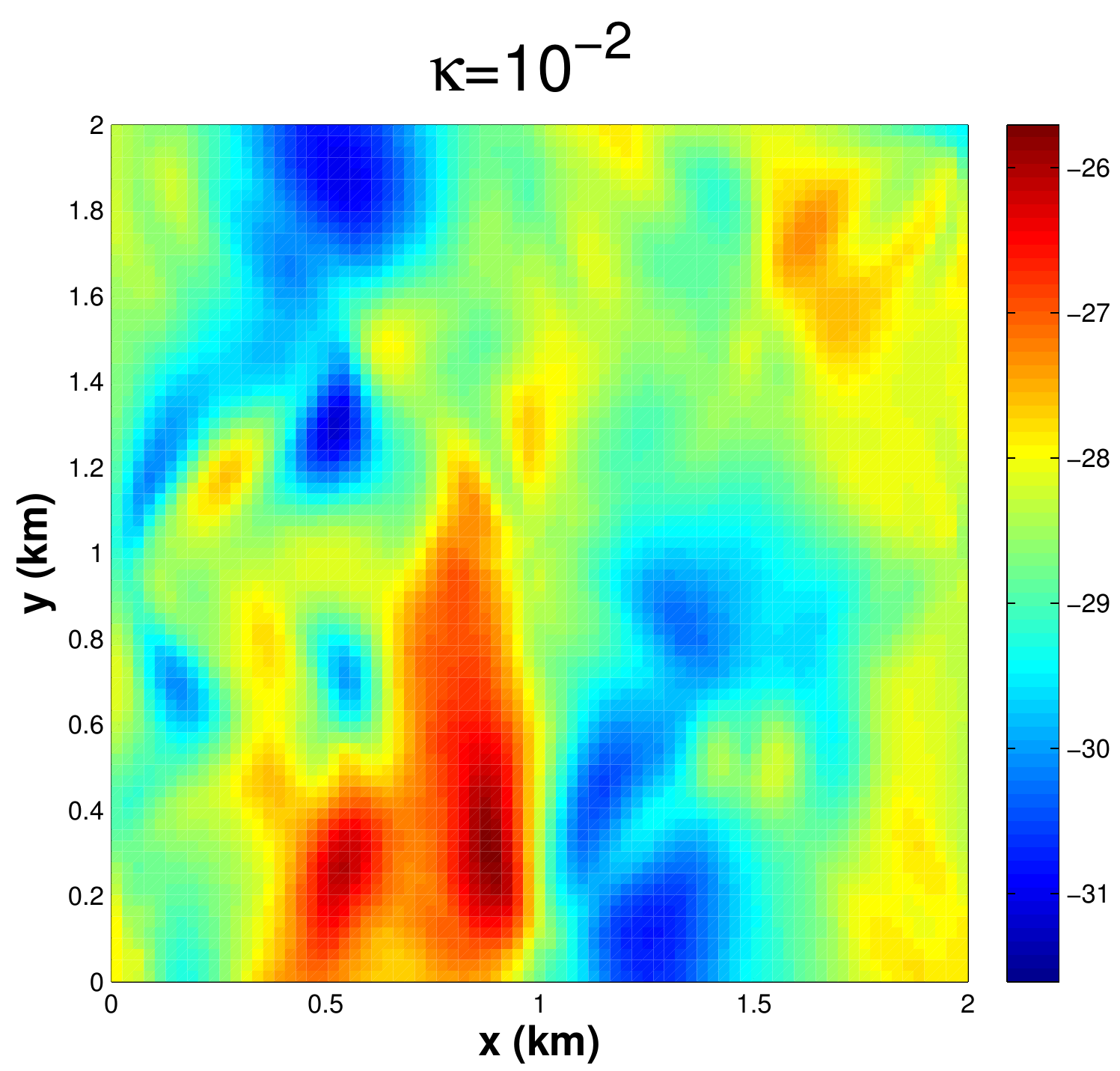}\\
\includegraphics[scale=0.25]{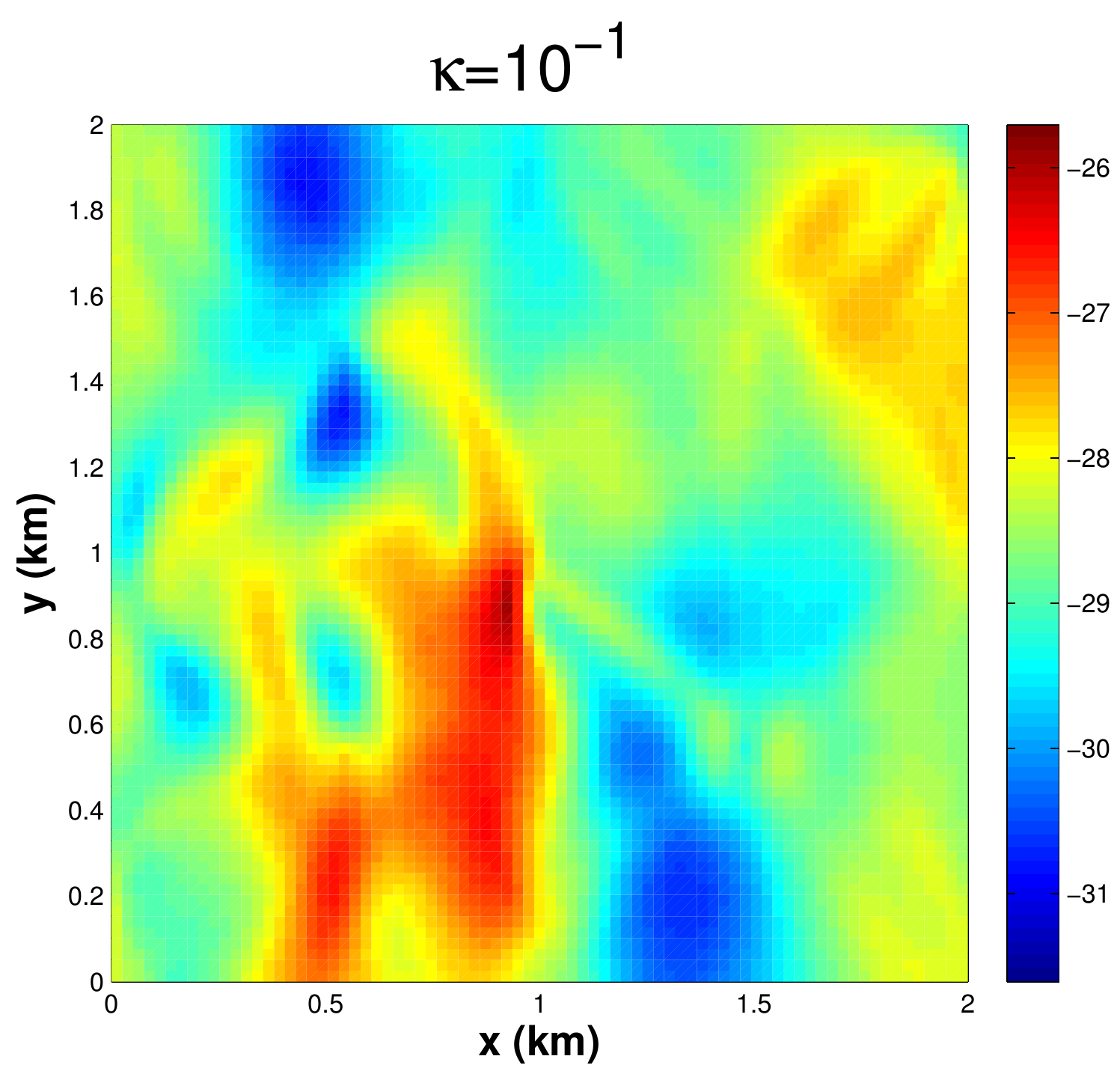}
\includegraphics[scale=0.25]{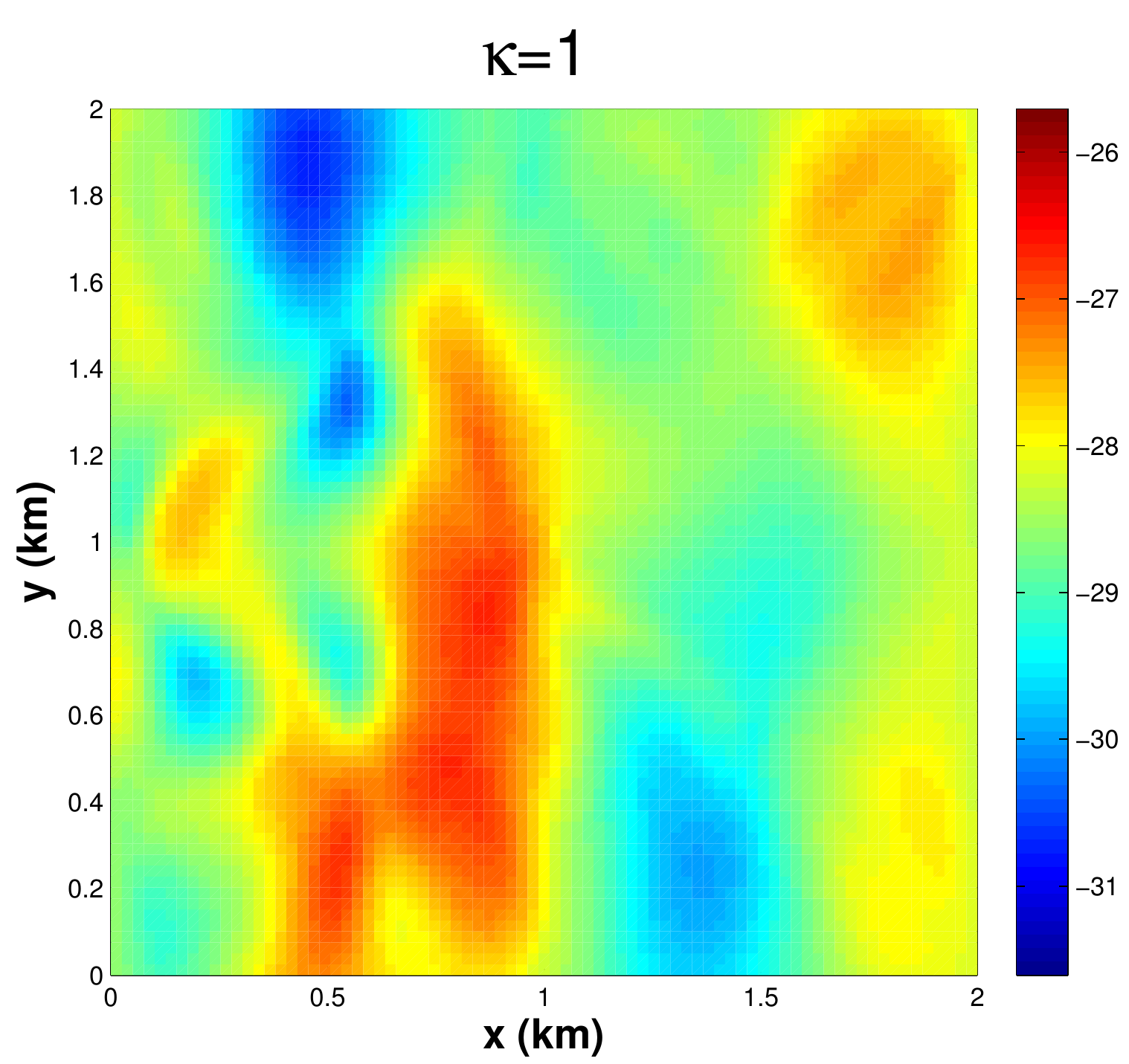}
\includegraphics[scale=0.25]{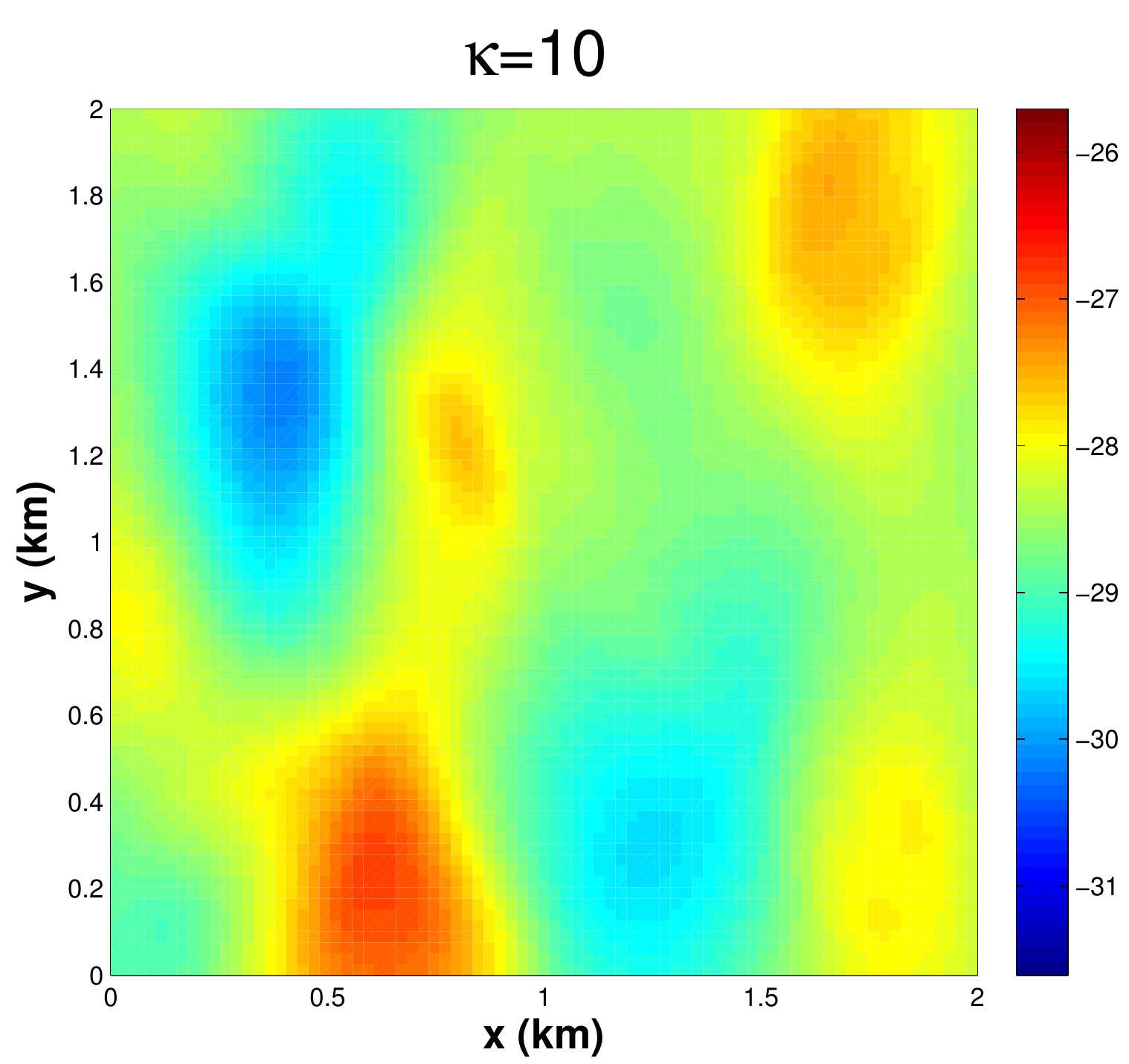}\\
\includegraphics[scale=0.25]{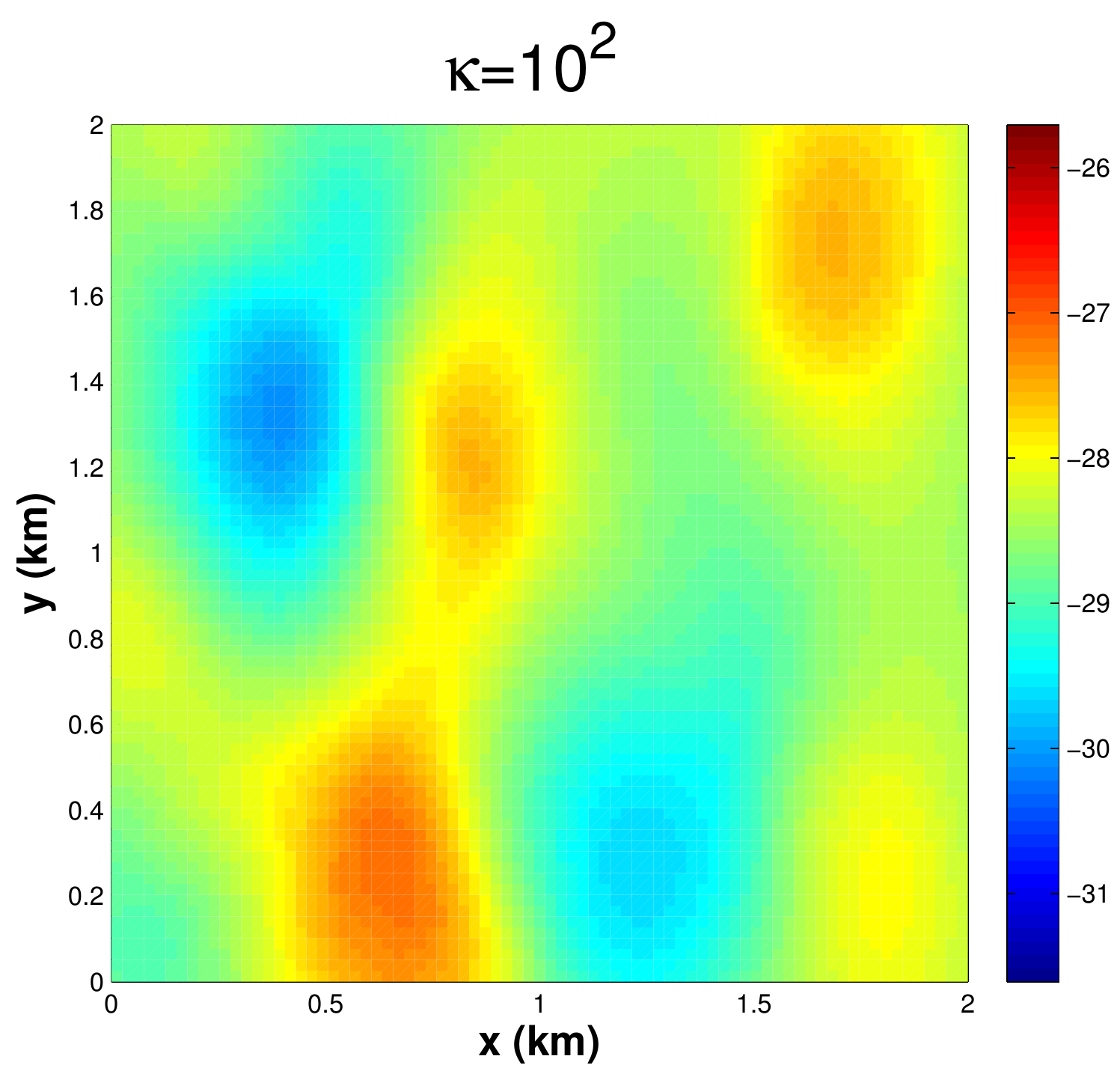}
\includegraphics[scale=0.25]{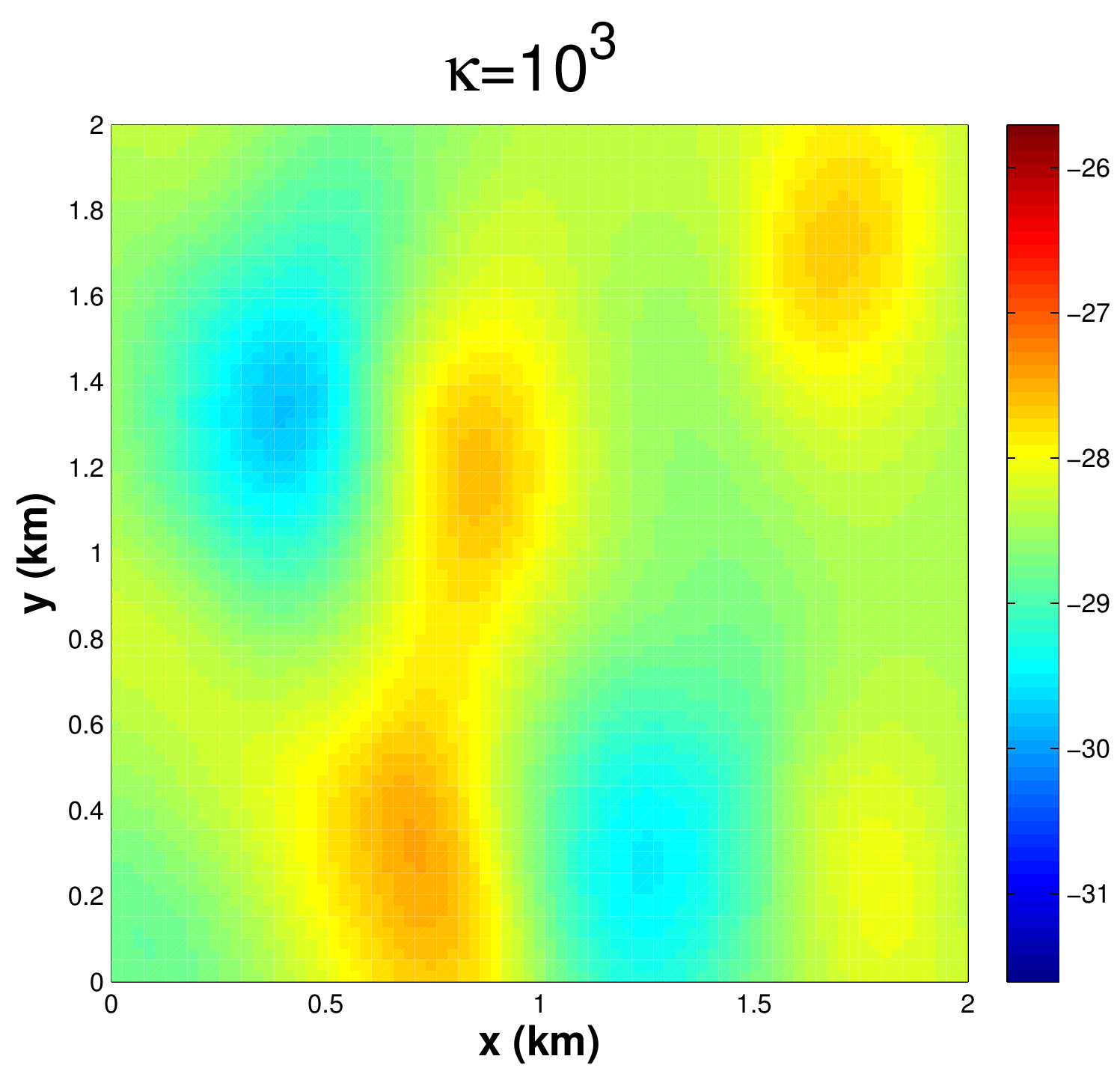}
\includegraphics[scale=0.25]{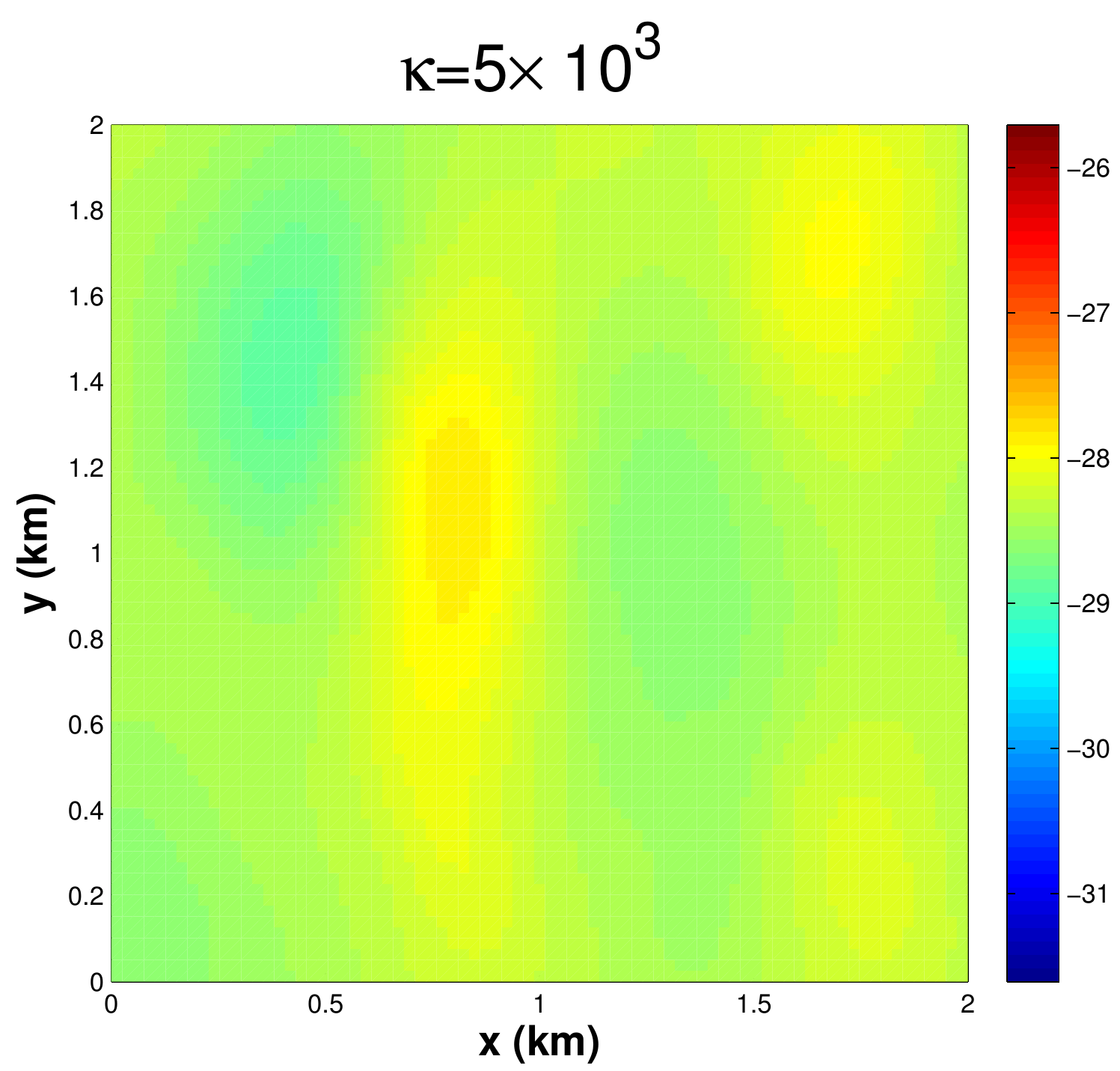}
\caption{Log-permeability estimates obtained with different $\kappa$'s in (\ref{eq:sa}) (i.e. the standard approach for history matching) [$(\log{\textrm{m}^2})$]}   \label{Figure4}
\end{center}
\end{figure}

\begin{figure}
\includegraphics[scale=0.225]{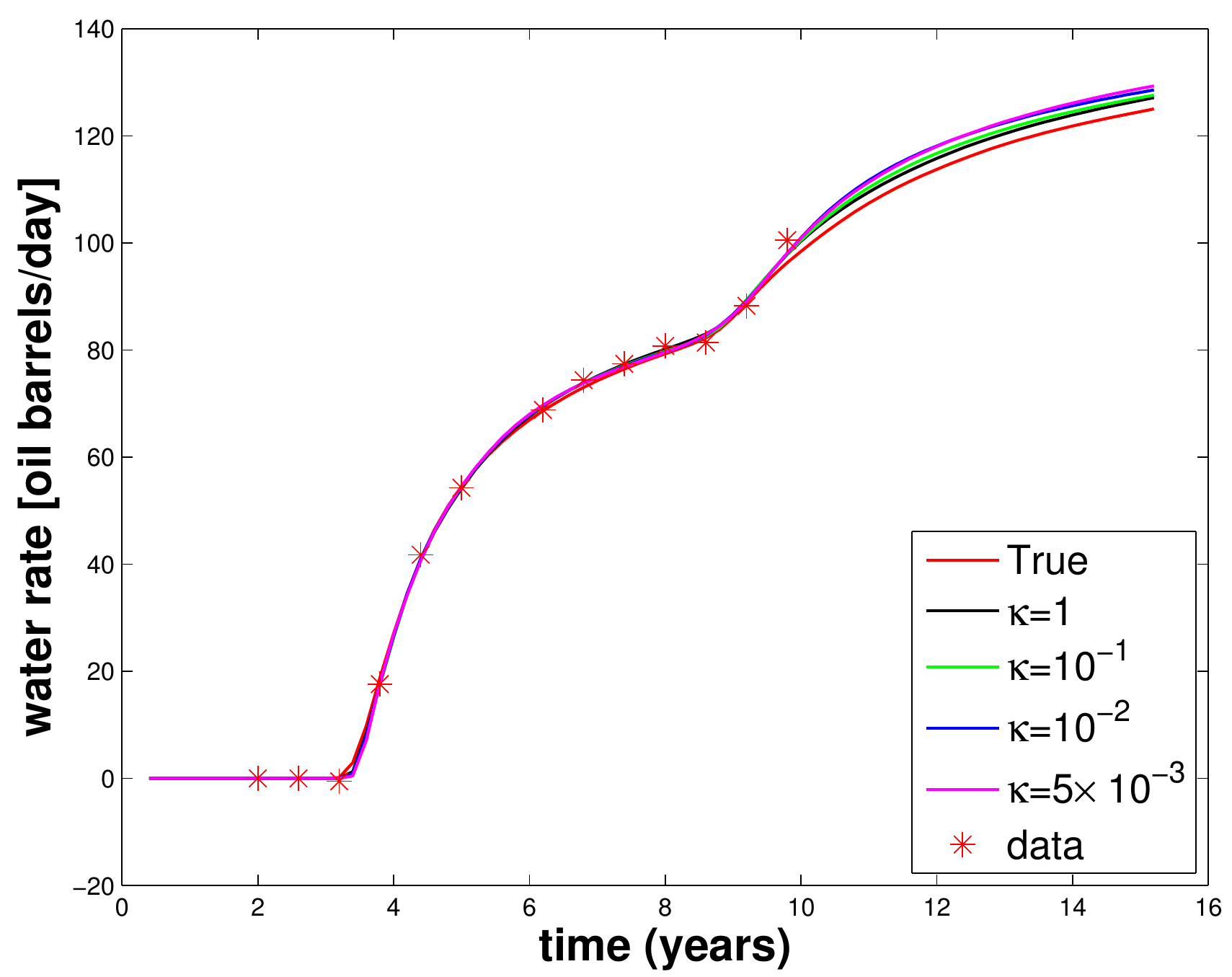}
\includegraphics[scale=0.225]{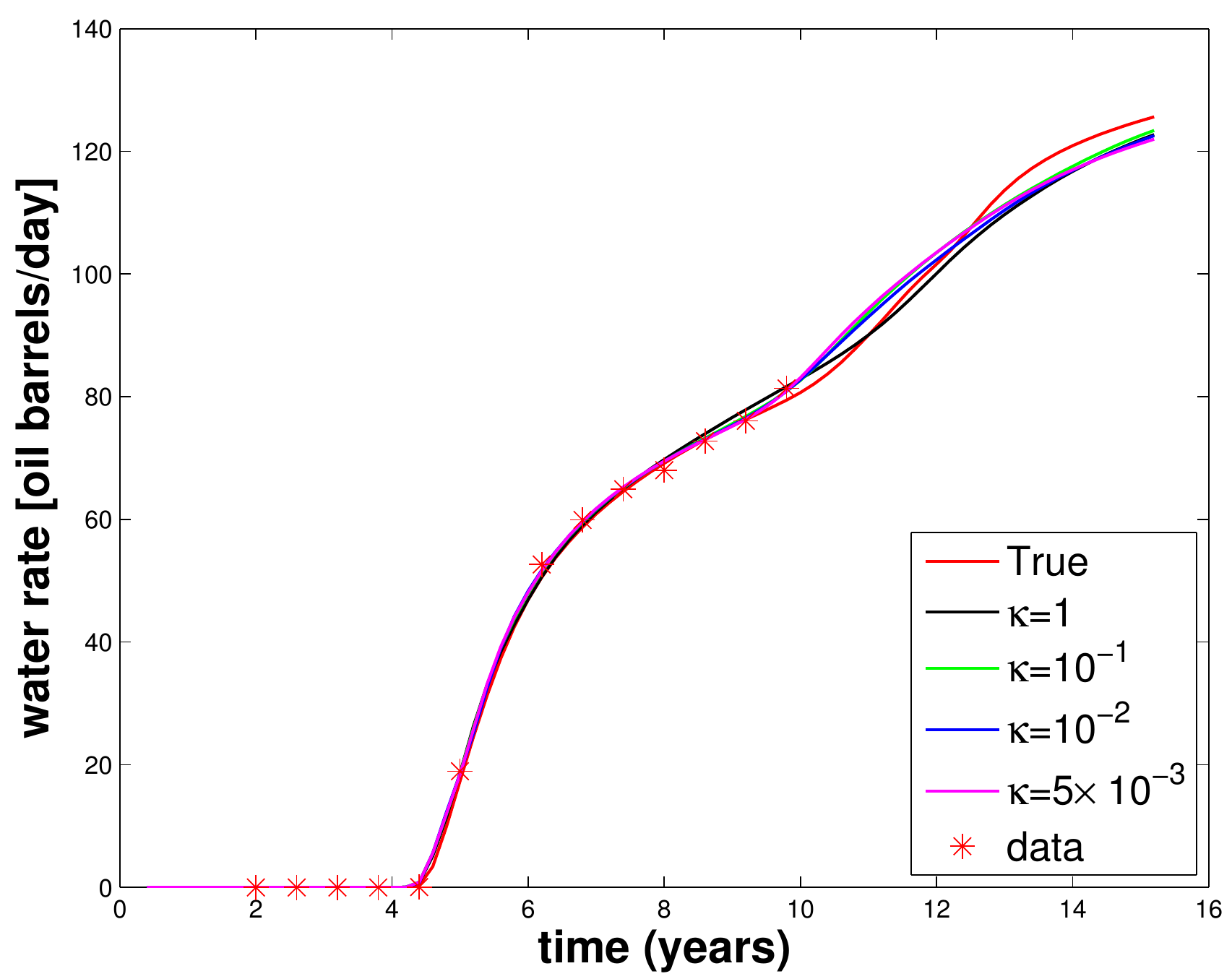}
\includegraphics[scale=0.225]{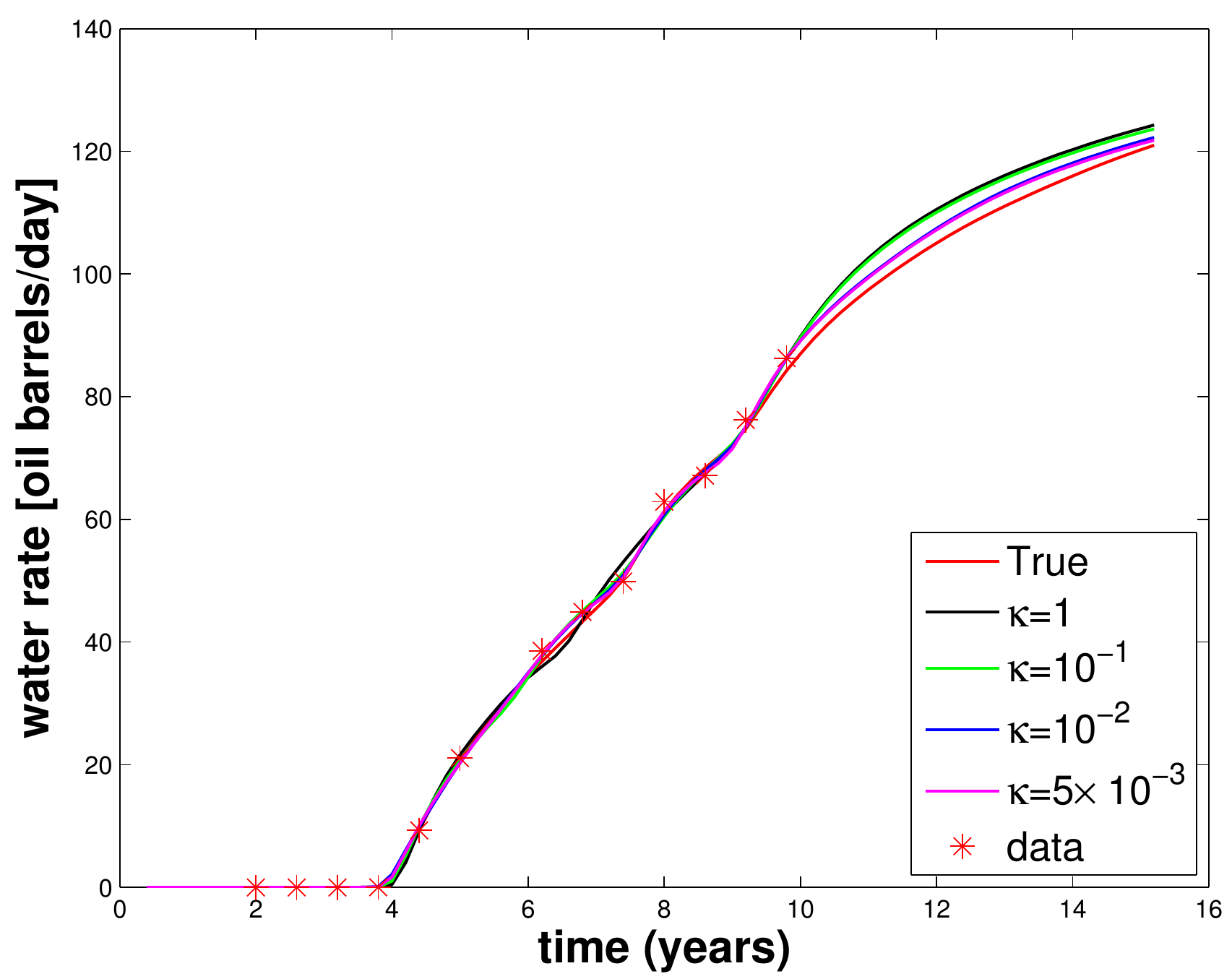}\\
\includegraphics[scale=0.225]{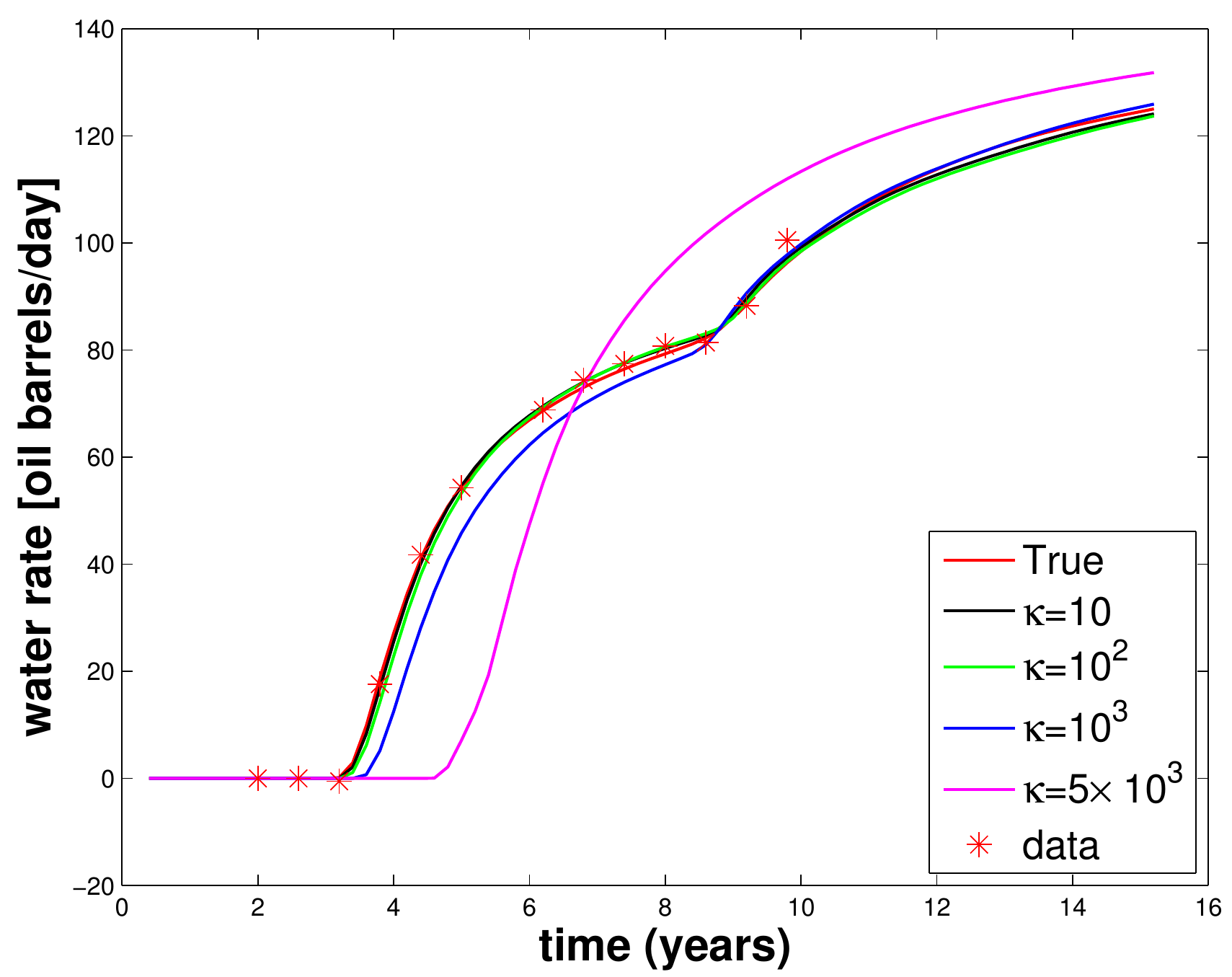}
\includegraphics[scale=0.225]{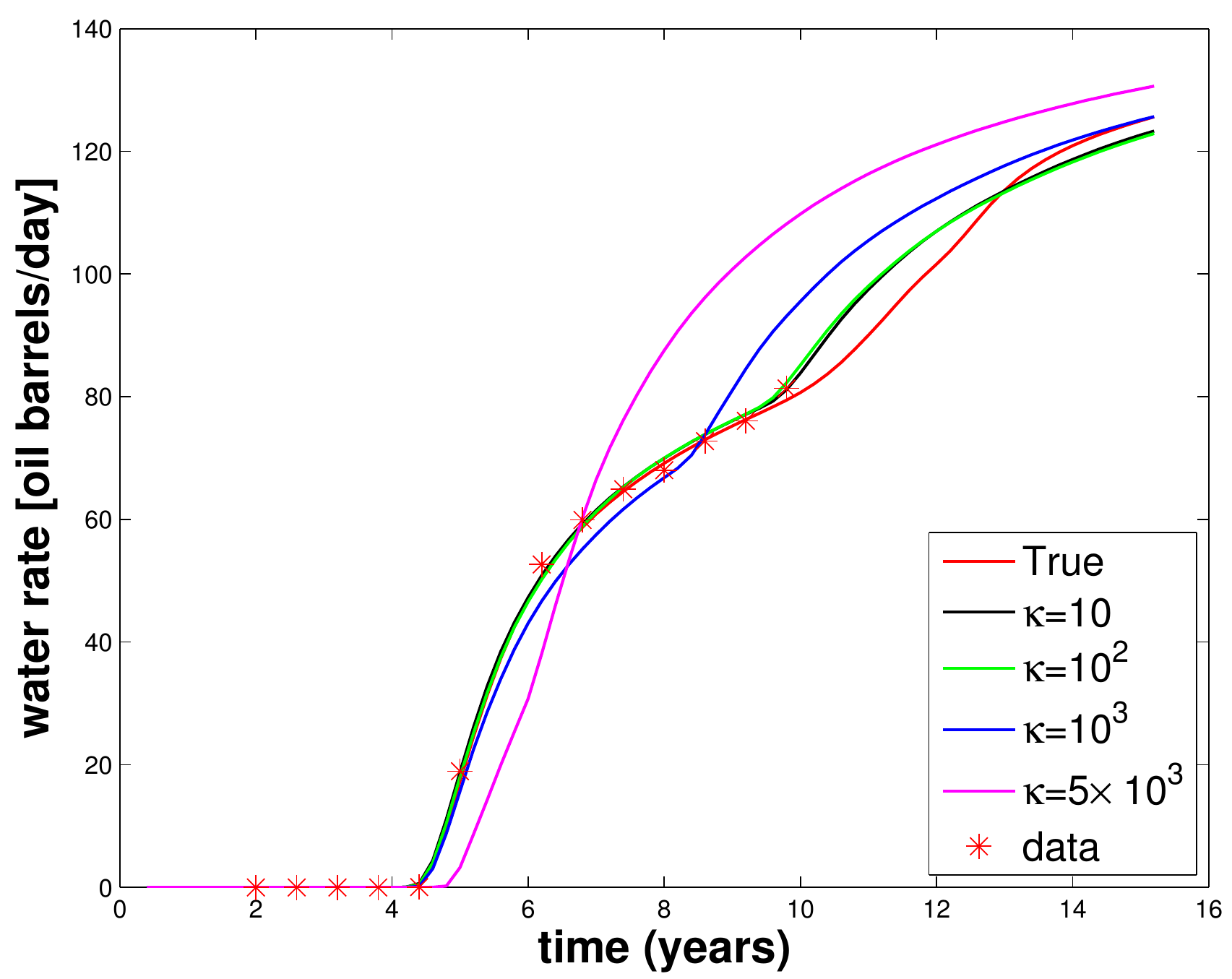}
\includegraphics[scale=0.225]{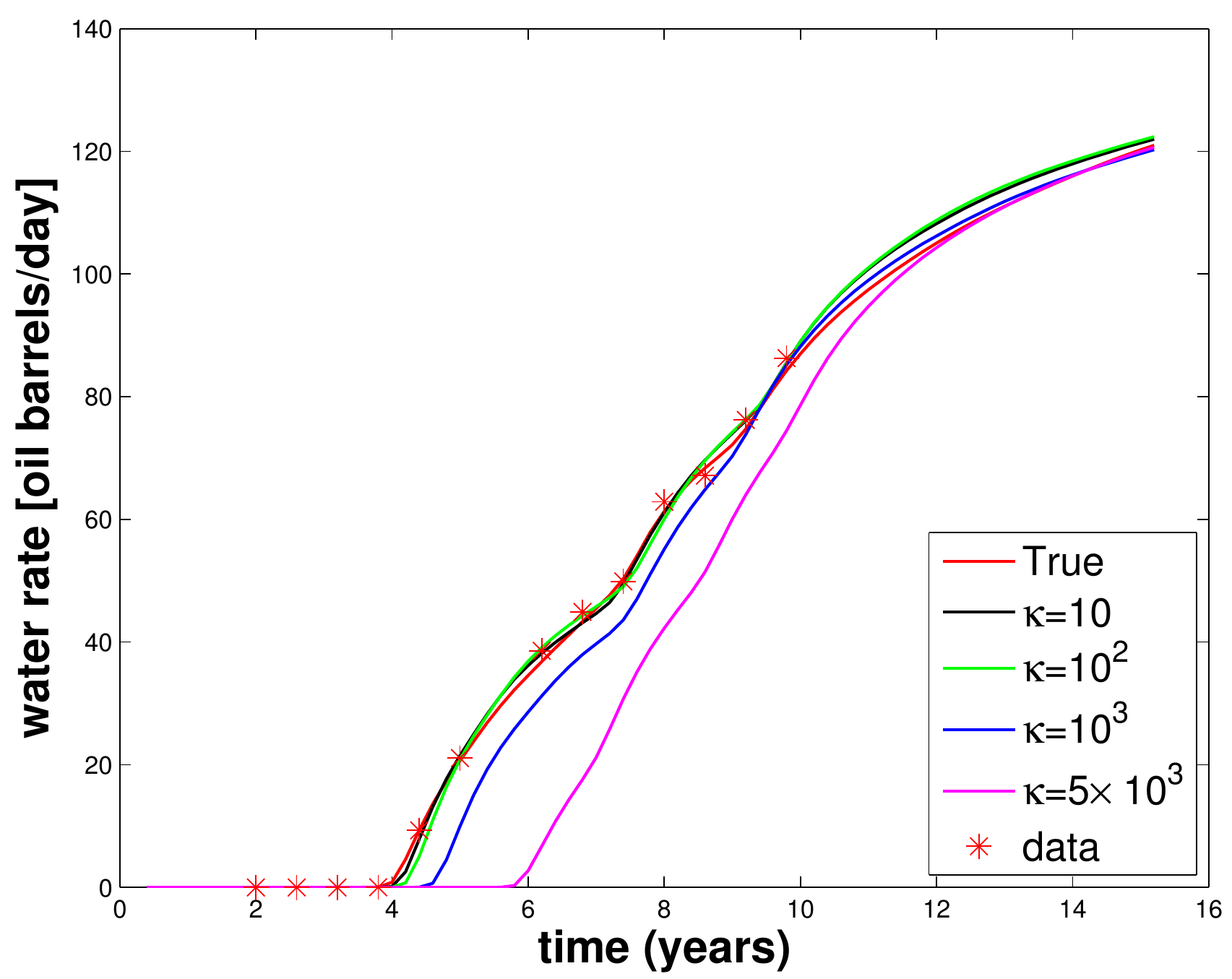}
\caption{Water rates [bbl/day]. From left to right: Wells $P_{2}$, $P_{4}$ and $P_{5}$. Top: Experiments for $\kappa\leq 1$. Bottom: Experiments for $\kappa\ge 10$ } 
\label{Figure4B}
\end{figure}

\begin{figure}
\includegraphics[scale=0.22]{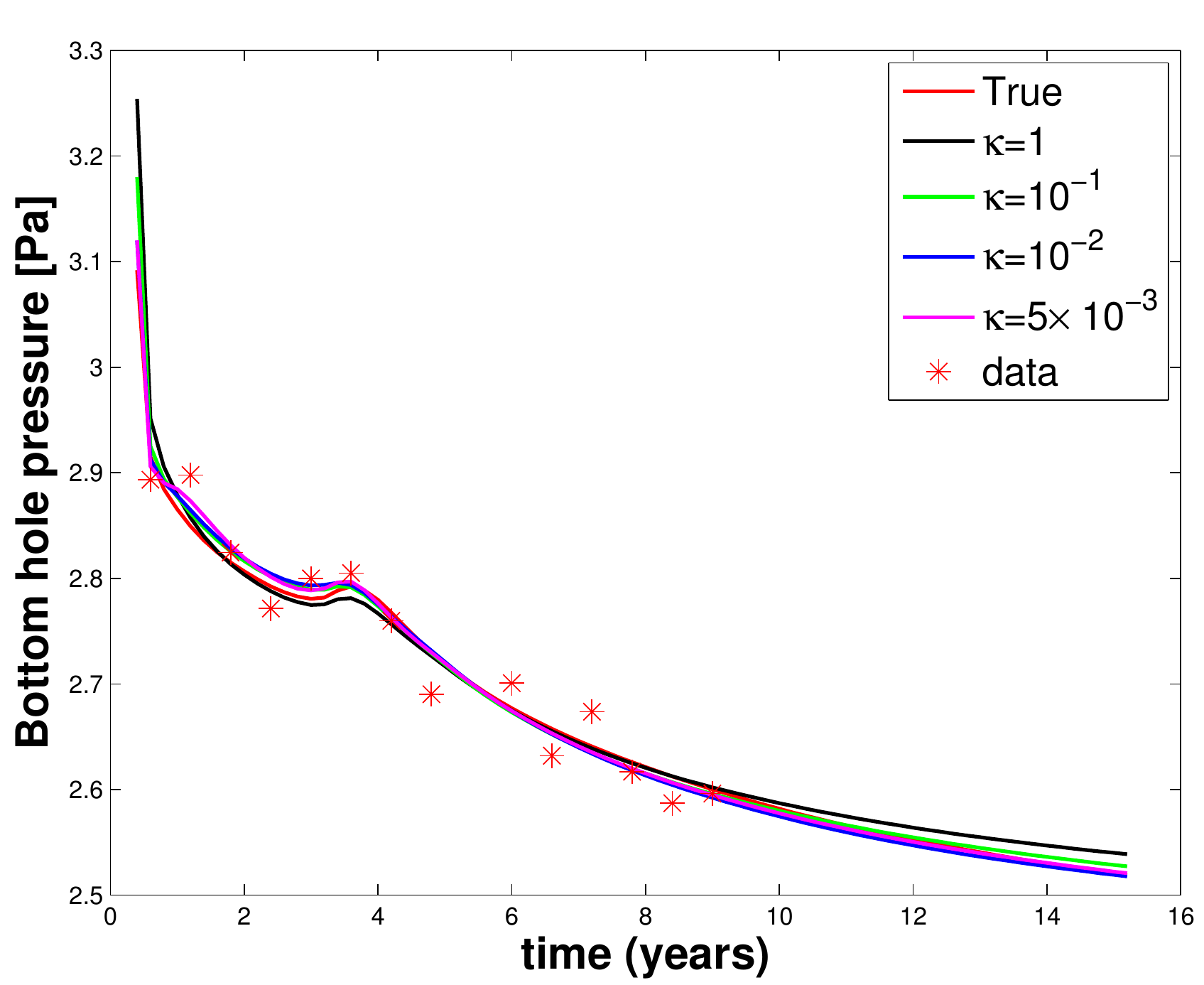}
\includegraphics[scale=0.22]{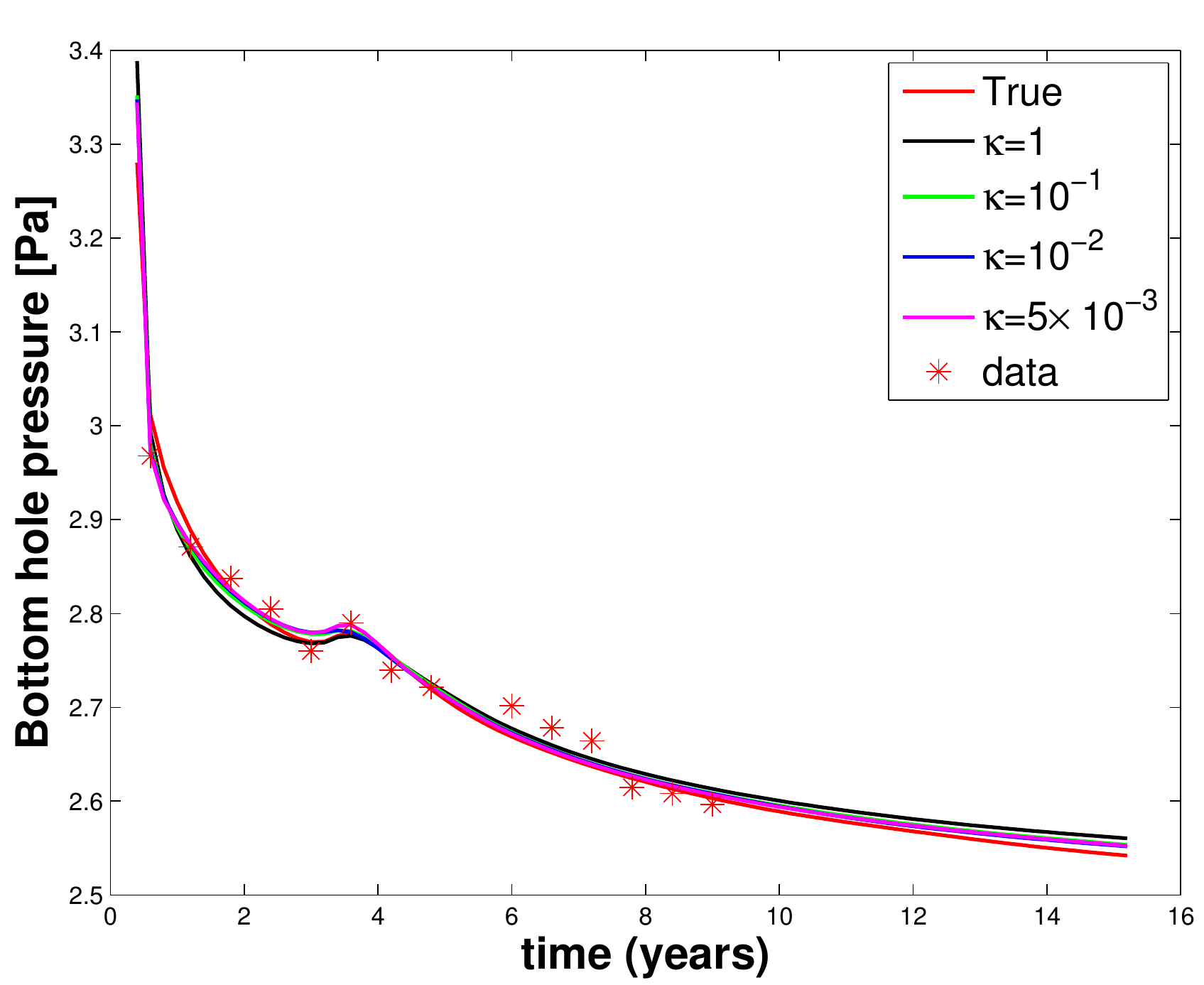}
\includegraphics[scale=0.22]{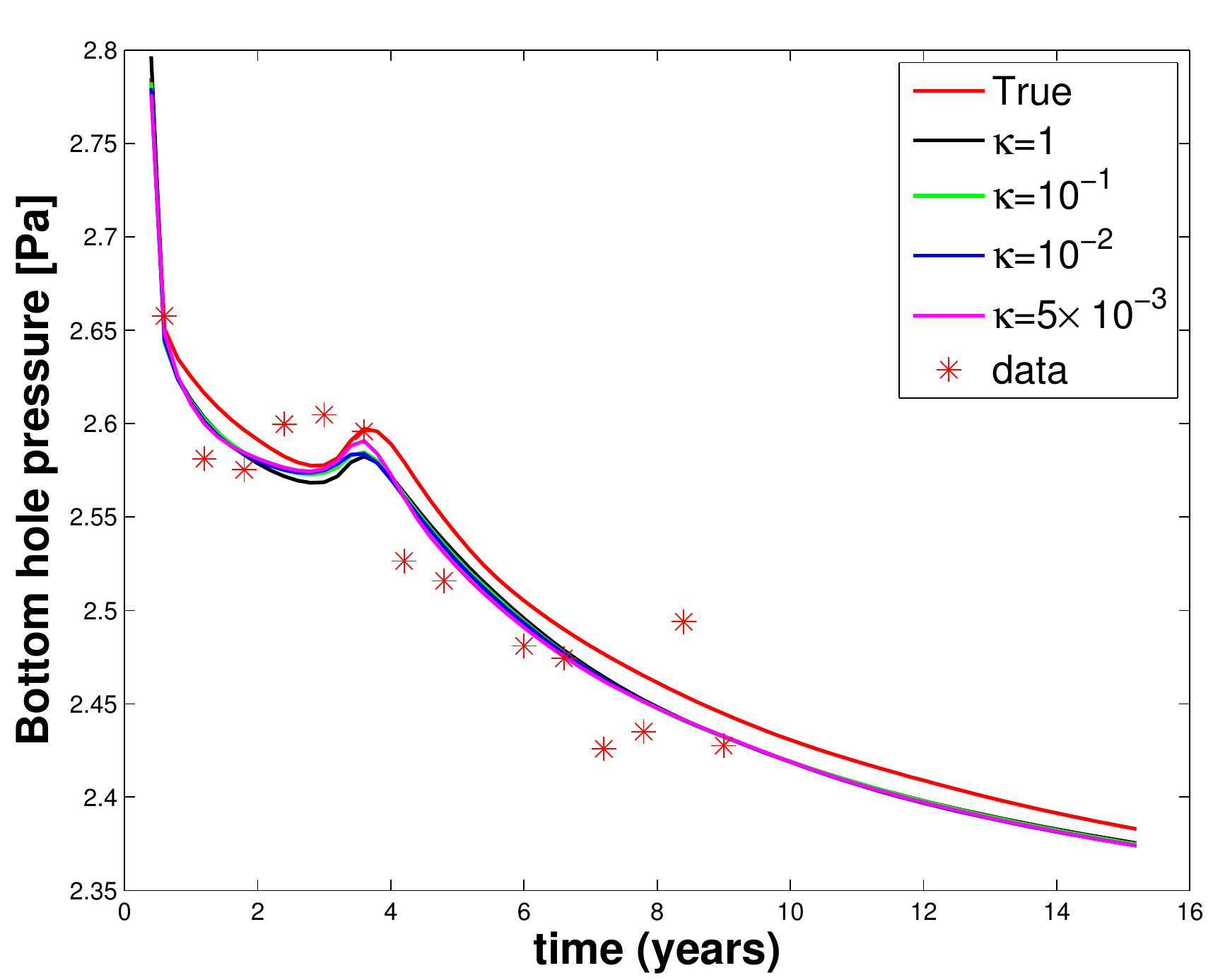}\\
\includegraphics[scale=0.22]{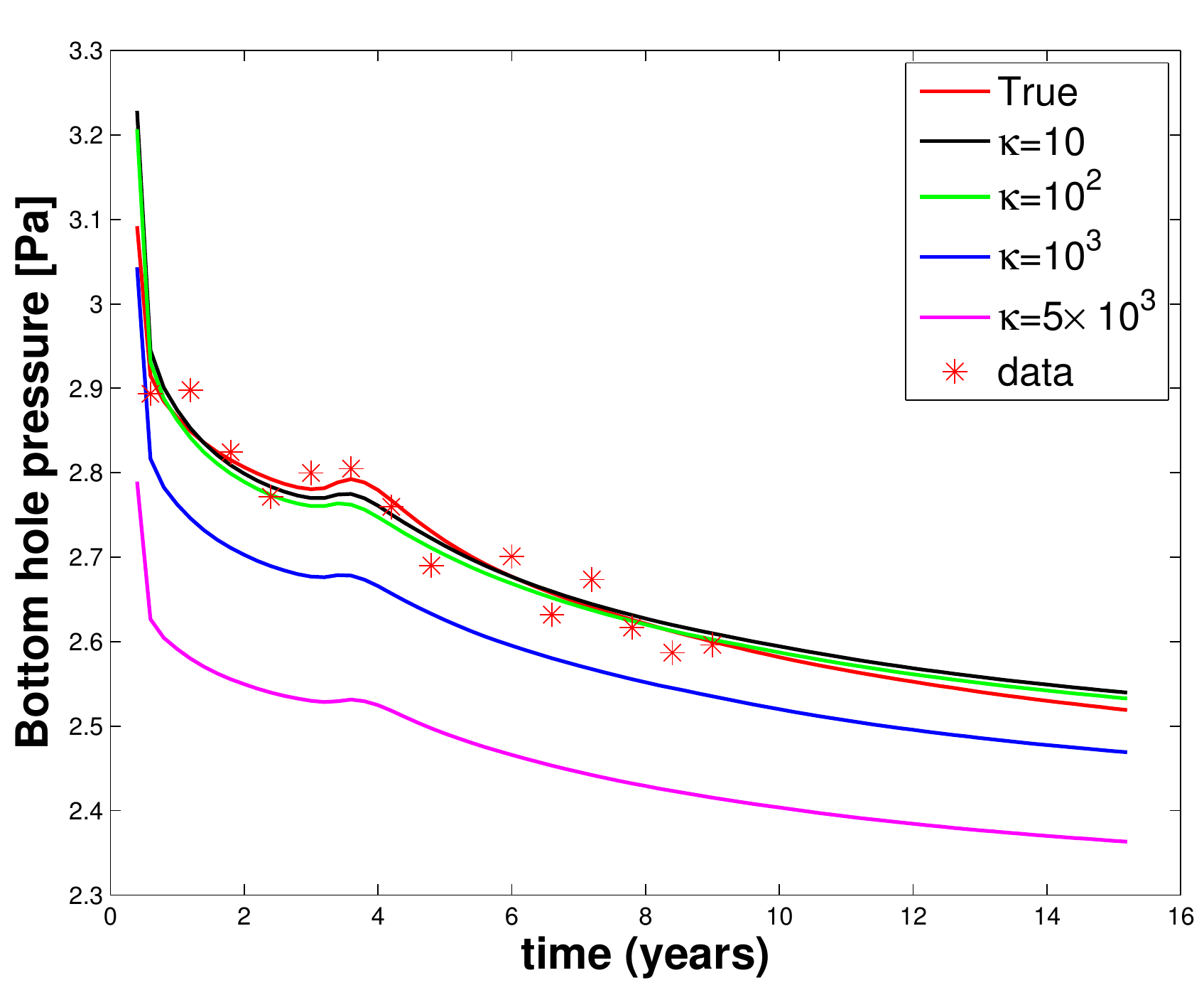}
\includegraphics[scale=0.22]{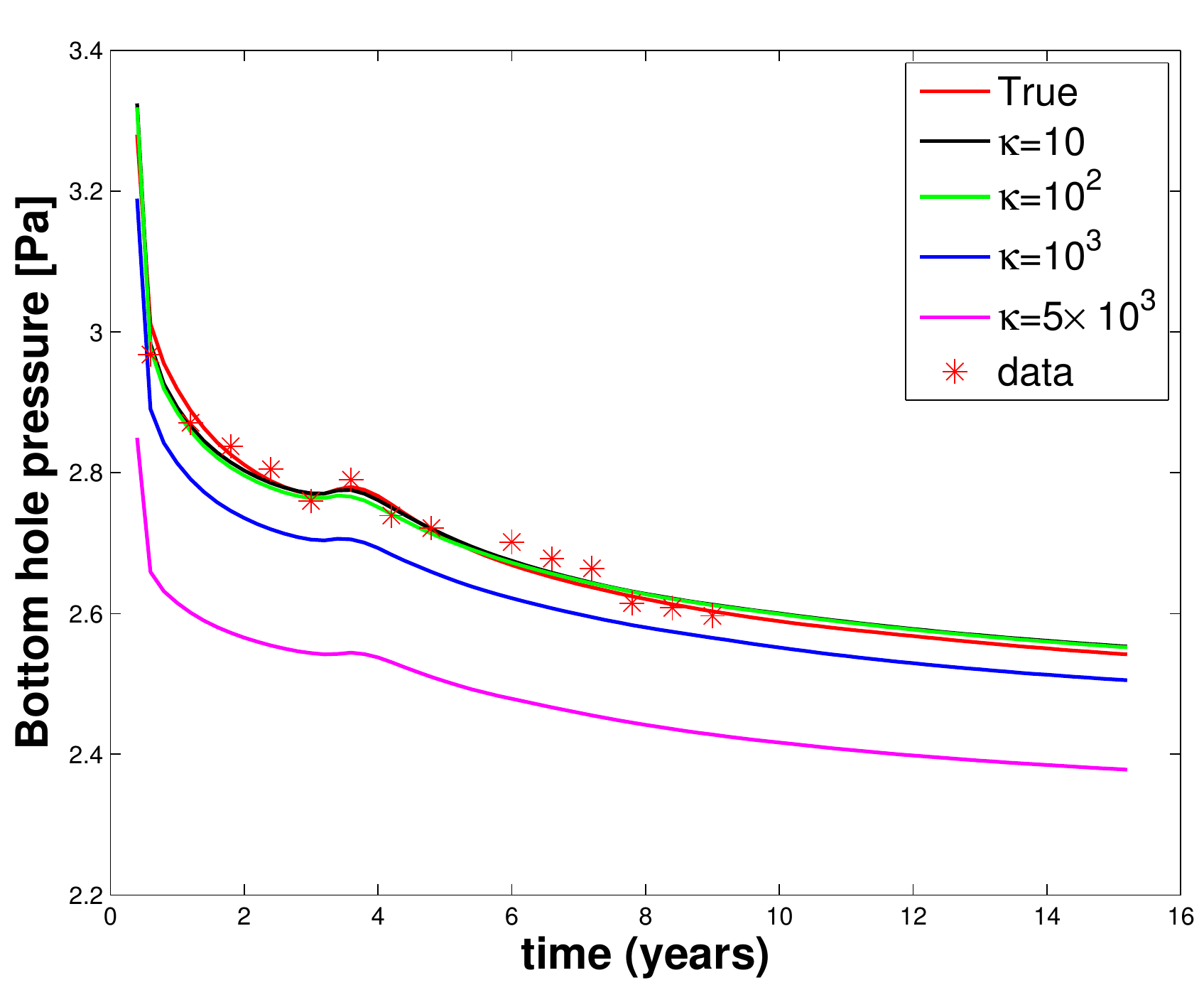}
\includegraphics[scale=0.22]{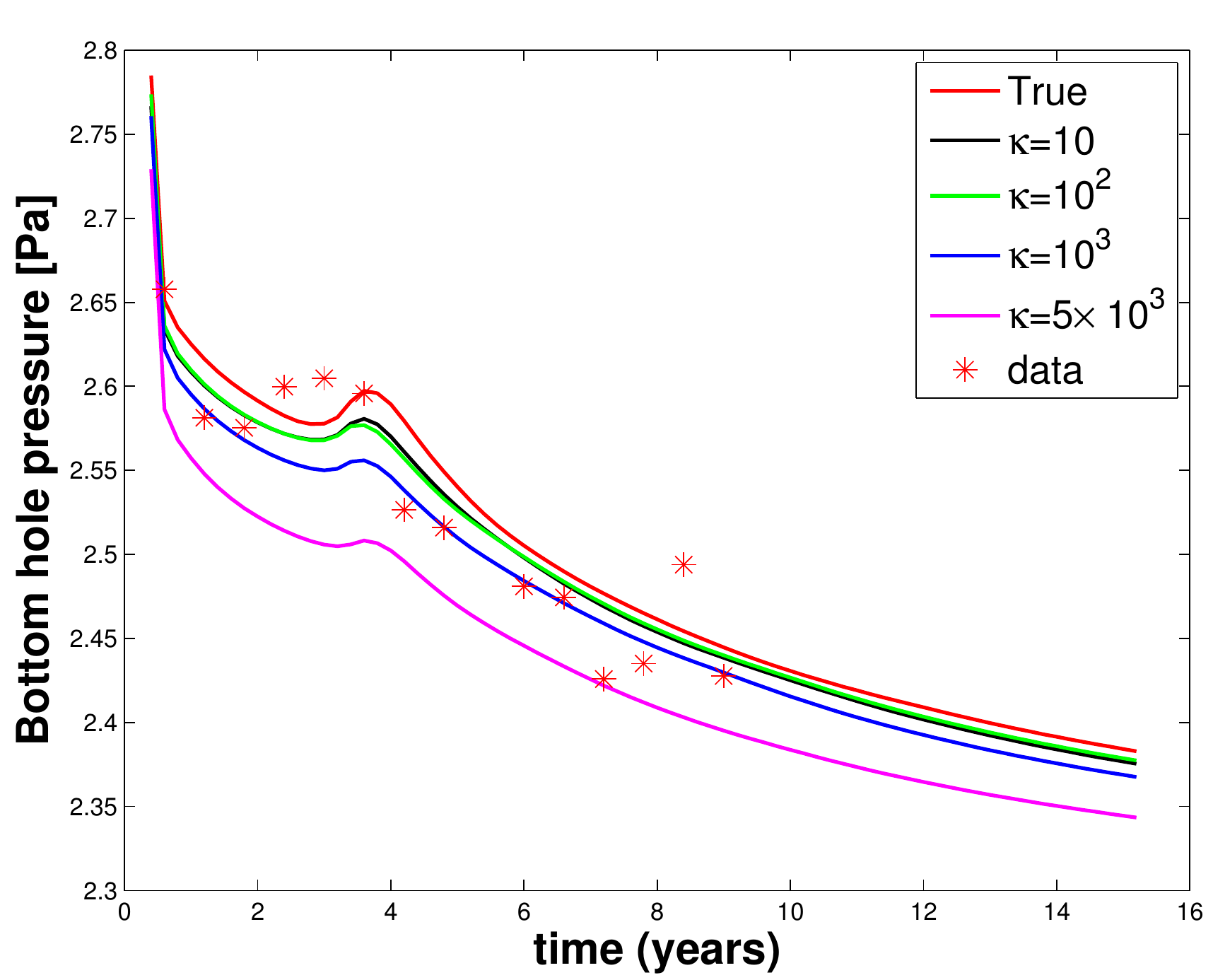}
\caption{Bottom hole pressure [Pa]. From left to right: Wells $I_{2}$, $I_{3}$ and $I_{4}$. Top: Experiments for $\kappa\leq 1$. Bottom: Experiments for $\kappa\ge 10$ } 
\label{Figure4C}
\end{figure}

We now consider the minimization of (\ref{eq:saB}) for $\kappa\ge 10$. The effect of larger regularization in (\ref{eq:saB}) is observed in Figure \ref{Figure3} (right) where, after some number of iterations, the error is indeed stabilized. In Figure \ref{Figure3} (left) we show the associated objective functional.  In this case, the stopping criteria (\ref{eq:3.9C}) and (\ref{eq:3.9D}) are both satisfied. However, for larger $\kappa$'s less accurate estimates are obtained. The estimates of the log-permeability obtained for all $\kappa$'s after 35 iterations of the standard method are displayed in Figure \ref{Figure4}. For small $\kappa$, the lack of stability in the computations is reflected in very large values of the log-permeability fields which is consistent with the results reported in \cite{Li,Oliver}. On the other hand, for large $\kappa$, the lack of fidelity of the corresponding estimates can be visually appreciated for $\kappa> 10^{2}$. From Figure \ref{Figure2} (right) and  Figure \ref{Figure3} (right) we conclude that the correct choice of $C$ (i.e. $\kappa=1$) does not lead to the optimal estimate in terms of the error with respect to the truth. In fact, from all the experiments, $\kappa=10^2$ provides the minimal error with respect to the truth. Similar to the previous set of experiments, in Figure \ref{Figure4B} and Figure \ref{Figure4C} we show model predictions during the 10 years of history matching and the prediction time of 5 years. Note that for small $\kappa$ ($\kappa\leq 1$), all the estimates provide a good data match even though the quality of the corresponding geologic properties (see Figure \ref{Figure4}) is severely degraded for small $\kappa$. In contrast, the lack of fidelity for larger values of $\kappa$ corresponds to poor estimates of the data match as we expected.

In contrast to the standard approach (see Figure \ref{Figure4}), even for small $\kappa$ in (\ref{eq:3.10}) the regularizing LM scheme produce stable estimates of the the true log-permeability (see Figure \ref{Figure7}). In addition, the geological constraint of $C$ is enforced with the regularizing LM scheme. However, by producing stable estimates of the minimization of (\ref{eq:1.1}) the LM scheme avoids the potential lack of fidelity of the standard approach due to the potential overestimation of the prior term. Moreover, as we indicated earlier, the computational cost per iteration of our implementation of the regularizing LM scheme is equivalent to the cost per iteration of the standard method of \cite{Li,Tavakoli,svdRML,Oliver} for minimizing (\ref{eq:sa}). From Figure \ref{Figure3} and \ref{Figure6} we observe that for $\kappa\ge 10^2$, the convergence for both approaches is achieved after 15 iterations. However, the accuracy in terms of the relative error seems to be outperformed by the regularizing LM scheme for larger values of $\kappa$. On the other hand, for the case with $\kappa < 1$, the convergence of the standard approach is not achieved due to the lack of stability reflected in the increase in the relative error. 

\section{Conclusions}\label{Conclu}

While the main contribution of this paper is the implementation of the regularizing LM scheme for history matching, our general aim is to promote further investigations of implementation based on well established theories for approximating stable solutions of history matching problems posed as the minimization of (\ref{eq:1.1}). Although the discussions and experiments of this paper are based on the LM method, the fundamental ideas can be applied to other techniques. In particular, in the iterative regularization literature, there are analogous gradient based (e.g. Landweber iteration, steepest descent) and quasi-Newton methods (BFGS, conjugate gradients) whose aim is to solve nonlinear inverse ill-posed problems such as the one presented here.  Moreover, similarly to the LM scheme, those gradient-based and quasi-Newton methods also share similarities with the corresponding ones for solving (\ref{eq:sa}) in the standard approach. Those similarities may lead to straightforward implementations of iterative regularization techniques when the standard technologies are already available.

Conducting history matching by minimizing (\ref{eq:sa}) has been often motivated from the Bayesian formulation for data assimilation. Under Gaussian assumptions, the minimizer of (\ref{eq:sa}), also called maximum a posteriori (MAP) estimate, maximizes the conditional  \textit{posterior} probability measure of the unknown given the observed data $y^{\eta}$ \cite{Oliver}. While our proposed history matching approach based on iterative regularization techniques is entirely deterministic, there is a potential use of these techniques within the context of Bayesian data assimilation for uncertainty quantification. This conjecture follows from the fact that some standard techniques that approximate the posterior distribution of the Bayesian framework are constructed by randomizing the solution to deterministic problems \cite{Iglesias6}. On the other hand, iterative regularization techniques can also be used in the context of facies identification as suggested in \cite{Iglesias4} where a geometric-based iterative regularization approach was applied for the estimation of geologic facies given data from an oil-water reservoir model similar to the one considered here. Iterative regularization provides then a broad spectrum of techniques that can be potentially used in the estimation of geologic properties in reservoir models.

\section{Appendix: The forward operator}\label{ReservoirModels}

We recall that $G$ is the forward operator that maps the geologic parameters to the production data. We briefly describe the forward operator that we use for the numerical experiments presented in Section \ref{sa} and Section \ref{Numerics}. We consider simplified two-dimensional models typical for testing history matching algorithms. The domain of the reservoir is denoted by $D$; the absolute permeability and porosity are denoted by $K$ and $\phi$ respectively. The interval $[0,T]$ ($T > 0$) is the time interval of interest for the flow simulation. For simplicity we assume that the only unknown parameter is $u=\log{K}$. Nevertheless, all the techniques and implementations that we describe in subsequent sections can be extended to include additional parameters (e.g. porosity).

We consider an  incompressible oil-water reservoir model initially saturated with oil and irreducible water. the water and oil phase are indexed by $\beta = w$ and $\beta = o$, respectively. We are interested in a waterflood process where water is injected at $N_I$ injection wells located at $\{x_I^{l}\}_{l=1}^{N_{I}}$ . Water and oil are produced at $N_P$ production wells located at$\{x_P^{l}\}_{l=1}^{N_{P}}$. Additionally, we assume that injection wells are operated under prescribed rates $\{q_{I}^l(t)\}_{l=1}^{N_I}$ while production wells are constrained to the total flow rate $\{q_{P}^{l}(t)\}_{l=1}^{N_P}$. The pressure $p(x,t)$ and the saturation $s(x,t)$ ($(x, t ) \in D \times [0, T ]$) are the state variables. From standard arguments it can be shown that $(s, p)$ is the solution to the following system \cite{Chen}
\begin{eqnarray}\label{eq:2.7}
-\nabla \cdot \lambda(s) e^{u}\nabla p= \sum_{l=1}^{N_{I}}q_{I}^l \delta(x-x_I^l)+\sum_{l=1}^{N_{w}}q_{P}^{l}\delta(x-x_P^l)\\
\phi \frac{\partial s}{\partial t} -\nabla \cdot \lambda_w(s) e^{u}\nabla p= \sum_{l=1}^{N_{I}}q_{I}^l \delta (x-x_I^l)+\sum_{l=1}^{N_{w}}\frac{\lambda_{w}}{\lambda}q_{P}^{l}\delta(x-x_P^l)\label{eq:2.7B}
\end{eqnarray}
in $D\times(0,T]$, where $\delta(x-x_P^l)$ and $\delta(x-x_I^l)$ are the (possibly mollified) Dirac deltas.  In (\ref{eq:2.7})-(\ref{eq:2.7B}), $\lambda_w(s)$ and $\lambda(s)$ denote the water and total mobility defined by
\begin{eqnarray}\label{eq:2.8}
\lambda_w(s)=\frac{k_{rw}(s)}{\mu_{w}},\qquad \lambda(s)=\frac{k_{ro}(s)}{\mu_{o}}+\lambda_w(s)
\end{eqnarray}
where $k_{r\gamma}(s)$ and $\mu_{\gamma}$ denote the relative permeability and the viscosity of the $\gamma$-phase fluid, respectively. Furthermore, we assume that
\begin{eqnarray}
k_{rw}(s)=a_{w}\Bigg[\frac{s-s_{iw}}{1-s_{iw}-s_{or}}\Bigg]^2,\qquad 
k_{ro}(s)=a_{o}\Bigg[\frac{1-s-s_{or}}{1-s_{iw}-s_{or}}\Bigg]^2
\end{eqnarray}
where $a_{w},a_{o}\in (0,1]$, $s_{iw}$ is the irreducible water saturation and $s_{or}$ is the residual oil saturation. We
additionally prescribe initial conditions for pressure and water saturation
\begin{eqnarray}\label{eq:2.9}
p=p_0, \qquad s=s_{0}\qquad \textrm{in }   D\times\{0\} 
\end{eqnarray}
For simplicity, no-flow boundary conditions are prescribed on the reservoir boundary
\begin{eqnarray}
- e^{u}\lambda(s)\nabla p\cdot \mathbf{n}&=&0 ~~~~~~~~~~\textrm{on }  \partial D\times (0,T]\\
- e^{u}\lambda_w(s)\nabla p\cdot \mathbf{n}&=&0 ~~~~~~~~~~\textrm{on }  \partial D\times (0,T]\label{eq:2.10}
\end{eqnarray}
Let us assume that there are $N_{M}$ measurement times denoted as before $\{t_{n}\}_{n=1}^{N_{M}}$. We assume measurements of bottom-hole pressure are collected at the injection wells at $\{t_{n}\}_{n=1}^{N_{M}}$. This, according to Peacemen well-model \cite{Chen} is defined by 
\begin{eqnarray}\label{eq:2.11}
M_{n}^{l,I}(p,s)=\frac{q_{I}^l(t_{n})}{\omega^{l}\lambda(s(x_{I}^{l},t_{n}))}+p(x_{I}^l,t_{n})
\end{eqnarray}
for $l=1,\dots,N_{I}$ and $n=1,\dots, N_M$. Analogously, we consider measurements of water and oil rates at the production wells
\begin{eqnarray}\label{eq:2.12}
M_{n}^{l,P_{w}}(p,s)=\frac{\lambda_{w}(s(x_{P}^{l},t))}{\lambda(s(x_{P}^{l},t_{n}))}q_{P}^{l}(t_{n}),\qquad 
M_{n}^{l,P_{o}}(p,s)=\frac{\lambda_{o}(s(x_{P}^{l},t))}{\lambda(s(x_{P}^{l},t_{n}))}q_{P}^{l}(t_{n})
\end{eqnarray}
for $l=1,\dots, N_{P}$ and $n=1,\dots, N_M$. In (\ref{eq:2.12}), $\lambda_{o}=\lambda-\lambda_{w}$. Let us define the $2N_{P}+N_{I}$-dimensional vector 
\begin{eqnarray}\label{eq:2.13}
M_{n}(p,s)=(M_{n}^{1,I}(p,s), \dots, M_{n}^{N_{I},I}(p,s),M_{n}^{1,P_{w}}(p,s), \dots, M_{n}^{N_{P},P_{w}}(p,s),M_{n}^{1,P_{o}}(p,s), \dots, M_{n}^{N_{P},P_{o}}(p,s))\nonumber\\
\end{eqnarray}
that contains the number of measurements from wells at a given time. The total number of measurements is $N=[2N_{P}+N_{I}]N_{M}$ and the forward map $G:X\to \mathbb{R}^{N}$ is then given by 
expression
\begin{eqnarray}\label{eq:2.14}
G(u)=(M_1(p,s), \dots, M_{N_{M}}(p,s))
\end{eqnarray}
which comprises the production data obtained from production and injection wells at the measurement times.

\begin{acknowledgements}
The authors would like to thank Andrew Stuart for helpful discussions and his generous feedback on the content and structure of the manuscript. The second author acknowledges the support of the Department of Energy (DOE grant number is DE-SC0009286). 

\end{acknowledgements}

\bibliographystyle{plain}
\bibliography{ITER_bib}   

\end{document}